\input amstex
\documentstyle{amsppt}
\NoBlackBoxes

\TagsOnRight

\def\cal{\Cal}
\def\AA{{\cal A}}

\def\HH{{\cal H}}
\def\MM{{\cal M}}
\def\NN{{\cal N}}
\def\JJ{{\cal J}}
\def\UU{{\cal U}}

\def\LL{{\cal L}}

\def\UU{{\cal U}}
\def\VV{{\cal V}}
\def\Z{{\Bbb Z}}
\def\C{{\Bbb C}}
\def\R{{\Bbb R}}

\def\Q{{\Bbb Q}}
\def\e{{\epsilon}}

\def\n{\noindent}
\def\part{{\partial}}
\def\dudtau{{\part u\over \part \tau}}
\def\dvdtau{{\part v\over \part \tau}}
\def\dudt{{\part u\over \part t}}
\def\dvdt{{\part v\over \part t}}

\rightheadtext{Spectral invariants} \leftheadtext{
Yong-Geun Oh }

\topmatter
\title
Spectral invariants, analysis of the Floer moduli space and
geometry of the Hamiltonian diffeomorphism group
\endtitle
\author
Yong-Geun Oh\footnote{Partially supported by the NSF Grant \#
DMS-9971446 \& DMS-0203593, by \# DMS-9729992 in the Institute for
Advanced Study, by the Vilas Associate Award of the University of
Wisconsin and by a grant of the Korean Young Scientist Prize
\hskip6.5cm\hfill}
\endauthor
\address
Department of Mathematics, University of Wisconsin, Madison, WI
53706, ~USA \& Korea Institute for Advanced Study, Seoul, Korea;
oh\@math.wisc.edu
\endaddress

\abstract  In this paper, we apply spectral invariants,
constructed in [Oh5,8], to the study of Hamiltonian
diffeomorphisms of  closed symplectic manifolds $(M,\omega)$.
Using spectral invariants, we first construct an invariant
norm called the {\it spectral norm} on the Hamiltonian diffeomorphism
group and obtain several lower bounds for the spectral norm in
terms of the $\e$-regularity theorem and the symplectic area of
certain pseudo-holomorphic curves. We then apply spectral
invariants to the study of length minimizing properties of certain
Hamiltonian paths {\it among all paths}. In addition to the
construction of spectral invariants, these applications rely
heavily on the {\it chain level Floer theory} and on some
existence theorems with energy bounds of pseudo-holomorphic
sections of certain Hamiltonian fibrations with prescribed
monodromy. The existence scheme that we develop in this paper in
turn relies on some careful geometric analysis
involving {\it adiabatic degeneration} and {\it thick-thin
decomposition} of the Floer moduli spaces
which has an independent interest of its own.

We assume that $(M,\omega)$ is {\it strongly semi-positive}
throughout, which will be removed in a sequel.
\endabstract


\endtopmatter
\document

\quad MSC-class: 53D35, 53D40

\bigskip

\centerline{\bf Contents}
\medskip

\n \S1. Introduction and statement of the main results
\smallskip

\n \S2. Preliminaries
\smallskip
\quad 2.1. The action functional and the action spectrum
\smallskip
\quad 2.2. Filtered Floer homology
\smallskip
\quad 2.3. Spectral invariants and Floer fundamental cycles
\smallskip
\quad 2.4. Comparison of two Cauchy-Riemann equations
\bigskip

\n{\it PART I: APPLICATIONS OF SPECTRAL INVARIANTS}
\medskip
\n \S3. Definition of the spectral norm
\smallskip

\n \S4. The $\e$-regularity theorem and its consequences
\smallskip

\n \S5. Proof of nondegeneracy
\smallskip

\n \S6. The spectral norm, the ball area and the homological area
of $\phi$
\smallskip

\n \S7. The homological area and the generating function of $\phi$
\bigskip

\n{\it PART II: ADIABATIC DEGENERATION AND  THICK-THIN}
\smallskip
\hskip0.25in {\it DECOMPOSITION}
\medskip

\n \S8. Proof of the fundamental existence theorem
\smallskip
\quad 8.1. Pants products and Hamiltonian fibrations
\smallskip
\quad 8.2. Heuristic discussion of the proof of Theorem 5.4
\smallskip
\quad 8.3. Construction of solutions I: analysis of thick parts
\smallskip
\quad 8.4. Construction of solutions II: analysis of thin parts
\smallskip
\quad 8.5. Construction of solutions III: wrap-up
\smallskip

\n \S9. Local Floer complex and its homology
\smallskip
\quad 9.1. The thin part of the Floer boundary operator
\smallskip
\quad 9.2. The thin part of the pants product
\smallskip
\quad 9.3. The canonical Floer fundamental cycle
\smallskip
\quad 9.4. Proof of Theorem 7.1
\smallskip

\n \S10. Thick-thin decomposition of the Floer moduli space
\smallskip

\n Appendix : Construction of flat connections

\bigskip

\head
{\bf \S 1. Introduction and the statement of the main results}
\endhead

This is a sequel to the papers [Oh5,8] in which we constructed
some invariants called {\it spectral invariants} of
(time dependent) Hamiltonian functions or of associated
Hamiltonian paths. In the present paper we provide
several applications of these spectral invariants to the study of
the geometry of Hamiltonian diffeomorphisms.

We start by reviewing the definition of the Hofer norm $\|\phi\|$ of
a Hamiltonian diffeomorphism $\phi$. For a given Hamiltonian
diffeomorphism $\phi \in {\cal H}am(M,\omega)$, we write $F
\mapsto \phi$ if the (time-dependent) Hamiltonian function
$F:[0,1] \times M \to \R$ satisfies $\phi^1_F = \phi$. We will always
consider a normalized Hamiltonian $F$, that is satisfying
$\int_M F_t d\mu = 0$. The composed flow
$t \mapsto \phi_F^t \circ \phi_H^t$ is generated by
the product Hamiltonian $F\# H$
$$
F\# H(t, x) = F(t, x) + H(t, (\phi_F^t)^{-1}(x)) \tag 1.1
$$
and the inverse flow $t \mapsto (\phi_H^t)^{-1}$ is generated by
the inverse Hamiltonian $\overline H$
$$
\overline H(t,x) = - H(t, \phi_H^t(x)). \tag 1.2
$$
We recall the following standard definitions :
$$
\aligned
E^-(H) & =  \int_0^1 - \min H_t \, dt, \quad
E^+(H) = \int_0^1 \max H_t \, dt \\
\|H\| & =  E^-(H) + E^+(H).
\endaligned
\tag 1.3
$$
The Hofer norm $\|\phi\|$ is defined by
$$
\|\phi\| = \inf_{H\mapsto \phi} \|H\|.
$$
Noting that $E^+(H) = E^-(\overline H)$, $\|H\|$ and $\|\phi\|$ can
be rewritten as
$$
\align
\|H\| & =  E^-(H) + E^-(\overline H) \tag 1.4 \\
\|\phi\| & =  \inf_{H \mapsto \phi}(E^-(H) + E^-(\overline H)).
\tag 1.5
\endalign
$$
Next we recall the basic properties of the particular
spectral invariant $\rho(H;1)$ corresponding to the
quantum cohomology class $1 \in QH^*(M)$. First, $\rho(H;1)$ satisfies the
general inequality
$$
-E^+(H) \leq \rho(H;1) \leq E^-(H) \tag 1.6
$$
(see Theorem 3.1), and a fortiori the weaker inequality
$$
-\|H\| \leq \gamma(H) \leq \|H\| \tag 1.7
$$
for the function $\gamma: C^\infty_m([0,1] \times M) \to \R$
defined by
$$
\gamma(H) = \rho(H;1) + \rho(\overline H;1). \tag 1.8
$$
A rather interesting stronger lower bound of $\gamma$
comes from the general triangle inequality
$$
\rho(H;a) + \rho(F;b) \geq \rho(H\# F;a\cdot b) \tag 1.9
$$
for spectral invariants where $a \cdot b$ is the quantum
product of $a, \, b \in QH^*(M)$ : Restricting to $a=b=1$, and
using the identity $H \# \overline H = \underline 0$
and the normalization axiom $\rho(\underline 0;1) = 0$, we derive
$$
\gamma(H) = \rho(H;1) + \rho(\overline H;1) \geq \rho(\underline
0;1) = 0.
$$
Here $\underline 0$ denotes the zero function.

This leads us to the definition of a non-negative function
$\gamma: Ham(M,\omega) \to \R_+$ defined by
$$
\gamma(\phi):= \inf_{H \mapsto \phi} \Big(\rho(H;1) +
\rho(\overline H;1)\Big). \tag 1.10
$$
We would like to point out the similarity between the expressions
(1.5) and (1.10). However it should be noted that
the nonnegativity of (1.5) is automatic while
the nonnegativity of (1.10) is a nontrivial consequence of the triangle
inequality (1.9), which is closely tied to the {\it positivity}
phenomenon in symplectic topology manifested through pseudoholomorphic
curves.

{\it Throughout the paper, we will assume the spectrality axiom
$$
\rho(H;1) \in \text{Spec}(H) \tag 1.11
$$
for nondegenerate Hamiltonians $H$.} This
then implies the homotopy invariance of $H \mapsto \rho(H;1)$ which
is in turn needed to prove Lemma 5.2 (5.5).
The spectrality axiom of $\rho(H;a)$ for arbitrary $H$ and $a$
was proved in the rational $(M,\omega)$ in [Oh8],
and in the irrational $(M,\omega)$ in the
preprint (math.SG/0406449) (for nondegenerate Hamiltonians).
Therefore this assumption is likely to be superfluous.
However, {\it until the proof of the spectrality axiom
for the irrational case in the preprint (math.SG/0406449) is confirmed,
we will assume that $(M,\omega)$ is rational in the present paper}.

The following theorem provides the general properties of
$\gamma: Ham(M,\omega) \to \R_+$.

\proclaim{Theorem I} $\gamma : Ham(M,\omega) \to \R_+$
satisfies the following properties :
\roster \item $\phi= id$ if and
only if $\gamma(\phi) = 0$ \item $\gamma(\eta\phi \eta^{-1}) =
\gamma(\phi)$ for any symplectic diffeomorphism $\eta$ \item
$\gamma(\psi\phi) \leq \gamma(\psi) + \gamma(\phi)$ \item
$\gamma(\phi^{-1}) = \gamma(\phi)$ \item $\gamma(\phi) \leq
\|\phi\|$.
\endroster
In particular, $\gamma$ defines a symmetric
(i.e., $\gamma(\phi) = \gamma(\phi^{-1}$))
invariant norm on $Ham(M,\omega)$.
\endproclaim
All but (1), nondegeneracy of $\gamma$, are immediate consequences of
the general properties of $\rho(H;1)$.

A similar construction was previously carried out by Viterbo
[V] and by the author [Oh4] for the cotangent bundle
considering Hamiltonian deformations of the zero section,
and by Schwarz [Sc] for Hamiltonian
diffeomorphisms on the (symplectically) aspherical symplectic
manifold $(M,\omega)$. In all these cases, $\gamma$ was defined by
$$
\gamma(H) = \rho(H;1) - \rho(H;\mu) \tag 1.12
$$
$\mu$ being the volume class in $H^*(M)$. The
quantity (1.12) coincides with the expression (1.8) in the
aspherical $(M,\omega)$ because in that case we have the
additional identity
$$
\rho(\overline H;1) = - \rho(H;\mu)
$$
[V], [Sc]. But Polterovich observed [Po2] that this latter
identity fails in the non-exact case due to the quantum
contributions. It turns out that for the general case
the right way of defining an invariant norm satisfying
the triangle inequality is (1.8) and (1.10), not (1.12).
As illustrated by Ostrover [Os] in the
aspherical case, $\gamma$ is different from Hofer's norm.

In a way, our nondegeneracy proof for $\gamma$ is a combination
of the techniques Hofer [Ho] and Lalonde-McDuff [LM] used to prove
nondegeneracy of Hofer's norm. Recall Hofer used the classical critical
point theory argument in his proof for $\C^n$ while Lalonde-McDuff
used the method of pseudoholomorphic curves via general Nonsqueezing
Theorem. What we present is a more sophisticated critical point theory for
general $(M,\omega)$ via Floer theory, which is also naturally tied to
pseudoholomorphic curves and quantum cohomology through deformations
of Hamiltonian functions. However for the lower estimate of
$\gamma(\phi)$, instead of the {\it ball area} of $\phi$ (see Definition 4.1),
we use other geometric invariants of $\phi$ which naturally arise
in our chain level arguments and better suits
Floer theory. These invariants, denoted by
$A(\phi)$ and $A(\phi;1)$ respectively, depend only on $\phi$ itself,
independently of its generating Hamiltonian. Our definition of
$A(\phi)$ resembles that of a similar invariant used by the author
in [Oh2,3] and Chekanov [Ch] in the context of Lagrangian submanifolds.

Now we will provide a precise description of the analytic invariant
$A(\phi)$ here in the introduction partly because it is needed
to explain other main results of the present paper. It also
has a second purpose : {\it we would like to
emphasize the precise details needed in various constructions
of the Floer theory and spectral invariants} which are essential
for the applications studied in the present paper. These
details have not been emphasized enough in the literature but
we believe that they deserve more attention and scrutiny
in the future.

Let $\phi$ be a Hamiltonian diffeomorphism that has only a finite
number of fixed points (e.g., non-degenerate ones). We denote by
$J_0$ a compatible almost complex structure on $(M,\omega)$ and by
$\JJ_\omega$ the set of all such structures. For a given pair $(\phi,J_0)$,
we consider paths $J': [0,1] \to \JJ_\omega$ with
$$
J'(0) = J_0, \quad J'(1) = \phi^*J(0) \tag 1.13
$$
and denote the set of such paths by
$$
j_{(\phi,J_0)}.
$$
One can extend this family to $\R$ by the equation
$$
J'(t+1) = \phi^*J'(t).
$$
For each given $J' \in j_{(\phi,J_0)}$, we define the constants
$$
\aligned A_S(\phi,J_0;J') = \inf \{\omega([u]) & \mid  u: S^2 \to
M \text{ non-constant and} \\
& \text{ satisfying $\overline \part_{J'_t}u = 0$ for some $t \in [0,1]$}\}
\endaligned
$$
and
$$
A_S(\phi,J_0) = \sup_{J'\in j_{(\phi,J_0)}} A_S(\phi,J_0;J').
$$
As usual, we set $A_S(\phi,J_0;J') = \infty$ if there is no
$J'_t$-holomorphic sphere for any $t \in [0,1]$ e.g., as in the weakly
exact case. The positivity of $A_S(\phi,J_0;J')$ and so $A_S(\phi,J_0)$
is a consequence of the one parameter version of the uniform $\e$-regularity
theorem (see [SU], [Oh1]). We give the proof of positivity in section 4.

Next for each given $J'\in j_{(\phi,J_0)}$, we consider the
equation
$$
\cases {\part v \over \part \tau} + J'_t {\part v \over
\part t} = 0 \\
\phi(v(\tau,1)) = v(\tau,0), \quad \int_{\R \times [0,1]} |{\part
v \over \part \tau }|_{J'_t}^2 < \infty
\endcases
\tag 1.14
$$
in $v: \R \times [0,1] \to M$. Here the norm $|\cdot |_{J_t'}$ is
the norm induced by the Riemannian metric
$$
g_{J_t'} := \omega(\cdot, J_t' \cdot).
$$
The choice (1.13) of $J'$ enables us to interpret solutions of
(1.14) as pseudo-holomorphic sections of the mapping cylinder of
$\phi$ with respect to a suitably chosen almost complex structure
on the mapping cylinder. See subsection 2.3 for a more detailed
discussion of this.

We now define the  constant
$$
A_D(\phi,J_0;J'): = \inf_{v} \Big\{ \int v^*\omega,  \mid
\text{$v$ non-constant solution of (1.14)}\Big\}
$$
for each $J \in j_{(\phi,J_0)}$. In section 4, we will prove that
if $\phi$ is nondegenerate, $A_D(\phi,J_0;J') > 0$. Combining
$A_S(\phi,J_0;J')$ and $A_D(\phi,J_0;J')$, we define
$$
A(\phi,J_0;J') = \min\{A_S(\phi,J_0;J'),A_D(\phi,J_0;J')\}.
$$
Finally we set
$$
A(\phi,J_0) : = \sup_{J' \in j_{(\phi,J_0)}} A(\phi,J_0;J')
$$
and
$$
A(\phi) = \sup_{J_0\in \JJ_\omega} A(\phi,J_0). \tag 1.15
$$
Note that when $(M,\omega)$ is weakly exact, in which case
$A_S(\phi,J_0;J') = \infty$, $A(\phi,J_0)$ is reduced to
$$
A(\phi,J_0) = \sup_{J'\in j_{(\phi,J_0)}}\{ A_D(\phi,J_0;J') \}.
$$
It is now clear from the definition of $A(\phi,J_0)$
that $A(\phi,J_0) > 0$ and so we have $A(\phi)> 0$.
At this moment, however, we  note
that a priori $A(\phi,J_0)$ (and so $A(\phi)$) could be $+\infty$.

In section 5, we prove the following inequality, which
on the one hand gives rise to nondegeneracy of $\gamma$ and
on the other hand proves $A(\phi) < \infty$.

\proclaim{Theorem II} Let $\phi$ be nondegenerate and
$A(\phi)$ be the constant defined in (1.15). Then  we have
$$
\gamma(\phi) \geq A(\phi). \tag 1.16
$$
In particular $A(\phi)$ is finite and $\gamma(\phi) > 0$.
\endproclaim

A priori, Theorem II does not guarantee positivity of
$\gamma(\phi)$ for {\it degenerate} $\phi$: Recall that for the
definition of $A(\phi)$ we assumed that $\phi$ is nondegenerate.
However Theorem I (2) implies that the null-set
$$
\text{null}(\gamma) = \{\phi \in Ham(M,\omega) \mid
\gamma(\phi) = 0 \}
$$
is a normal subgroup of $Ham(M,\omega)$. Theorem II shows that
$\text{null}(\gamma)$ is a proper subgroup of $Ham(M,\omega)$
and so it must be a trivial subgroup thereof by Banyaga's theorem on
the simpleness of $Ham(M,\omega)$ [Ba]. This proves
the nondegeneracy of $\gamma$.

Exploiting the fact that the
definition of $\gamma$ involves only the identity class $1 \in
QH^*(M)$, we refine the definition of $A(\phi)$ to get a stronger
invariant denoted by $A(\phi;1)$, which we call the
{\it homological area} of $\phi$.
$A(\phi;1)$ always satisfies $A(\phi;1) \geq A(\phi)$
and is more computable than $A(\phi)$.
Referring to section 6 for the precise definition
of $A(\phi;1)$, we state

\proclaim{Theorem III} For any nondegenerate Hamiltonian
diffeomorphism $\phi$, we have
$$
\gamma(\phi) \geq A(\phi;1).
$$
\endproclaim

To provide some intuitive meaning to $A(\phi;1)$, we provide
one corollary of Theorem III which is
an immediate translation of the inequality $A(\phi;1)\leq \gamma(\phi)$
from the definition of $A(\phi;1)$.

\proclaim{Corollary} Let $\phi$ be nondegenerate and let
$q \in M \setminus \text{Fix }\phi$.
Then for any choice of $J_0 \in \JJ_\omega$ and
$J' \in j_{(\phi,J_0)}$, one of the following alternative holds:
\roster
\item there exists a non-constant $J'_t$-holomorphic sphere
passing through $q$ for some $t \in [0,1]$ and of area less
than or equal to $\gamma(\phi)$
\item (1.14) has a non-constant solution passing through
$q$ and of area less than or equal to $\gamma(\phi)$.
\endroster
\endproclaim

The inequality $A(\phi) \leq \gamma(\phi)$ in Theorem II
also implies a similar existence result but {\it without the phrase}
``passing through $q$'' which illustrates the fact that
the definition of $A(\phi;1)$ is tied to the fundamental
class $[M] = 1^\flat$. Therefore
$A(\phi;1)$ is a much better invariant than $A(\phi)$ for a
nondegenerate $\phi$. In fact when a sequence $\{\phi_i\}$ of
nondegenerate Hamiltonian diffeomorphisms converges to, say, a
{\it degenerate} $\phi \neq id$ in the $C^0$-topology,
$A(\phi_i;1)$ may fail to converge to zero together with $A(\phi_i)$.
In section 6, we will provide an example of a Hamiltonian
diffeomorphism $\phi$ on $S^2$ for which $0 < A(\phi) <
A(\phi;1)$ and also construct a sequence of nondegenerate
$\phi_i$ for which $A(\phi_i)$ converges to 0 but $A(\phi_i;1)$
remains constant. We conjecture in general that $A(\phi;1)$ is
continuous in $C^\infty$-topology (maybe even in
$C^0$-topology) and so can be extended to arbitrary smooth
Hamiltonian diffeomorphisms.

Noting $\|\phi\| \geq \gamma(\phi)$,
we also derive the following corollary.

\proclaim{Corollary}
The Hofer norm $\|\phi\|$ satisfies
$$
\|\phi\| \geq A(\phi;1). \tag 1.17
$$
\endproclaim

Using the spectral norm $\gamma(\phi)$ and the homological area
$A(\phi;1)$ of $\phi$, both of which depend only on $\phi$ and not
on its generating Hamiltonian function, we study the length
minimizing property of certain Hamiltonian paths {\it among all
paths} with their end points at the identity and $\phi$. For
simplicity of exposition, we will restrict to the case where $\phi$
is $C^1$-close to the identity in this paper,  leaving the proof
for the more general engulfable Hamiltonian diffeomorphisms
to a sequel [Oh10].

For the $C^1$-small case (or more generally for the engulfable
case), we can represent the subset $\hbox{graph}(\phi) \subset (M,
-\omega) \times (M, \omega)$ in the following form
$$
\Phi(\operatorname{graph}\phi) = \operatorname{graph} dS_\phi.
$$
where $\Phi: \UU \subset M \times M \to \VV \subset T^*\Delta$ is
a Darboux-Weinstein chart of the diagonal Lagrangian submanifold $\Delta
\subset (M - \omega) \times (M,\omega)$. The
function $S_\phi$ is unique up to the addition of a constant and
is called the {\it generating function} of $\phi$. For a function
$f: M \to \R$, we define
$$
osc(f) = \max f - \min f.
$$

\proclaim{Theorem IV} Suppose that $\phi$ is sufficiently
$C^1$-small. Then we have
$$
\gamma(\phi) = A(\phi;1) = osc(S_\phi)= \|\phi\|.
$$
\endproclaim

An immediate corollary of this theorem is the following result by
McDuff [Mc].

\proclaim{Corollary [Proposition 1.8, Mc]} Let $\phi$ be
sufficiently $C^1$-small and let $\phi^t$ be the Hamiltonian path
determined by the equation
$$
\Phi(\operatorname{graph}\phi^t) = \operatorname{graph} t dS_\phi.
$$
Then the path $t\in [0,1] \mapsto \phi^t \in Ham(M,\omega)$ is
a Hamiltonian path which is length minimizing among all paths from
the identity to $\phi$.
\endproclaim

The proof of Theorem IV, more specifically the inequality
$A(\phi;1) \geq \hbox{\rm osc}(S_\phi)$, is also a consequence of
an existence theorem for the equation (1.14) with a suitably chosen
asymptotic condition. Besides Theorem 5.4 this
existence proof heavily relies on the thick-thin decomposition
of the Floer moduli space. Such a decomposition exists when $H$
is sufficiently $C^2$-small. Study of the persistence of such a
decomposition will enable us to extend Theorem IV to the
cases where $\phi$ is not necessarily $C^1$-small. We will study
further properties of the homological area $A(\phi;1)$ and the
relevance of $A(\phi;1)$ to the energy-capacity inequality elsewhere.

Now the organization of the paper is in order. In section 2, we
briefly review various points needed throughout the paper
about chain level Floer theory
and spectral invariants of Hamiltonian diffeomorphisms. We
also set up our conventions concerning the definition of the
action functional and the gradings of quantum
cohomology and Floer homology. We then carefully review two
different set-ups of Floer homology from the literature as well
as the correspondences between them : one is the dynamical version of
perturbed Cauchy-Riemann equations [Fl3] and the other is
the mapping cylinder version [Fl1, DS].
A careful study of this correspondence is one of the crucial
elements in our construction and for the proof of nondegeneracy
of the norm $\gamma$. A similar correspondence was
used by the author in [Oh3,4] in the context of the
`open string' version of Floer homology i.e., in the context of
Lagrangian submanifolds.

After this review the paper is divided into two parts. Part I,
consisting of sections 3 - 7, contains the main applications of spectral
invariants mentioned above. Part II, consisting of sections 8 - 10,
contains key analysis needed in the proofs of the main
theorems stated in part I. Section 8 contains the proof of
Theorem 5.4, which relies on the analysis of
{\it adiabatic degeneration} in subsection 8.4.
In relation to the Fukaya category of the cotangent bundle,
similar analysis, but without quantum contributions, was carried out
in [FOh] in the context of degeneration of the moduli space of
pseudoholomorphic discs to the moduli space of graph flows.
The proof of the fundamental existence theorem
involves a rather delicate geometric analysis that is
not common in the literature, especially in the literature of
symplectic geometry. As far as we know,
[FOh] is the only other paper in the literature of symplectic geometry
which deals with analysis of adiabatic degeneration in details.
The paper [DS] by Dostoglou and Salamon studies
a different kind of adiabatic degeneration in relation to
the Atiyah-Floer conjecture, whose nature is somewhat
different from ours. In addition to this degeneration analysis,
proof of the lower estimate of $\gamma(\phi)$ in Theorem
III relies also on the chain level argument of the Floer theory
that is developed in [Oh5].

Sections 9 and 10 give the proof
of {\it thick-thin decomposition} of the Floer moduli spaces
needed for the proof of Theorem IV. This is the analog to
the similar decomposition result used by the author
[Oh2] in the study of the Maslov
class obstruction to the Lagrangian embeddings via the Floer
homology of Lagrangian submanifolds.
The analysis developed in part II have other applications
which will be presented elsewhere. The appendix provides an explicit
construction of a flat connection for the triple
$(H,\widetilde H\# \e f; \e f)$ of Hamiltonians
that is used in section 8.

To make the main stream of ideas transparent in this paper without
getting bogged down with technicalities related to transversality
issues of various moduli spaces, we will assume in the present
paper that $(M,\omega)$ is strongly semi-positive in the sense of
[Se] and [En1]. A closed symplectic manifold is called {\it
strongly semi-positive} if there is no spherical homology class $A
\in \pi_2(M)$ such that
$$
\omega(A) > 0, \quad2-n \leq c_1(A) < 0. \tag 1.18
$$
Under this condition the transversality problem concerning
various moduli spaces of pseudo-holomorphic curves is standard. We
will not mention this generic transversality question at all in
the main body of the paper unless it is absolutely necessary. The
scheme of the proof should go through on arbitrary closed
symplectic manifolds in the general framework of virtual moduli
cycles . In a sequel to the present paper we will incorporate
this enhanced machinery to generalize all the results  to
arbitrary closed symplectic manifolds.

For those who are familiar with the construction of spectral
invariants detailed in [Oh8] and want to get only the idea behind
the nondegeneracy proof for $\gamma$, we recommend them to read
section 5 and the heuristic discussion in subsection 8.2. The
nondegeneracy proof is based on the existence result of certain
pseudoholomorphic curves as used in [Oh4] : In [Oh4] the author
uses the existence result to prove positivity of the quantity
(1.12) when $\phi_H^1 \neq id$ on the cotangent bundle $T^*N$. To
prove this existence theorem in [Oh4], the author uses the cap
action realized by the degenerate pants in the context of Morse
homology and then proves that the cap action isomorphism $\cap \mu
: H_n(N) \to H_0(N)$ coincides with the standard pants product
$1*\mu = \mu$ in the Floer homology after applying the Poincar\'e
duality. The author's derivation of the latter fact in turn is
based on the analytical result from [FOh] mentioned above. In the
present paper, we use the nontrivial quantum product relation
`$1*1=1$' instead which suits better our definition (1.10) of
$\gamma$ and also produces such geometric lower bounds as those
obtained in Theorem III and IV, in addition to the nondegeneracy result for
$\gamma$.

The main results of the present work were presented in our paper
entitled ``Mini-max theory, spectral invariants and geometry of
the Hamiltonian diffeomorphism group'' (ArXiv:math.SG/0206092),
containing fewer details and some minor errors, which has
circulated since June, 2002. We isolate and expand the applications
of the spectral invariants in this paper with more
clarification and corrections of some details, leaving the
construction of them to a separate paper [Oh8].
Other applications
of the spectral invariants to the study of length minimizing
properties of Hamiltonian paths are given by the author [Oh7].

During the submission of the present paper, there
appeared a book by McDuff and Salamon
(``$J$-holomorphic Curves and Symplectic Topology'', AMS, Providence,
2004) in which the authors present a description of
the spectral invariants and a nondegeneracy proof in the rational case
based on the approach relying on the so called `PSS-isomorphism'
[PSS]. As we indicated in the above mentioned preprint and in [Oh8],
such an approach should be possible eventually.
However, to make the approach well-founded, it remains to write
out a proof of the isomorphism property of the PSS-map which is
used in various constructions carried out in [PSS] and in the book
mentioned above.

We thank the Institute for Advanced Study  for the excellent
environment and hospitality during our participation in the
2001-2002 program ``Symplectic Geometry and Holomorphic Curves''.
Much of the present work was finished during our stay in IAS. We
thank D. McDuff for useful communications at IAS. The final
writing of the previous version of the paper was carried out while
visiting the Korea Institute for Advanced Study in Seoul. We thank
KIAS for providing an uninterrupted quiet environment and time
for writing the paper and an excellent atmosphere for research.

Thanks to M. Entov and L. Polterovich for enlightening discussions
on spectral invariants and for explaining their applications
[En2], [EnP] of spectral invariants to study of the
Hamiltonian diffeomorphism group, and Y. Ostrover for explaining
his work from [Os] to us during our visit to Tel-Aviv University.
We also thank P. Biran and L. Polterovich for their invitation to
Tel-Aviv University and hospitality during our visit. We are very
much indebted to the referees for their careful reading of the
previous versions of the paper and for many helpful
questions, corrections and suggestions, which have both reinforced
the content of the paper and led to much improvement in the
presentation and readability of the paper. We also thank P. Spaeth
for his careful proofreading of the previous
version of the paper and many helpful comments in the presentation
and in the English of the paper, and N. Kieserman for the
final touch in English of the introduction of the paper.

\head{\bf \S 2. Preliminaries}
\endhead

\n {\it 2.1. The action functional and the action spectrum}
\smallskip

Let $(M,\omega)$ be any closed symplectic manifold and
$\Omega_0(M)$ be the set of contractible loops. We first recall
the construction from [HoS] of a covering space of $\Omega_0(M)$,
which we denote by $\widetilde \Omega_0(M)$.

Let $(\gamma,w)$ be a pair consisting of $\gamma \in \Omega_0(M)$
and a disc $w$ bounding $\gamma$. We say that $(\gamma,w)$ is {\it
$\Gamma$-equivalent} to $(\gamma,w^\prime)$ if
$$
\omega([w'\# \overline w]) = 0 \quad \text{and }\quad c_1([w'\#
\overline w]) = 0 \tag 2.1
$$
where $\overline w= w\circ c$,  $c: D^2 \to D^2$ is the complex
conjugation and $w'\# \overline w$ is the obvious glued sphere.
Here $\Gamma$ denotes the group
$$
\Gamma = {\pi_2(M)\over \text{ker\ } (\omega|_{\pi_2(M)}) \cap
\text{ker\ } (c_1|_{\pi_2(M)})}. \tag 2.2
$$
We denote by $[\gamma,w]$ the
$\Gamma$-equivalence class of $(\gamma,w)$ and by $\pi: \widetilde
\Omega_0(M) \to \Omega_0(M)$ the canonical projection. We also
call $\widetilde \Omega_0(M)$ the $\Gamma$-covering space of
$\Omega_0(M)$. The period group of $(M,\omega)$ is defined by
$$
\Gamma_\omega: = \{\omega(A) \mid A \in \pi_2(M) \} = \omega(\Gamma)
\subset \R. \tag 2.3
$$
It is either a discrete or a countable dense subset of $\R$.

Unless otherwise stated, we will always consider one-periodic
Hamiltonian functions $H: [0,1] \times M \to \R$ normalized by
$$
\int_M H_t\, d\mu = 0
$$
where $d\mu$ is the Liouville measure of $(M,\omega)$. We denote
by
$$
C^\infty_m([0,1] \times M)
$$
the set of such Hamiltonian functions, where ``$m$" stands for the
word ``mean zero".

When a periodic normalized Hamiltonian $H:  (\R/\Z) \times M  \to
\R$ is given, we define the action functional $\AA_H: \widetilde
\Omega_0(M) \to \R$ by
$$
\AA_H([\gamma,w])= -\int w^*\omega - \int H(t, \gamma(t))dt. \tag
2.4
$$
We denote by $\text{Per}(H)$ the set of periodic orbits of $X_H$.
\medskip

\definition{Definition 2.1}  We define the {\it
action spectrum} of $H$, denoted as $\hbox{\rm Spec}(H) \subset
\R$, by
$$
\hbox{\rm Spec}(H) := \{\AA_H(z,w)\in \R ~|~ [z,w] \in
\widetilde\Omega_0(M), z\in \text {Per}(H) \},
$$
i.e., the set of critical values of $\AA_H: \widetilde\Omega(M)
\to \R$.
\enddefinition

The following was proven in [Oh5].

\proclaim{Proposition 2.2} $\hbox{\rm Spec}(H)$ is a measure
zero subset of $\R$ for any $H$.
\endproclaim

We say that two Hamiltonians $F$ and $G$ are equivalent and denote
$F\sim G$ if they are connected by a one-parameter family of
Hamiltonians $\{F^s\}_{0\leq s\leq 1}$ such that $F^s \mapsto
\phi$ for all $s \in [0,1]$. We write $[F]$ or $[\phi,F]$ for the
equivalence class of $F$. Under the canonical one-one
correspondence between $H$ and its Hamiltonian path $\phi_F: t
\mapsto \phi_F^t$, $[F]$ corresponds to the path-homotopy class of
the Hamiltonian path $\phi_H: [0,1] \to Ham(M,\omega)$ with
fixed end points. If we denote by $\widetilde{Ham}(M,\omega)$
the set of path-homotopy classes of the based path at the identity
of $Ham(M,\omega)$, then the set of $[F]$ coincides with
$\widetilde{Ham}(M,\omega)$.

The following lemma was proven in the aspherical case in [Sc,
Po1]. We refer the reader to [Oh6] for complete details of its
proof in the general case.

\proclaim{Proposition 2.3}
If $F\sim G$, we have
$$
\text{\rm Spec}(F) = \text{\rm Spec}(G)  \tag 2.5
$$
as a subset of $\R$.
\endproclaim

\definition{Definition 2.4}
For given $h \in \widetilde{Ham}(M,\omega)$, we define
$$
\text{\rm Spec}(h) = \text{\rm Spec}(H)
$$
for a (and so any) $H$ with $[H] = h$.
\enddefinition

\medskip
\n {\it 2.2. Filtered Floer homology}
\smallskip

Suppose that $\phi \in Ham(M,\omega)$ is  nondegenerate in the
sense of Lefshetz fixed point theory: the derivative $T_p\phi:
T_pM \to T_pM$ has no eigenvalue one at each fixed point $p \in
M$. We call a Hamiltonian $H: S^1 \times M \to \R$
nondegenerate if $\phi_H^1 = \phi$ is  nondegenerate Hamiltonian
diffeomorphism.

For each nondegenerate $H:S^1 \times M \to \R $, we consider the
free $\Q$ vector space over
$$
\text{Crit}\AA_H = \{[z,w]\in \widetilde\Omega_0(M) ~|~ z \in
\text{Per}(H)\}. \tag 2.6
$$
To be able to define the Floer boundary operator correctly, we
need to complete this vector space downward with respect to the
real filtration provided by the action $\AA_H([z,w])$ of the
element $[z, w]$ of (2.6). More precisely, as in [Fl2, HoS], we
give the following definition.

\definition {Definition 2.5}
Consider the formal sum
$$
\beta = \sum _{[z, w] \in \text{Crit}\AA_H} a_{[z, w]}
[z,w], \, a_{[z,w]} \in \Q \tag 2.7
$$
\roster
\item We call those $[z,w]$ with $a_{[z,w]}
\neq 0$ {\it generators} of the sum $\beta$ and write
$$
[z,w] \in \beta.
$$
We also say that $[z,w]$ {\it contributes} to $\beta$ in that case.
\item We define the {\it support} of $\beta$ by
$$
\text{supp}(\beta): = \{ [z,w] \in \text{Crit }\AA_H \mid
a_{[z,w]} \neq 0 \text{ in the sum (2.7)}\}.
$$
\item We call the formal sum $\beta$
a {\it Novikov Floer chain} (or simply a {\it Floer chain}) if
it satisfies
$$
\#\Big(\text{supp}(\beta) \cap \{[z,w] \mid \AA_H([z,w])
\geq \lambda \}\Big) < \infty \tag 2.8
$$
for any $\lambda \in \R$. We denote by $CF_*(H)$
the set of Floer chains.
\endroster
\enddefinition

Note that $CF_*(H)$ is a graded $\Q$-vector space which is
infinite dimensional in general, unless $\pi_2(M) = 0$.
As in [Oh5], we introduce the following notion which is a crucial
concept for the mini-max argument that we carry out in this paper.

\definition{Definition 2.6}~  Let $\beta$ be a Floer chain in
$CF_*(H)$. We define the {\it level} of a chain $\beta$ and
denote it by
$$
\lambda_H(\beta) =\max_{[z,w]} \{\AA_H([z,w]) ~|~[z,w] \in
\text{supp}(\beta) \}
$$
if $\beta \neq 0$, and set $\lambda_H(0) = -\infty$ as usual.
\enddefinition

We briefly review construction of basic operators in the Floer
homology theory [Fl3] and their effects on the level of a given
Floer chain.

Let $J = \{J_t\}_{0\leq t \leq 1}$ be a periodic family of
compatible almost complex structures on $(M,\omega)$. For each
given pair $(H,J)$, we consider the perturbed Cauchy-Riemann equation
$$
\cases
\frac{\part u}{\part \tau} + J\Big(\frac{\part u}{\part t}
- X_H(u)\Big) = 0\\
\lim_{\tau \to -\infty}u(\tau) = z^-, \,  \lim_{\tau \to
\infty}u(\tau) = z^+. \\
\endcases
\tag 2.9
$$
This equation, when lifted to $\widetilde \Omega_0(M)$, is the
{\it negative} gradient flow of $\AA_H$ with respect to the
$L^2$-metric on $\widetilde \Omega_0(M)$ induced by the family of
metrics on $M$
$$
g_{J_t} = (\cdot, \cdot)_{J_t}: = \omega(\cdot, J_t\cdot);
$$
this $L^2$-metric is defined by
$$
\langle \xi, \eta \rangle_J: = \int_0^1 \ (\xi, \eta)_{J_t}\, dt.
$$
A simple computation then shows that the $L^2$-gradient vector
field of $\AA_H$ on $\widetilde \Omega_0(M)$, denoted by
$\hbox{grad}_J \AA_H$, is given by the formula
$$
(\hbox{grad}_J\AA_H)([\gamma,w]) = J(\dot \gamma -X_H(\gamma))
$$
from which we also derive that the gradient is projectable to
$\Omega_0(M)$.

We also denote by $\|v\|_{J_0}$ the norm of the vector $v$ with
respect to the metric $g_{J_0}$, i.e.,
$$
\|v\|^2_{J_0} = (v,v)_{J_0} = \omega(v, J_0 v) \tag 2.10
$$
for $v \in TM$.

Now suppose that $(H,J)$ is regular in the Floer theoretic sense.
For each given $[z^-,w^-]$ and $[z^+,w^+]$, we define the moduli space
$$
\MM(H,J; [z^-,w^-],[z^+,w^+])
$$
to be the set of solutions $u$ of (2.9) having finite energy and
satisfying the homotopy condition
$$
w^-\# u \sim w^+. \tag 2.11
$$
As a path in $\widetilde \Omega_0(M)$, we write the asymptotic
condition of $u$ as
$$
u(-\infty) = [z^-,w^-], \quad u(\infty) = [z^+,w^+].
$$
The Floer boundary map $\part =
\part_{(H,J)}: CF_*(H) \to CF_*(H)$, when the Conley-Zehnder
indices of $[z^-,w^-]$ and $[z^+,w^+]$ satisfy
$$
\mu_H([z^-,w^-]) - \mu_H([z^+,w^+]) = 1,
$$
is defined by its matrix coefficient
$$
\langle \part([z^-,w^-]), [z^+,w^+]\rangle : = \#\Big(\MM(H,J;
[z^-,w^-],[z^+,w^+])\Big)
$$
where $\#\Big(\MM(H,J; [z^-,w^-],[z^+,w^+])\Big)$ is the algebraic
number of the elements in $\MM(H,J; [z^-,w^-],[z^+,w^+])$, and is
set to be zero otherwise.

We refer to subsection 2.4 for more discussion about the energy of
the map $u$. The Floer boundary map
$$
\part = \part_{(H,J)}:CF_*(H) \to CF_*(H)
$$
has degree $-1$ and satisfies $\part\circ \part = 0$.

\definition{Definition 2.7} We say that a Floer chain $\beta \in CF(H)$
a  is {\it Floer cycle} if $\part \beta = 0$ and a {\it Floer
boundary} if $\beta = \part \delta$ for a Floer chain $\delta$.
Two Floer chains $\beta, \, \beta'$ are said to be {\it
homologous} if $\beta' - \beta$ is a boundary.
\enddefinition

When we are given a family $(\HH,j)$ with $\HH = \{H^s\}_{0\leq s
\leq 1}$ and $j = \{J^s\}_{0\leq s \leq 1}$, the chain
homomorphism
$$
h_{(\HH,j)}: CF_*(H^0) \to CF_*(H^1)
$$
is defined by considering the non-autonomous equation
$$
\cases \frac{\part u}{\part \tau} + J^{\rho(\tau)}\Big(\frac{\part
u}{\part t}
- X_{H^{\rho(\tau)}}(u)\Big) = 0\\
\lim_{\tau \to -\infty}u(\tau) = z^-, \, \lim_{\tau \to
\infty}u(\tau) = z^+
\endcases
\tag 2.12
$$
with the condition
$$
w^-\# u \sim w^+. \tag 2.13
$$
Here $\rho$ is a cut-off function of the type $\rho :\R \to [0,1]$,
$$
\align
\rho(\tau) & = \cases 0 \, \quad \text {for $\tau \leq -R$}\\
                    1 \, \quad \text {for $\tau \geq R$}
                    \endcases \\
\rho^\prime(\tau) & \geq 0
\endalign
$$
for some $R > 0$ (Here the monotonicity of $\rho$ is not essential
but provides a better estimate in general). $h_{(\HH,j)}$ has
degree 0 and satisfies
$$
\part_{(H^1,J^1)} \circ h_{(\HH,j)} = h_{(\HH,j)} \circ
\part_{(H^0,J^0)}.
$$
Two such chain maps $h_{(\HH^1,j^1)}, \, h_{(\HH^2,j^2)}$ are also
known to be chain homotopic [Fl3].

The following upper estimate of the action change can be proven by
the same argument as that of the proof of [Theorem 7.2, Oh3].

\proclaim{Proposition 2.8} Let $H$ and $K$ be any Hamiltonians,
not necessarily nondegenerate, $j = \{J^s\}_{s \in [0,1]}$ be any
given homotopy and $\HH^{lin} = \{H^s\}_{0\leq s\leq 1}$ be the
linear homotopy $H^s = (1-s)H + sK$. Suppose that (2.12) has a
solution with finite energy and satisfies (2.13). Then we have the
identity
$$
\align \AA_K([z^+,w^+]) & - \AA_H([z^-,w^-]) \\
& = - \int \Big|\dudtau \Big|_{J^{\rho_1(\tau)}}^2 -
\int_{-\infty}^\infty \rho'(\tau)\int_0^1\Big(K(t,u(\tau,t)) -
H(t,u(\tau,t))\Big)
\, dt\,d\tau  \tag 2.14\\
& \leq - \int \Big|\dudtau \Big|_{J^{\rho_1(\tau)}}^2 + \int_0^1
-\min_{x \in M} (K_t - H_t) \, dt \tag 2.15\\
& \leq  \int_0^1 -\min_{x \in M} (K_t - H_t) \, dt. \tag 2.16
\endalign
$$
\endproclaim

The following is an immediate corollary of Proposition 2.8.

\proclaim{Corollary 2.9}  Let $H$ and $K$ be nondegenerate
Hamiltonians, and let
$$
[z^+,w^+], \quad [z^-, w^-]
$$
be  critical points of $\AA_H$ and
$\AA_K$ respectively.  Suppose that there exists a family
$\{J^s\}_{0 \leq s \leq 1}$ such that the equation (2.12) for the
linear homotopy $\HH^{lin}$ has a solution $u$ satisfying (2.13)
and
$$
u(-\infty) = [z^-, w^-], \quad u(\infty) = [z^+, w^+]
$$
where $u$ is considered as a gradient trajectory on $\widetilde \Omega_0(M)$.
Then we have
$$
\AA_K([z^+,w^+]) - \AA_H([z^-,w^-]) \leq  \int_0^1 -\min_{x \in M}
(K_t - H_t)\,dt.
$$
In particular, we have
$$
\lambda_K(h_{(\HH,j)}(\beta)) \leq \lambda_H(\beta) + \int_0^1
-\min_{x \in M} (K_t - H_t)\,dt \tag 2.17
$$
where $h_{(\HH,j)}(\beta)$ is the image chain of $\beta$ under
the chain map $h_{(\HH,j)}$.
\endproclaim

\definition{Remark 2.10}
\roster \item Geometrically, what this corollary proves is that
under the Floer chain map, the level of Floer chains cannot
increase more than the change of the corresponding Hamiltonians,
more precisely, not more than $E^-(K - H)$. On the other hand,
there is no such lower bound for the amount of the decrease of the
level in general. Having a lower bound  is a manifestation of some
nontrivial homological property of Floer cycles.

\item We would like to remark that a similar calculation proves
there is also a uniform upper bound $C(\HH,j)$ for the chain map
over general homotopies $(\HH,j)$ or for the Floer chain homotopy
maps. Existence of such an upper estimate is crucial for the
construction of these maps. This upper estimate depends on the
choice of homotopy $(\HH,j)$ and is related to the curvature
estimates of the relevant Hamiltonian fibration (see [Po1],
[En1]). Whether the upper bound
$C(\HH,j)$ is the same as the upper bound (2.17) for the linear
homotopy $\HH^{lin}$ is closely tied to the minimality of the
Hamiltonian path $\HH$ with respect to the Hofer metric on
$Ham(M,\omega)$. (See [Oh7].)
\endroster
\enddefinition

Now we recall that $CF_*(H)$ is not only a $\Q$-module but also a
$\Lambda$-module where $\Lambda = \Lambda_\omega$ is the Novikov ring
of $(M,\omega)$: each $A \in \Gamma$ acts on
${\text Crit}\AA_H$ and so on $CF_*(H)$ by ``gluing a
sphere''
$$
A: [z,w] \to [z, w\# A].
$$
Then $\partial$ is $\Lambda$-linear and induces the standard Floer
homology $HF_*(H;\Lambda)$ with coefficients in $\Lambda$. (See
[HoS] for a detailed explanation of the Novikov ring
$\Lambda_\omega$ and its action on $CF_*(H)$.) However the action
does {\it not} preserve the filtration we defined above. Whenever
we talk about filtration, we will always presume that the relevant
coefficient ring is $\Q$. For a given nondegenerate $H$ and
$\lambda \in \R \setminus
\hbox{Spec}(H)$, we define the relative chain group
$$
CF_k^\lambda(H): = \{ \beta \in CF_k(H) \mid \lambda_H(\beta) <
\lambda\}.
$$

There is one more important aspect of the filtered Floer homology;
the behavior under symplectic conjugation
$$
\phi \mapsto \eta\phi \eta^{-1}; \quad Ham(M,\omega) \to Ham(M,\omega)
$$
for any symplectic diffeomorphism $\eta: (M,\omega) \to
(M,\omega)$. We summarize the basic facts relevant to the filtered
Floer homology here.

\roster \item When $H \mapsto \phi$, $\eta_*H \mapsto \eta \phi
\eta^{-1}$. \item If $H$ is nondegenerate, $\eta_*H$ is also
nondegenerate. \item If $(H,J)$ is regular in the Floer theoretic
sense, then so is $(\eta_*J, \eta_*H)$. \item There exists a
natural bijection $\eta_*: \Omega_0(M) \to \Omega_0(M)$ defined by
$$
\eta_*([z,w]) = ([\eta\circ z, \eta\circ w])
$$
under which we have the identity
$$
\AA_H([z,w]) = \AA_{\eta_*H}(\eta_*[z,w]).
$$
\item The $L^2$-gradients of the corresponding action functionals
satisfy
$$
\eta_*(\hbox{grad}_J\AA_H)([z,w]) =
\hbox{grad}_{\eta_*J}(\AA_{\eta_*H})(\eta_*([z,w])).
$$
\item If $u: \R \times S^1 \to M$ is a solution of perturbed
Cauchy-Riemann equation for $(H,J)$, then $\eta_*u = \eta\circ u$
is a solution for the pair $(\eta_*J, \eta_*H)$. In addition, all
the Fredholm properties of $(J,H,u)$ and $(\eta_*J,
\eta_*H,\eta_*u)$ are the same.
\endroster
These facts imply that conjugation by $\eta$ induces the canonical
filtration-preserving chain isomorphism
$$
\eta_*: (CF_*^\lambda(H), \part_{(H,J)}) \to (CF_*^\lambda
(\eta_*H), \part_{(\eta_*H,\eta_*J)})
$$
for any $\lambda \in \R \setminus \hbox{Spec}(H) = \R \setminus
\hbox{Spec}(\eta_*H)$.
In particular it induces a filtration preserving isomorphism
$$
\eta_*: HF_*^\lambda(H, J) \to HF_*^\lambda (\eta_*H,\eta_*J)
$$
in homology.

We close this section by fixing our grading convention for
$HF_*(H)$. Following the convention used in [Oh8], we will always
grade the complex $CF_*(H)$ by the Conley-Zehnder index,
denoted by $\mu_H([z,w])$ of its generators $[z,w]$.

To make a comparison of this grading with the gradings of quantum
or Morse cohomology, we recall that solutions of the {\it
negative} gradient flow equation of $-f$, i.e., of the {\it
positive} gradient flow of $f$
$$
\dot\chi - \text{grad } f(\chi) = 0
$$
correspond to the {\it negative} gradient flow of the action
functional $\AA_{\e f}$. This gives rise to the relation between
the Morse indices $\mu_{-\e f}^{Morse}(p)$ and the Conley-Zehnder
indices $\mu_{\e f}([p,\widehat p])$ where $\widehat p: D^2 \to
M$ is the constant disc $\widehat p(z) \equiv p$:
$$
\mu_{\e f}([p,\widehat p])  = \mu_{-\e f}^{Morse}(p)-n \tag
2.18
$$
in our convention [Oh8]. On the other hand, obviously we have
$$
\mu_{-\e f}^{Morse}(p) -n = (2n-\mu_{\e f}^{Morse}(p)) -n = n -
\mu_{\e f}^{Morse}(p). \tag 2.19
$$
Our grading convention of $CF_k(\e f)$ is
$$
k = \mu_{\e f}([p,\widehat p])\tag 2.20
$$
which related to the Morse index of $p \in \text{Crit }\e f$ by
(2.20). Now we relate this grading convention to the usual grading
convention of the quantum cohomology ring $QH^*(M)$.

This grading convention is translated to the standard grading of
the quantum cohomology via the duality map
$$
\flat: QH^{n-k}(M) \to QH_{n+k}(M)
$$
Under this grading
convention, we have the following grading preserving isomorphisms
$$
QH^{n-k}(M) \to QH_{n+k}(M) \to HQ_{n+k}(-\e f) \to HF_k(\e f)
\to HF_k(H). \tag 2.21
$$
(See [section 4.1, Oh8] for more explanations.) Recall that
when a nondegenerate Hamiltonian $H$ is given, we can represent
this Floer homology class by a Floer Novikov cycle
$\alpha \in CF_*(H)$. We denote the corresponding
Floer homology class by $a_H^\flat$ and mostly suppressing
$H$ from the notation as long as the meaning is clear.

Here we recall that $HQ_*(-\e f)$ is the Morse homology group
associated to the complex
$$
(CM_*(-\e f) , \part^{Morse})\otimes \Lambda
$$
a Morse theoretic realization of the quantum chain complex and
coincides with the Floer complex $(CF_{*-k}(\e f),
\part)$. This fact is one of the important ingredients
in the mini-max theory via the chain level Floer theory. We refer
to [Oh5,7-9] for the detailed illustration of previous usages of
this important fact, first observed by Floer himself [Fl3].

We have also shown in [Oh8] that with this grading convention the
pants product, denoted by $*$, has degree $-n$:
$$
*: HF_k(H) \otimes HF_\ell(K) \to HF_{(k+\ell)-n}(H\# K). \tag 2.22
$$
This suggests that it is also natural to consider the {\it shift
functor} $[n]$ and the shifted complex $CF[n](H)$ by
$$
CF[n]_k(H):= CF_{k+n}(H), \quad \part_{[n]} := (-1)^n\part
$$
and the corresponding homology $HF[n](H)$. In this shifted
complex, the pants product $*$ has degree zero, i.e., it satisfies
$$
*: HF[n]_k(H) \otimes HF[n]_\ell(K) \to HF[n]_{(k+\ell)}(H\# K).
$$
Finally we mention that our grading convention is consistent with
the quantum product
$$
\cdot : QH^a(M) \otimes QH^b(M) \to QH^{a+b}(M)
$$
having degree $0$, under the graded isomorphism (2.21).
\medskip

\noindent{\it 2.3. Spectral invariants and Floer fundamental cycles}
\smallskip

In this section, we briefly recall basic properties of the
spectral invariants constructed in [Oh5,9].

First let a nondegenerate $H$ be given and
choose a one-periodic $J$ such that the pair $(H,J)$ is
regular in the standard Floer theoretic sense. We will call any such
pair {\it Floer regular}.  For each given pair $(H,J)$ and $0\neq a
\in QH^*(M)$, we define
the following mini-max value of the action functional $\AA_H$
$$
\rho((H,J);a) = \inf_{\alpha}\{\lambda_H(\alpha) \mid \alpha \in
\ker\part_{(H,J)} \subset CF_*(H)\, \text{with }\, [\alpha] = a^\flat
\} \tag 2.23
$$
Finiteness of $\rho((H,J);a)$, i.e., $\rho((H,J);a) > -\infty$
is proven in [Oh8] and the numbers satisfy
$$
\rho((H,J);a) = \rho((H,J');a) \tag 2.24
$$
for two different choices of generic $J,\, J'$, which
is a consequence of (2.17).
We denote the common value by $\rho(H;a)$. Then it is proved in
[Oh8] that $\rho$ satisfies the following list of axioms.

\proclaim{Theorem 2.11 [Theorem I, Oh8]} Let $(M,\omega)$ be
arbitrary closed symplectic manifold. For any given quantum
cohomology class $0 \neq a \in QH^*(M)$, we have a continuous
function denoted by
$$
\rho =\rho(H;a): C_m^\infty([0,1] \times M) \times QH^*(M) \to \R
$$
such that they satisfy the following axioms: Let $H, \, F  \in
C_m^\infty([0,1] \times M)$ be smooth Hamiltonian functions and $a
\neq 0 \in QH^*(M)$. Then $\rho$ satisfies the following axioms:

\roster \item {\bf (Projective invariance)} $\rho(H;\lambda a) =
\rho(H;a)$ for any $0 \neq \lambda \in \Q$. \item {\bf
(Normalization)} For $a = \sum_{A \in \Gamma} a_A q^{-A} $, we
have $\rho(\underline 0;a) = v(a)$ where $\underline 0$ is the
zero function and
$$
v(a): = \min \{\omega(-A) ~|~  a_A \neq 0 \} = - \max \{\omega(A)
\mid a_A \neq 0 \}.
$$
is the (upward) valuation of $a$. \item {\bf (Symplectic
invariance)} $\rho(\eta_*H ;a) = \rho(H ;a)$ for any symplectic
diffeomorphism $\eta$ \item {\bf (Triangle inequality)} $\rho(H \#
F; a\cdot b) \leq \rho(H;a) + \rho(F;b) $ \item {\bf
($C^0$-continuity)} $|\rho(H;a) - \rho(F;a)| \leq \|H \# \overline
F\| = \|H - F\|$ where $\| \cdot \|$ is the Hofer's pseudo-norm on
$C_m^\infty([0,1] \times M)$. In particular, the function $\rho_a:
H \mapsto \rho(H;a)$ is $C^0$-continuous.
\endroster
\endproclaim

It is important to note that $\rho_a$ is defined for each {\it
function} $H$, whether it is nondegenerate or not, and
continuously extends to $C^0$ functions. Recall that the Lie
algebra of $Ham(M,\omega)$ as a vector space is `universal';
that is, independent of the symplectic form $\omega$. This enables
one to study the change of invariants $\rho(H;a)$ under the change
of symplectic form  for each {\it fixed} $H$.
This will be a subject of future study.

Here the symplectic invariance (3) follows from the discussion on
symplectic conjugation in subsection 2.2, which in particular
implies that $\eta_*$ induces a canonical filtration-preserving
bijection between the Floer cycles of $(H,J)$ and
$(\eta_*H,\eta_*J)$. This in particular implies that the mini-max
values $\rho(H;a)$ and $\rho(\eta_*H;a)$ are the same.

We next state the following theorem which will be used later
in the proof of the identity (5.5) for the class $a =1$.

\proclaim{Theorem 2.12 (Homotopy Invariance)} Assume that
the spectrality axiom (1.16) holds for nondegenerate
Hamiltonians. Let $F$ and $K$ be
arbitrary Hamiltonians satisfying $F \sim K$. Then we have
$$
\rho(F;a) = \rho(K;a).
$$
In particular, the function $F \mapsto \rho(F;a)$ pushes
down to a continuous function $\rho_a: \widetilde{Ham}(M,\omega)
\to \R$ in the natural topology on $\widetilde{Ham}(M,\omega)$.
\endproclaim
In [Oh8], we gave a proof of this homotopy invariance based on
the {\it spectrality axiom} for the rational $(M,\omega)$.
See the recent preprint (math.SG/0406449) of the author for
the proof of spectrality axiom for the irrational case.

To understand the value given by (2.23) in practice, it is crucial
to find a good Floer fundamental cycle $\alpha \in CF_n(H)$ for a
given regular $(H,J)$. An important task is to construct such a cycle.
One good way of constructing such a cycle is to transfer a
corresponding Morse cycle via the Floer chain map, which has been
frequently used by the author [Oh7-9]. We explain this procedure
now. First, we fix a generic Morse function $f$ and represent the
fundamental cycle $1^\flat = [M]$ by a Morse cycle of
$(-\e f,g_{J_0})$
$$
\alpha = \sum_j a_j [x_j,\widehat x_j], \quad a_j \in \Q
$$
with $x_j \in \text{Crit}_{2n} (-\e f)$ and $\widehat x_j$ is the
constant bounding disc. This particular fundamental cycle has the
important property that the cycle cannot be pushed down in the
sense of [Oh4]. We recall the following definition from [Oh7]
which will be useful for our later discussions.

\definition{Definition 2.13} Let $(H,J)$ be a Floer regular pair
so that the Floer complex $(CF_*(H), \part_{(H,J)})$ is
well-defined. Let $\alpha$ be a Floer cycle of $H$ and $a \in
QH^*(M)$ be the corresponding quantum cohomology class with
$[\alpha] = a^\flat$. We call the Floer cycle $\alpha$ {\it tight}
if it realizes the mini-max value, i.e.,
$$
\lambda_H(\alpha) = \rho(H;a)
$$
and we call the corresponding critical value a {\it homologically
essential critical value} of $\AA_H$.
\enddefinition

Using the fact that, when $\e > 0$ is sufficiently small, this
Morse cycle also provides a Floer cycle of $\e f$ that represents
the quantum cohomology class $1$, we transfer $\alpha$ to
$$
\alpha_H = h_{(\HH,j)}(\alpha) \tag 2.25
$$
by the chain map $h_{(\HH,j)}: CF_n(\e f) \to CF_n(H)$
along a homotopy $(\HH,j)$ to the cycle $h_{(\HH,j)}(\alpha_H)$.
This will then represent a Floer fundamental cycle for $H$
by the invariance of the Floer
homology under the homotopy of Hamiltonians.

In this paper and others [Oh5,7-9], we have often used the linear
homotopy $\HH$ defined by
$$
\HH: s \mapsto (1-s) \e f + s H \tag 2.26
$$
to represent a Floer fundamental cycle of the given Hamiltonian
$H$ in this way.  Following [Oh7], we call this particular
fundamental cycle a {\it canonical Floer fundamental cycle}.
We emphasize that {\it this cycle, however, depends on the initial
Morse function $f$} (and also on the choice of the homotopy
$(\HH,j)$). In general the transferred cycles will
{\it not} be tight. Understanding under what condition of $H$
the transferred cycle
remains tight provides much information on the length minimizing
property of Hofer's geodesics. (See [Oh7] for further discussions
on this issue.)

\medskip

\noindent{\it 2.4. Comparison of two Cauchy-Riemann equations}
\smallskip

In this section, we explain the relation between Floer's standard
perturbed Cauchy-Riemann equation for $u:\R \times S^1 \to M$
$$
\cases
\dudtau + J_t\Big(\dudt - X_H(u)\Big) = 0 \\
\int |{\part v \over \part \tau }|_{J_t}^2 < \infty,
\endcases
\tag 2.27
$$
and its mapping cylinder version $v: \R \times \R \to M$
$$
\cases {\part v \over \part \tau} + J'_t {\part v \over
\part t} = 0 \\
\phi(v(\tau,t+1)) = v(\tau,t),
\quad \int |{\part v \over \part \tau }|_{J'_t}^2 < \infty
\endcases
\tag 2.28
$$
where $\phi = \phi_1^H$.  We often restrict $v$ to $\R \times
[0,1]$ and consider it as a map from $\R \times [0,1]$ that
satisfies $\phi(v(\tau,1)) = v(\tau,0)$. The present author
previously exploited a similar correspondence in [Oh3,4] in the
`open string' context of Lagrangian submanifolds for the same
purpose, and call the former version of Floer homology the {\it
dynamical} and the latter {\it geometric}. We do the same here.

Before explaining the correspondence, we would like to point out
one subtle difference between the case of the `open string'
version [Oh3,4] and the current `closed string' version. In the
point of view of transforming the equations (2.27), (2.28) from
one to the other, we {\it must} consider a periodic family
$\{J_t\}_{0\leq t \leq 1}$ for (2.27), which was not necessary for
the open string version.  Then for any given solution
$u=u(\tau,t): \R \times S^1 \to M$, we `open up' $u$ along
$t=0\equiv 1$ and define the map
$$
v: \R \times [0,1] \to M
$$
by
$$
v(\tau,t) = (\phi_H^t)^{-1}(u(\tau,t)) \tag 2.29
$$
and then extend to $\R$ so that $\phi(v(\tau, t+1)) = v(\tau, t)$.
A simple computation shows that the map satisfies (2.28), {\it
provided} the family $J' =\{J'_t\}_{0 \leq t \leq 1}$ is defined
by
$$
J'_t = (\phi_H^t)^*J_t
$$
for the given periodic family $J_t$ used for the equation (2.27).
By definition, this family $J'$ of almost complex structure satisfies
$$
J'(t+1) = \phi^*J'(t). \tag 2.30
$$
The condition
$$
\phi(v(\tau,t+1)) = v(\tau,t)
$$
enables us to consider the map
$$
v: \R \times \R \to M
$$
as a pseudo-holomorphic section of the `mapping cylinder'
$$
E_\phi: = \R \times M_\phi = \R \times \R \times M/(\tau,t+1,\phi(x))
\sim (\tau,t,x)
$$
where $M_\phi$ is the mapping circle defined by
$$
M_\phi:=  \R \times M/(t+1, \phi(x)) \sim (t, x).
$$
Note that the product symplectic form $d\tau \wedge dt + \omega$
on $\R \times \R \times M$ naturally projects to $E_\phi$ since
$\phi$ is symplectic, and so $E_\phi$ has the structure of a
Hamiltonian fibration. In this setting, $v: \R \times \R \to M$
can be regarded as the section $s: \R \times S^1 \to E_\phi$
defined by
$$
s(\tau,t) = [\tau,t,v(\tau,t)]
$$
which becomes a pseudo-holomorphic section of $E_\phi$ for a
suitably defined almost complex structure. (See [section 3, Oh9]
for the precise formulations.)

One advantage of the mapping cylinder version over the more
standard dynamical version (2.27) is that its dependence on the
Hamiltonian $H$ is much weaker than in the latter. Indeed, this
mapping cylinder version can be put into the general framework of
Hamiltonian fibrations with given fixed monodromy of the fibration
at infinity as in [En1]. This framework turns out to be essential
to prove the triangle inequality of the spectral invariants. (See
section 5).

Another important ingredient is the formula for the energy: For
the given {\it periodic} family $J=\{J_t\}_{0 \leq t \leq 1}$, we
define the energy of the map $u: \R \times S^1 \to M$ by
$$
E_J(u) = \frac{1}{2} \int_{\R \times S^1}
\Big(\Big|\dudtau\Big|_{J_t}^2 + \Big|\dudt - X_H(u)\Big|_{J_t}^2
\Big)dt\,d\tau
$$
for (2.27). Similarly for the given $J' = \{J'_t\}_{0 \leq t \leq 1}
\in j_{(J_0,\phi)}$,
we define the energy of the map $v: \R \times [0,1] \to M$ by
$$
E_{J'}(v) = \frac{1}{2} \int_{\R \times [0,1]}
\Big(\Big|\dvdtau\Big|_{J'_t}^2 + \Big|\dvdt\Big|_{J'_t}^2\Big) dt\,d\tau.
$$
This energy is the {\it vertical} part of the energy of the
section $s:\R \times S^1 \to E_\phi$ defined above with respect to
a suitably chosen almost complex structure $\widetilde J$ on
$E_\phi$.  (See [section 3, Oh9] for more explanation.) Note that
because of (2.29)-(2.30), one can replace the  domain of
integration $\R \times [0,1]$ by any {\it fundamental domain} of
the covering projection
$$
\R \times \R \to \R \times (\R/\Z)
$$
without changing the integral. The choice of $\R \times [0,1]$ is
one such  choice.

The following identity plays an important role in the proof of the
nondegeneracy of the invariant norm  we construct later. The proof
is a straightforward computation and left to the readers.

\proclaim{Lemma 2.14} Let $J = \{J_t\}_{0 \leq t \leq 1}$ be a
periodic family and define $J' = \{J'_t\}_{0\leq t\leq 1}$ by
$$
J'_t = (\phi_H^t)^*J_t.
$$
Let $u: \R \times S^1 \to M$ be any smooth map and $v: \R \times [0,1]
\to M$ be the map defined by
$$
v(\tau,t) = (\phi_H^t)^{-1}(u(\tau,t)).
$$
Then we have
$$
E_J(u) = E_{J'}(v). \tag 2.31
$$
\endproclaim
\medskip

\centerline{\it PART I: APPLICATIONS OF SPECTRAL INVARIANTS}
\medskip

\head{\bf \S 3. Definition of the spectral norm}
\endhead

In this section, we prove Theorem II in the introduction.
We first recall
$$
\gamma(H) = \rho(H;1) + \rho(\overline H;1) \geq 0. \tag 3.1
$$
We now relate $\rho(H;1)$ with the negative Hofer length $E^-(H)$.

\proclaim{Theorem 3.1} For any $H$ and $0 \neq a \in QH^*(M)$, we
have
$$
-E^+(H) + v(a) \leq \rho(H;a) \leq E^-(H) + v(a). \tag 3.2
$$
In particular for any classical cohomology class $b \in H^*(M)
\hookrightarrow QH^*(M)$, we have
$$
-E^+(H) \leq \rho(H;b) \leq E^-(H) \tag 3.3
$$
for any Hamiltonian $H$.
\endproclaim
\demo{Proof} We first recall  the following general inequality
$$
\int_0^1 - \max(H-K)\, dt  \leq \rho(H;a) -
\rho(K;a) \leq  \int_0^1 - \min(H-K)\, dt. \tag 3.4
$$
proven in [Oh 9], which can be rewritten as
$$
\rho(K;a)  +   \int_0^1 - \max(H-K)\, dt \leq \rho(H;a) \leq
\rho(K;a)  +   \int_0^1 - \min(H-K)\, dt.
$$
Now let $K \to 0$ which results in
$$
\rho(0;a) + \int_0^1 -\max(H)\, dt \leq
\rho(H;a) \leq \rho(0;a)  +   \int_0^1 - \min(H)\, dt. \tag 3.5
$$
By the normalization axiom from subsection 2.3, we have $\rho(0;a)
= v(a)$ which turns (3.5) to
$$
v(a) - E^+(H) \leq \rho(H;a) \leq v(a) + E^-(H)
$$
for any $H$. (3.3) immediately follow from the definitions and the
identity $v(b) = 0$ for a classical cohomology class $b$.
This finishes the proof. \qed\enddemo

Applying (3.3) to $b = 1$, we derive $\rho(H;1) \leq E^-(H)$
and $\rho(\overline H;1) \leq E^-(\overline H)$
for arbitrary $H$, and hence
$$
\gamma(H) \leq \|H\|.  \tag 3.6
$$
The nonnegativity (3.1) leads us to the following definition.

\definition{Definition 3.2} We define $\gamma: Ham(M,\omega)
\to \R_+$ by
$$
\gamma(\phi):= \inf_{H \mapsto \phi} (\rho(H;1) + \rho(\overline
H;1))
$$
and call this the {\it spectral (pseudo)-norm}.
\enddefinition

The following theorem shows that $\gamma$ is indeed a norm on
$Ham(M,\omega)$ that is invariant under the conjugate action of
symplectic diffeomorphisms which in turn implies the biinvariance
of the norm.

\proclaim{Theorem 3.3} Let $\gamma$ be as above. Then $\gamma :
Ham(M,\omega) \to \R_+$  enjoys
the following properties. \roster \item $\phi= id$ if and only if
$\gamma(\phi) = 0$ \item $\gamma(\eta^{-1}\phi \eta) =
\gamma(\phi)$ for any symplectic diffeomorphism $\eta$ \item
$\gamma(\phi\psi) \leq \gamma(\phi) + \gamma(\psi)$ \item
$\gamma(\phi^{-1}) = \gamma(\phi)$ \item $\gamma(\phi)\leq
\|\phi\|$
\endroster
\endproclaim

In the remaining section, we will give the proofs of these
statements postponing the most non-trivial statement,
nondegeneracy, to section 5. We split our
proof of Theorem 3.3 item by item.

\demo{Proof of $(2)$} We recall the symplectic invariance of
spectral invariants $\rho(H;a) = \rho(\eta_*H;a)$. Applying this
to $a = 1$, we derive the identity
$$
\align
\gamma(\phi) & = \inf_{H \mapsto \phi}\Big(\rho(H;1) + \rho(\overline H;1)\Big)\\
& = \inf_{H \mapsto \phi}\Big(\rho(\eta_*H;1) +
\rho(\overline{\eta_*H};1)\Big) = \gamma(\eta^{-1}\phi \eta),
\endalign
$$
which finishes the proof.
\qed \enddemo

\demo{Proof of  $(3)$} We first recall the triangle inequality
$$
\rho(H\# K;1) \leq \rho(H;1) + \rho(K;1) \tag 3.7
$$
and
$$
\rho(\overline K \# \overline H;1) \leq \rho(\overline K;1) + \rho(\overline H;1).
\tag 3.8
$$
Adding up (3.7) and (3.8), we have
$$
\aligned
\Big(\rho(H\# K;1)  & +\rho(\overline H \# \overline K;1) \Big) \\
& \leq \Big(\rho(H;1) +\rho(\overline H;1) \Big)+  \Big(\rho(K;1)
+ \rho(\overline K;1)\Big).
\endaligned
\tag 3.9
$$
Now let $H \mapsto \phi$ and $K \mapsto \psi$.
Because $H\# K$ generates $\phi\psi$ and $
\overline K \# \overline H$ generates $\psi^{-1}\phi^{-1} = (\phi\psi)^{-1}$,
we have
$$
\gamma(\phi\psi) \leq \rho(H\# K;1)  +\rho(\overline H \# \overline K;1)
$$
and hence
$$
\gamma(\phi\psi) \leq  \Big(\rho(H;1) +\rho(\overline H;1) \Big)
+  \Big(\rho(K;1)
+ \rho(\overline K;1)\Big)
$$
from (3.9). By taking the infimum of the right hand side over all
$H \mapsto \phi$ and $K \mapsto \psi$, (3) is proved. \qed\enddemo

\demo{Proof of  $(4)$}
The proof immediately follows from the observation that the definition of
$\gamma$ is symmetric over the map $\phi \mapsto \phi^{-1}$.
\qed\enddemo

\demo{Proof of $(5)$}  This is an immediate consequence of (3.6)
with its infimum taken over $H \mapsto \phi$. \qed\enddemo

It now remains to prove nondegeneracy of $\gamma$, which we
will do in section 5 modulo the proof of an existence theorem,
Theorem 5.4. Proof of the latter will be given in section 8.

\head{\bf 4. The $\e$-regularity theorem and its consequences}
\endhead

As in the proof of nondegeneracy of Hofer's norm in [Ho]
and [LM], we will prove nondegeneracy of $\gamma$ by comparing
it with a geometric invariant of $\phi$ independent of
its generating Hamiltonians. In [Ho] and [LM], the authors use
the maximal Gromov area $\lambda > 0$ of symplectic balls
$B(\lambda)$ and then prove the lower bound of $\|H\|$
for $H$s which satisfy $\phi_H^1(B(\lambda)) \cap B(\lambda) = \emptyset$.
Here we denote by $B(\lambda)$ the image of any symplectic
embedding $g: B^{2n}(r) \subset \C^n \to (M,\omega)$, and
call the image $g(0)$ the center of the symplectic embedding
$g$ or $B(\lambda)$. It seems to be useful to formalize this
invariant of general Hamiltonian diffeomorphisms in the following way.

\definition{Definition 4.1} Let $\phi$ be an arbitrary
Hamiltonian diffeomorphism. We define
$$
bArea(\phi):=\sup_{B} \{\lambda=\pi r^2 \mid \phi(B(\lambda))
\cap B(\lambda) =\emptyset \}
$$
and call it the {\it ball-area} of $\phi$.
\enddefinition

For the lower estimate of $\gamma(\phi)$, we will use another
geometric invariant of $\phi$, instead of the ball area, which much
better suits our chain level Floer theory.
The rest of this section is occupied by the description
of this invariant and the proof of its nontriviality.

Suppose $\phi$ is a nondegenerate Hamiltonian diffeomorphism.
Denote by $J_0$ a compatible almost complex structure on
$(M,\omega)$. For given such a pair $(\phi,J_0)$, we consider paths $J':
[0,1] \to \JJ_\omega$ with
$$
J'(0) = J_0, \quad J'(1) = \phi^*J_0 \tag 4.1
$$
and define $j_{(\phi,J_0)}$ to be the set of such paths. We extend
$J'$ to $\R$ so that
$$
J'(t+1) = \phi^*J'(t).
$$

For each given $J' \in j_{(\phi,J_0)}$, we define the constant
$$
\aligned A_S(\phi,J_0;J') = \inf \Big\{\omega([u]) & \mid  u: S^2
\to M \text{
non-constant and} \\
& \text{ satisfying $\overline \part_{J'_t}u = 0$ for some $t \in
[0,1]$} \Big\}.
\endaligned
\tag 4.2
$$
A priori it is not obvious whether $A_S(\phi,J_0;J')$ is not zero.
We will prove this using the so called {\it $\e$-regularity
theorem}.  This kind of $\e$-regularity theorem is one of the
fundamental ingredients, first introduced by Sachs and Uhlenbeck
[SU], for the study of compactness and regularity theory in the
geometric analysis of harmonic maps, minimal surfaces and others
including the Gromov compactness theorem of pseudo-holomorphic
maps. We state a parametric version of the $\e$-regularity theorem
stated in [Oh1] in the context of pseudo-holomorphic curves.

Let $g$ be any almost K\"ahler metric of $(M,\omega)$, i.e., $g =
\omega(\cdot, J_{ref}\cdot)$ for some fixed almost complex structure
$J_{ref}$ compatible  to $\omega$. Denote by $|\cdot|$ the norm associated to
this fixed metric. {\it We will always assume that such an almost
K\"ahler metric $g$ is given and all the norms used are in
terms of this fixed metric, unless otherwise said}.
We denote the derivative of $u$ by $Du$.

\proclaim{Lemma 4.2 ($\e$-Regularity Theorem) [Proposition 3.3,
Oh1]} We denote by $D= D^2(1)$ the unit open disc . Let $J_0$ be
any almost complex structure and let $u: D \to M$ be a
$J_0$-holomorphic map. Then there exists some $\e = \e(g, J_0)
> 0$ such that if $\int_D |Du|^2 < \e$, then for any smaller disc
$D' = D^2(r)$ with $\overline D' \subset D$, we have
$$
\|Du\|_{\infty, D'}: = \max_{z\in D'}|Du(z)| \leq C
$$
where $C >0$ depends only on $g, \, \e, J_0$ and $D'$, not on
$u$. Furthermore, the same $C^1$-bound holds for any compact
family $K$ of compatible almost complex structures this time with
$\e=\e(g,K)$ and $C=C(g, \e, K, D')$ depending on $K$.
\endproclaim

An immediate corollary of this $\e$-regularity theorem is the
following uniform $C^1$-estimate of pseudo-holomorphic curves.

\proclaim{Corollary 4.3}  Let $J' \in j_{(\phi,J_0)}$. Then there
exists an $\e = \e(J') > 0$ such that if $\omega(u)  < \e$, then
we have
$$
\|Du\|_\infty:=\max_{z\in S^2}|Du(z)| \leq C \tag 4.3
$$
for any $J'_t$-holomorphic sphere $u: S^2 \to M$ and for any $t
\in [0,1]$ where $C=C(\e, J')$ does not depend on $u$.
\endproclaim
\demo{Proof} Let $u$ be $J_{t_0}$-holomorphic for $t_0 \in [0,1]$.
We choose an open covering of $S^2$ by a finite number of discs
$\{D_\alpha\}$ with a fixed size. Obviously we have
$$
\frac{1}{2}\int_{D_\alpha} |Du|^2 \leq \frac{1}{2} \int_{S^2}
|Du|^2. \tag 4.4
$$
On the other hand, since $M$ and $[0,1]$ are compact, there exists
a uniform constant $K> 0$ such that
$$
\frac{1}{\sqrt{K}} |\cdot| \leq |\cdot|_{J'_{t}} \leq \sqrt{K}
|\cdot| \tag 4.5
$$
for all $t \in [0,1]$, where we recall $|\cdot|$ is the norm
measured in terms of the metric $g$ associated to the fixed
$J_{ref}$. Therefore combining (4.4) and (4.5), we derive
$$
\frac{1}{2}\int_{D_\alpha} |Du|^2 \leq \frac{K}{2} \int_{S^2}
|Du|_{J_{t_0}}^2 = K \omega(u). \tag 4.6
$$
Here we used the well-known identity
$$
\frac{1}{2} \int_{S^2} |Du|_{J_{0}}^2 = \omega(u)
$$
for $J_0$-holomorphic map $u: S^2 \to M$ where $J_0$ is compatible
to $\omega$. We emphasize that for this identity to hold the norm
$|\cdot|_{J_0}$ and the pseudo-holomorphic map $u$ must be
associated to the {\it same} compatible almost complex structure
$J_0$. This is the reason why we need to use (4.5) to relate the
harmonic energy $\frac{1}{2} \int_{S^2} |Du|^2$ of
$J_{t_0}$-holomorphic map $u$, which is measured by the reference
metric $g = \omega(\cdot, J_{ref}\cdot)$, to the symplectic area
$\omega(u)$.

Since $\int_{D_\alpha} |Du|^2$ is invariant under the conformal
change of metric of the domain, we may assume the radius of
$D_\alpha$ is one for all $\alpha$. Choose $\e > 0$ so that the
$C^1$-bound in Lemma 4.2 holds for $K \e$. Then we choose a
smaller discs $D'_\alpha$ so that $\overline D'_\alpha \subset
D_\alpha$ and $\{D'_\alpha\}$ still covers $S^2$ for each
$\alpha$. Then we apply Lemma 4.2 for each $\alpha$ and obtain
(4.3). \qed\enddemo

With Corollary 4.3 at our disposal, we prove the following strict
positivity result. We note that we cannot rule out the possibility
$A_S(\phi,J_0;J') = \infty$ as in the weakly exact case where
$A_S(\phi,J_0;J') = \infty$ for all choices of $(J_0;J')$.

\proclaim{Proposition 4.4} Let $\phi,$ $J_0$ and $J'$ be as above.
Then we have
$$
A_S(\phi,J_0;J')> 0.
$$
\endproclaim
\demo{Proof} Suppose $A(\phi, J_0;J') = 0$. Then there exists a
sequence $t_j \in [0,1]$ and a sequence of non-constant maps $u_j:
S^2 \to M$ such that $u_j$ is $J_{t_j}$-holomorphic and
$$
\omega(u_j) = E_{J_{t_j}}(u_j) \to 0
$$
as $j \to \infty$. By choosing a subsequence of $t_j$, again
denoted by $t_j$, we may assume that $t_j \to t_\infty \in [0,1]$
and so $J_{t_j}$ converges to $J_{t_\infty}$ in the
$C^\infty$-topology. We choose a sufficiently large $N \in \Z_+$ so
that
$$
E_{J_{t_j}}(u_j) = \omega(u_j) < \e(J')
$$
for all $j \geq N$, where $\e(J')$ is the constant $\e$ provided
in Corollary 4.3. Then we have the uniform $C^1$-bound
$$
0 < \|Du_j\|_\infty \leq C(\e, J').
$$
The Ascoli-Arzela theorem then implies that there exists a
subsequence, again denoted by $u_j$, such that $u_j$ converges
uniformly to a continuous map $u_\infty: S^2 \to M$. Recalling
that all the $u_j$ are $J_{t_j}$-holomorphic and $J_{t_j}$
converges to $J_{t_\infty}$ in the $C^\infty$-topology, the
standard boot-strap argument implies that $\{u_j\}$ converges to
$u_\infty$ in the $C^1$ topology (and so  in the
$C^\infty$-topology). However we have
$$
E_{J_{t_\infty}}(u_\infty) = \lim_{j \to \infty}E_{J_{t_j}}(u_j) = 0
$$
and hence $u_\infty$ must be a constant map, say $u_\infty\equiv x
\in M$. Therefore $\{u_j\}$ converges to the point $x$ in the
$C^\infty$-topology. In particular, if $j$ is sufficiently large,
then the image of $u_j$ is contained in a ({\it contractible})
Darboux neighborhood of $x$. Therefore we must have $\omega([u_j])
= 0$ and in turn
$$
E_{J_{t_j}}(u_j) = 0
$$
for all sufficiently large $j$, because  $E_{J_{t_j}}(u)=
\omega(u) $ holds for any $J_{t_j}$-holomorphic curve $u$. This
contradicts the assumption that $u_j$ is non-constant. This
finishes the proof. \qed\enddemo

Next for each given $J'\in j_{(\phi,J_0)}$, we consider the
equation of $v: \R \times \R \to M$
$$
\cases {\part v \over \part \tau} + J'_t {\part v \over
\part t} = 0 \\
\phi(v(\tau,t+1)) = v(\tau,t), \quad \int_{\R \times [0,1]}
|{\part v \over \part \tau}|_{J'_t}^2 < \infty.
\endcases
\tag 4.7
$$

\proclaim{Lemma 4.5}  Assume that $\phi$ is nondegenerate. For any
solution $v$ of (4.7), we consider the paths
$$
v(\tau): [0,1] \to M, \quad v(\tau)(t):= v(\tau,t).
$$
Then $v(\tau)$ uniformly converges to a fixed point $p \in
\hbox{Fix }\phi$ as $\tau \to \pm\infty$.  We simply say
$$
v(\pm\infty) \in \text{Fix }\phi. \tag 4.8
$$
\endproclaim
\demo{Proof} We will prove this only as $\tau \to \infty$,
since the case as $\tau \to -\infty$ will be the same.
The finite energy condition
$$
\int_{\R \times [0,1]} \Big|{\part v \over \part
\tau}\Big|_{J'_t}^2 < \infty
$$
implies that there exists a sequence $\tau_j \to \infty$ such that
$$
\int_0^1 \Big|{\part v \over \part \tau}(\tau_j,
\cdot)\Big|_{J'_t}^2\, dt \to 0
$$
as $j \to \infty$.  Together with the nondegeneracy assumption on
$\phi$ and the choice of $J'$, this implies an exponential decay
of $|{\part v \over \part \tau}(\tau_j, \cdot)|_{J'_t}^2$ as $\tau
\to \infty$. (See [Fl1, SZ] for such a proof.) Therefore we have
$$
\int_0^1 \Big|{\part v \over \part \tau}(\tau,
\cdot)\Big|_{J'_t}^2\, dt \to 0 \quad \hbox{as  $\tau \to
\infty$}.
$$
Since $M$ is compact, using the compactness of the Sobolev
embedding $W^{1,2}([0,1]) \hookrightarrow C^0([0,1])$ and then the exponential
decay property, we derive that $v(\tau): [0,1] \to M$ converges
uniformly on $[0,1]$ to a constant point $p \in M$. Due to the
boundary condition $\phi(v(\tau, t+1)) = v(\tau,t)$, we have
$\phi(v(\tau,1)) = v(\tau,0)$ which immediately implies that the
point $p$ must satisfy $\phi(p) = p$, i.e., $p$ is a fixed point
of $\phi$. This finishes the proof. \qed\enddemo

Now it is a crucial matter to produce a non-constant solution of
(4.7). For this purpose, using the fact that $\phi \neq id$, we
choose a symplectic ball $B(\lambda)$ such that
$$
\phi(B(\lambda)) \cap B(\lambda) = \emptyset \tag 4.9
$$
where $B(\lambda)$ is the image of a symplectic embedding
$g: B^{2n}(r) \to B(\lambda) \subset M$
of the standard Euclidean ball $B^{2n}(r) \subset \C^n$ of radius
$r$ with $\lambda = \pi r^2$. We then study (4.7) together with
$$
v(0,0) \in B(\lambda). \tag 4.10
$$
Because of (4.8) and (4.9), it follows
$$
v(\pm\infty) \in \text{Fix }\phi \subset M \setminus B(\lambda).
\tag 4.11
$$
Therefore any such solution cannot be constant.

We now define the constant
$$
A_D(\phi,J_0;J'): = \inf_{v} \Big\{ \int_{\R \times [0,1]}
v^*\omega\,  \Big| \, \text{$v$ non-constant solution of (4.7)}\Big\}
$$
for each $J' \in j_{(\phi,J_0)}$. Obviously
we have $A_D(\phi,J_0;J') \geq 0$.

Now we will prove $A_D(\phi,J_0;J') \neq 0$. We first derive the
following simple lemma.

\proclaim{Lemma 4.6} Suppose that $u: \R \times S^1 \to M$ is any
finite energy solution of
$$
\cases
\dudtau + J_t\Big(\dudt - X_H(u)\Big) = 0 \\
\int |{\part u \over \part \tau }|_{J_t}^2 < \infty.
\endcases
\tag 4.12
$$
that satisfies
$$
u(-\infty,t) = u(\infty, t). \tag 4.13
$$
Then $\int_{\R \times S^1} u^*\omega$ converges, and we have
$$
E_J(u) = \int_{\R \times S^1} u^*\omega. \tag 4.14
$$
\endproclaim

\demo{Proof} First note that when $H$ is nondegenerate,
any finite energy solution has well-defined asymptotic limits
$z^\pm = u(\pm\infty)$. Then we pick any bounding discs
$w^\pm$ of $z^\pm$ such that $w^+ \sim w^-\# u$. Now (4.14)
is an immediate consequence of (2.14) applied to
$\HH^{lin} \equiv H$, since we assume (4.13), i.e.,
$z^+ = z^-$.
\qed\enddemo

With this proposition, we are ready to prove positivity of
$A_D(\phi,J_0;J')$. Again a priori it is possible that
$A_D(\phi,J_0;J') = \infty$.

\proclaim{Proposition 4.7} Suppose that $\phi$ is nondegenerate,
and $J_0$ and $J' \in j_{(\phi,J_0)}$
as above. Then we have
$$
A_D(\phi,J_0;J') > 0.
$$
\endproclaim
\demo{Proof} We prove this by contradiction. Suppose
$A_D(\phi,J_0;J') = 0$ so that there exists a sequence of
non-constant maps $v_j: \R \times [0,1] \to M$ that satisfy (4.7)
and
$$
E_{J'}(v_j) \to 0 \quad \hbox{as $j \to \infty$}. \tag 4.15
$$
Therefore we have
$$
E_{J'}(v_j) < \e(J')
$$
for all sufficiently large $j$', where $\e(J')$ is the constant in
Lemma 4.2 and Corollary 4.3. In particular, the sequence $v_j$
cannot bubble off. This implies that $v_j$ locally uniformly
converge, and in turn that $v_j$ must (globally) uniformly
converge to a constant map because $E_{J'}(v_j) \to 0$. Since
there are only finitely many fixed points of $\phi$ by the
nondegeneracy hypothesis, by choosing a subsequence if necessary,
we conclude
$$
v_j(-\infty) = v_j(\infty) = p  \tag 4.16
$$
for all $j$'s for some $p \in \hbox{Fix }\phi$.
Now we fix any Hamiltonian $H: [0,1] \times M \to
\R$ that is zero near $t = 0 , \, 1$ and with $H \mapsto \phi$,
and consider the following maps
$$
u_j: \R \times S^1 \to M, \quad u_j(\tau,t): =
(\phi_H^t)(v_j(\tau,t)).
$$
It follows from (4.16) that
$$
u_j(-\infty,t) = u_j(\infty, t). \tag 4.17
$$
Furthermore for the family $J=\{J_t\}_{0 \leq t \leq 1}$
with
$$
J_t:= (\phi_H^t)_*(J'_t),
$$
the $u_j$'s satisfy the perturbed Cauchy-Riemann equation (4.12).

We note that (4.13) and the exponential convergence of $u_j(\tau)$
to $u_j(\pm\infty)$, as $\tau \to \pm\infty$ respectively, allows
us to compactify the maps $u_j$ and consider each of them as a
cycle defined over a torus $T^2$. Therefore  the integral $\int
u_j^*\omega$ depends only on the homology class of the
compactified cycles.

Now, because $v_j: \R \times [0,1] \to M$ uniformly converges to
the constant map $p \in \hbox{Fix }\phi$, the image of $u_j$ will
be contained in a tubular neighborhood of the closed orbit $z^p_H:
S^1 \to M$ of $\dot x = X_H(x)$ given by
$$
z_H^p(t) = \phi_H^t(p).
$$
In particular, $\int u_j^*\omega = 0$ because the cycle $[u_j]$ is
homologous to the one dimensional cycle $[z_H^p]$. Then Lemma 4.6
implies the energy $E_J(u_j) = 0$. But by the choice of $J$ above,
(2.31) implies $E_{J'}(v_j) = 0$, a contradiction to the
hypothesis that $v_j$ are non-constant. This finishes the proof of
Proposition 4.7. \qed\enddemo

We then define
$$
A(\phi,J_0;J') = \min\{A_S(\phi,J_0;J'),A_D(\phi,J_0;J')\}.
$$
Proposition 4.4 and 4.6 immediately imply
$$
A(\phi,J_0;J') > 0.
$$
However the finiteness $A(\phi,J_0;J') < \infty$ is a nontrivial
consequence of Theorem 5.4 below. Finally we define
$$
A(\phi,J_0) : = \sup_{J' \in j_{(\phi,J_0)}} A(\phi,J_0;J') \tag
4.18
$$
and
$$
A(\phi) = \sup_{J_0} A(\phi,J_0). \tag 4.19
$$
By definition, we have $A(\phi,J_0) > 0$ and so we have $A(\phi)
> 0$. However a priori it is not obvious whether they are finite,
which will be a consequence of the fundamental existence theorem,
Theorem 5.4, in the next section.

\head{\bf \S 5. Proof of nondegeneracy}
\endhead

With the definitions and preliminary studies of
the invariants $A(\phi,J_0;J')$ in section 4, this section is
occupied with the proof of nondegeneracy of the semi-norm
$$
\gamma: Ham(M,\omega) \to \R_+
$$
defined in section 3. As we pointed out in the introduction,
nondegeneracy follows from the following theorem.

\proclaim{Theorem 5.1} Suppose that $\phi$ is nondegenerate. Then
for any $J_0$ and $J' \in j_{(\phi,J_0)}$, we have
$$
\gamma(\phi) \geq A(\phi,J_0;J') \tag 5.1
$$
and hence
$$
\gamma(\phi) \geq A(\phi). \tag 5.2
$$
In particular, $A(\phi)$ is finite.
\endproclaim

The rest of this section is occupied by the proof of this
theorem.

Let $\phi$ be a nondegenerate Hamiltonian diffeomorphism with
$\phi \neq id$.
By the definition of $\gamma$, for any given $\delta > 0$, we can
find $H\mapsto \phi$ such that
$$
\rho(H;1) + \rho(\overline H;1) \leq \gamma(\phi)+ \delta. \tag
5.3
$$
For any Hamiltonian $H\mapsto $ know that $\overline H\mapsto
\phi^{-1}$. We will however use another Hamiltonian
$$
\widetilde H(t, x):= - H(1-t, x) \tag 5.4
$$
generating $\phi^{-1}$, which is more useful than $\overline H$,
at least in the study of duality and pants product.

\proclaim{Lemma 5.2} Let $H\mapsto \phi$.
Then $\widetilde H \mapsto \phi^{-1}$,
$[\phi^{-1},\overline H] = [\phi^{-1},\widetilde H]$.
In addition, if the (nondegenerate) spectrality axiom holds
for $a=1$, then we have
$$
\rho(\overline H;1) = \rho(\widetilde H;1). \tag 5.5
$$
\endproclaim
\demo{Proof} A direct calculation shows that the Hamiltonian
path generated by $\widetilde H$ is given by the path
$$
t \mapsto \phi_H^{(1-t)}\circ (\phi_H)^{-1}. \tag 5.6
 $$
The following composition of homotopies shows that this path is
homotopic to the path $t \mapsto (\phi_H^t)^{-1}$. Consider the
homotopy
$$
\phi_s^t:= \cases \phi_H^{s-t}\circ (\phi_H^s)^{-1}
& \quad \text{for }\, 0\leq t\leq s\\
(\phi_H^t)^{-1} & \quad \text{for }\, s\leq t\leq 1.
\endcases
$$
It is easy to check that the Hamiltonian path for $s=1$ is (5.6)
and the one for $s=0$ is $t \mapsto (\phi^t_H)^{-1}$ and satisfies
$\phi_s^1 \equiv (\phi_H)^{-1}$ for all $s\in [0,1]$. This finishes
the proof of $[\phi^{-1},\overline H] = [\phi^{-1},\widetilde H]$.
Now (5.5) follows from the homotopy invariance, Theorem 2.12.
\qed\enddemo

Lemma 5.2 now allows us to write $\gamma$ also as
$$
\gamma(\phi)= \inf_{H \mapsto \phi} (\rho(H;1) + \rho(\widetilde
H;1)).
$$
One advantage of using the representative $\widetilde H$ over
$\overline H$ is that the time reversal
$$
t \mapsto 1-t
$$
acting on the loops $z: S^1 \to M$ induces a natural one-one
correspondence between $\text{Crit}(H)$ and
$\text{Crit}(\widetilde H)$. Furthermore the space-time reversal
$$
(\tau,t) \mapsto (-\tau, 1-t) \tag 5.7
$$
acting on the maps $u: \R \times S^1 \to M$ induces a bijection
between the moduli spaces $\MM(H,J)$ and $\MM(\widetilde H,
\widetilde J)$ of the perturbed Cauchy-Riemann equations
corresponding to $(H,J)$ and $(\widetilde H, \widetilde J)$
respectively, where $\widetilde J_t= J_{1-t}$. This correspondence
reverses the flow and the corresponding actions satisfy
$$
\AA_{\widetilde H}([\widetilde z, \widetilde w]) = - \AA_H([z,w]).
\tag 5.8
$$
Here $[\widetilde z, \widetilde w]$ is the class corresponding to
$\widetilde z(t):= z(1-t)$ and $\widetilde w = w\circ c$ where
$c: D^2 \to D^2$ is the complex conjugation of $D^2 \subset \C$.

The following estimate of the action difference is an important
ingredient in our proof of nondegeneracy. The proof here is
similar to the analogous non-triviality proof for the Lagrangian
submanifolds studied in [\S 6-7, Oh4].

\proclaim{Proposition 5.3} Let $J_0$ be any  compatible almost
complex structure, $J' \in j_{(\phi,J_0)}$ and $J$ be the
one-periodic family $J_t = (\phi_H^t)_*J_t'$. Let $H$ be any
Hamiltonian with $H \mapsto \phi$. Consider the equation
$$
\cases
\dudtau + J_t\Big(\dudt - X_H(u)\Big) = 0 \\
u(-\infty) = [z_-,w_-], \, u(\infty) = [z_+,w_+] \\
w_- \# u \sim w_+, \quad u(0,0)=q \in B(\lambda)
\endcases
\tag 5.9
$$
for a map $u: \R \times S^1 \to M$. If (5.9) has a broken
trajectory solution (without sphere bubbles attached)
$$
u_1 \# u_2 \# \cdots \cdots \# u_N
$$
which is a connected union of solutions of (5.9) for $H$ that
satisfies the asymptotic condition
$$
\align
u_N(\infty) & = [z', w'], \, u_1(-\infty) = [z,w] \tag 5.10\\
u_j(0,0) & = q \,\,  \text{for some $j$}.
\endalign
$$
For some $[z,w] \in \text{Crit }\AA_H$ and $[\widetilde
z',\widetilde w'] \in \text{Crit }\AA_{\widetilde H}$, then we
have
$$
\AA_H(u(-\infty)) - \AA_H(u(\infty)) \geq A_D(\phi,J_0;J'). \tag
5.11
$$
\endproclaim
\demo{Proof} Suppose $u$ is such a solution. Opening up $u$ along
$t = 0$, we define a map $v: \R \times [0,1] \to M$ by
$$
v(\tau,t) = (\phi_H^t)^{-1}(u(\tau,t)).
$$
It is straightforward to check that $v$ satisfies (4.7) as
explained in subsection 2.4. Moreover we have
$$
\int \Big|{\part v \over \part \tau} \Big|^2_{J'_t} = \int
\Big|{\part u \over \part \tau} \Big|^2_{J_t} < \infty \tag 5.12
$$
from the energy identity (2.31).  Since $\phi(B(\lambda)) \cap
B(\lambda) = \emptyset$, we have
$$
v(\pm\infty) \in \hbox{Fix }\phi \subset M \setminus B(\lambda).
$$
On the other hand since $v(0,0) = u(0,0) \in B(\lambda)$, $v$
cannot be a constant map. In particular, we have
$$
\int \Big|{\part v \over \part \tau} \Big|_{J'_t}^2 = \int
v^*\omega \geq A_D(\phi, J_0;J'). \tag 5.13
$$
Combining (5.12) and (5.13), we have proven
$$
\AA_H(u(-\infty)) - \AA_H(u(\infty)) = \int \Big|{\part u \over
\part \tau} \Big|_{J_t}^2 \geq A_D(\phi, J_0;J').
$$
This finishes the proof. \qed\enddemo

This proposition demonstrates relevance of the existence result of
the equation (5.9) to Theorem 5.1. However we still need to
control the asymptotic condition (5.10) and establish some
relevance of the asymptotic condition to the inequality (5.1). For
this, we use (5.8) and impose the asymptotic condition
$$
u(-\infty) = [z,w], \, u(\infty) = [z',w']
$$
so that
$$
[z,w] \in \alpha_H, \quad [\widetilde z', \widetilde w'] \in
\beta_{\widetilde H} \tag 5.14
$$
for the {\it suitably chosen} fundamental Floer cycles $\alpha_H$ of
$H$ and $\beta_{\widetilde H}$ of $\widetilde H$.

We recall from (5.3) and (5.5) that we have
$$
\rho(H;1) + \rho(\widetilde H;1) \leq \gamma(\phi)+ \delta.
$$
By the definition (2.24) of $\rho$, there exist Floer cycles
$\alpha_H \in CF_n(H)$ and $\beta_{\widetilde H} \in
CF_n(\widetilde H)$ representing $1^\flat = [M]$ such that
$$
\aligned \rho(H;1) & \leq \lambda_H(\alpha_H) \leq \rho(H;1) + {\delta
\over 2} \\
\rho(\widetilde H;1) & \leq \lambda_{\widetilde
H}(\beta_{\widetilde H}) \leq \rho(\widetilde H;1) + {\delta \over
2}.
\endaligned
\tag 5.15
$$
Once we have these, by adding the two, we obtain
$$
\aligned
0 \leq \rho(H;1) + \rho(\widetilde H;1) & \leq
\lambda_H(\alpha_H) + \lambda_{\widetilde H}(\beta_{\widetilde H})\\
& \leq \rho(H;1) + \rho(\widetilde H;1) +  \delta.
\endaligned
\tag 5.16
$$
The fundamental cycles $\alpha_H$ and $\beta_{\widetilde H}$ that
satisfy (5.15)  will be used as the
asymptotic boundary condition required in (5.14).

The following is the fundamental existence theorem of the Floer
trajectory {\it with its asymptotic limits lying near the `top' of
the given Floer fundamental cycles} which will make the difference
$$
\AA_H([z,w]) - \AA_H([z', w'])= \AA_H([z,w]) + \AA_{\widetilde H}
([\widetilde z', \widetilde w'])
$$
as close to $\rho(H;1) + \rho(\widetilde H;1)$ as possible.
We refer readers to section 8 for its proof.

\proclaim{Theorem 5.4 [Fundamental Existence Theorem]}
Let $\phi$, $H$ and $(J_0;J')$ be as in
Proposition 5.3 and let $q \in M \setminus \text{\rm Fix}(\phi)$
be given. Choose $\delta$ so that
$$
0 < \delta < A_D(\phi,J_0;J'). \tag 5.17
$$
Then there exist some fundamental
cycles $\alpha_H$ of $(H,J)$ for $J_t = (\phi_H^t)_*J_t'$,
and $\beta_{\widetilde H}$ of $(\widetilde H, \widetilde J)$
such that they satisfy
$$
\aligned \lambda_H(\alpha_H) & \leq \rho(H;1) + \frac{\delta}{2} \\
\lambda_{\widetilde H}(\beta_{\widetilde H}) & \leq
\rho(\widetilde H;1) + \frac{\delta}{2}
\endaligned
\tag 5.18
$$
and have some generators
$[z,w] \in \alpha_H$ and $[\widetilde z', \widetilde w'] \in
\beta_{\widetilde H}$ that satisfy the following alternative :
\roster \item (5.9) has
a broken-trajectory solution (without sphere bubbles attached)
$$
u_1 \# u_2 \# \cdots \cdots \# u_N
$$
which is a connected union of Floer trajectories for $H$ that
satisfies the asymptotic condition
$$
u_N(\infty) = [z', w'], \, u_1(-\infty) = [z,w], \quad u_j(0,0) =
q \in B(\lambda)
$$
for some $1 \leq j \leq N$, (and hence
$$
\AA_H([z,w]) - \AA_H([z',w']) \geq A_D(\phi,J_0;J') \tag 5.19
$$
from Proposition 5.3) or,

\item there exists a $J_t'$-holomorphic sphere $v_\infty:S^2 \to
M$ for some $t \in [0,1]$ passing through the point $q \in
B(\lambda)$, and hence
$$
\AA_H([z,w]) - \AA_H([z', w']) \geq A_S(\phi,J_0;J') -\delta.
\tag 5.20
$$
\endroster
This in particular implies
$$
A(\phi,J_0;J') < \rho(H;1) + \rho(\widetilde H;1) + \delta<
\infty \tag 5.21
$$
for any $J'$ and $J_0$.
\endproclaim

\demo{Wrap-up of the proof of nondegeneracy} Let $\phi$ be a
nondegenerate Hamiltonian diffeomorphism. From the definition of
$\gamma(\phi)$, since $\delta$ and $H$ are arbitrary as long as
$H \mapsto \phi$, we immediately derive, from (5.21),
$$
A(\phi, J_0;J') \leq \gamma(\phi) \tag 5.22
$$
for all $J_0$ and $J' \in j_{(\phi,J_0)}$. Next by taking the
supremum of $A(\phi, J_0;J')$ over all $J_0$ and $J' \in
j_{(\phi,J_0)}$ in (5.22), we also derive
$$
A(\phi) \leq \gamma(\phi).
$$
This finishes the proof of Theorem 5.1 and so the proof of
nondegeneracy. \qed\enddemo

\head{\bf \S 6. The spectral norm, the ball area  and the
homological area of $\phi$}
\endhead

We first start with improving the lower bound (5.1). For this
purpose, we re-examine our set-up of the definition
$$
\gamma(\phi)= \inf_{H \mapsto \phi} (\rho(H;1) + \rho(\overline
H;1))= \inf_{H \mapsto \phi} (\rho(H;1) + \rho(\widetilde
H;1)).
$$
This involves only the class $1 \in QH^*(M)$ and the orbits
$[z,w]$ of $H$ and $[\widetilde z',\widetilde w']$ of $\widetilde
H$ of Conley-Zehnder indices
$$
\mu_H([z,w]) = n, \quad \mu_{\widetilde H}([\widetilde
z',\widetilde w']) = n. \tag 6.1
$$
\definition{Remark 6.1}
If we reverse the flow of $\widetilde H$ by considering the time
reversal map
$$
 [\widetilde z',\widetilde w'] \mapsto [z',w'] ; \,
\text{Crit }\AA_{\widetilde H} \to \text{Crit } \AA_H
$$
then the broken-trajectory solution of $u$ constructed in
Theorem 5.4 satisfies the asymptotic condition
$$
\aligned u(-\infty)  & = [z,w] \in CF_{n}(H) \quad
\text{i.e., $\mu_H([z, w]) =  n$} \\
u(\infty)  & = [z', w'] \in CF_{-n}(H) \quad\, \text{i.e.,
$\mu_H([z', w']) =- n$}
\endaligned
\tag 6.2
$$
with respect to our grading convention. We emphasize that the
time reversal map also replaces the
downward Novikov ring with the upward Novikov ring, which is a
source of complications in the non-aspherical case.
\enddefinition

Let $\phi$ be nondegenerate. We now introduce a stronger lower
bound of $\gamma(\phi$, denoted by $A(\phi;1)$.
This is a refined version of $A(\phi)$ which
correctly reflects the fact that the definition of $\gamma$
involves only the class $1 \in QH^*(M)$.

Let $J_0$ be a given almost complex structure and $J' \in
j_{(\phi,J_0)}$. Fix a point $q \in M \setminus \hbox{Fix }\phi$.
Let $\alpha_H$ and $\beta_{\widetilde H}$ be the Floer Novikov
cycles defined in section 5 with $[\alpha_H] =
[\beta_{\widetilde H}] = 1^\flat = [M]$. Consider all
broken-trajectory solutions
$$
u = u_1\# u_2 \# \cdots \# u_N
$$
where each $u_k$ itself is a cusp-curve which has the unique
principal component with a finite number of sphere bubbles.
Furthermore the chain of principal components satisfies
$$
\cases
\dudtau + J_t \Big(\dudt - X_H(u)\Big) = 0 \\
u_1(-\infty) = [z,w] \in \alpha_H, \, u_N(\infty) = [z', w'] \,
\text{with }\, [\widetilde z',
\widetilde w'] \in {\beta_{\widetilde H}}, \\
u_j (0,0) = q  \quad \text{for some $j=1, \cdots, N$
}
\endcases
\tag 6.3
$$
with $J_t = (\phi_H^t)_*J'_t$. Here one should regard $(0,0)$ as
the marked point in the domain of the broken-trajectory curves
as a one-marked `Floer-stable trajectory' in the standard
sense of [FOn], [Ru], [LT].

We define the energy of the broken-trajectory $u = u_1\# \cdots \# u_N$
by
$$
E_J(u):= \sum_{k=1}^N E_J(u_k)
$$
as usual. Here we define $E_J(u_k)$ by the sum of the energy of
the principal component and the symplectic areas of the sphere
bubbles. We would like to recall that there is no loss of
energy under the bubbling process and that the symplectic area of
such a $J_{t}$-holomorphic sphere bubble has the same value as the
harmonic energy with respect to the Riemannian metric associated
to $g_{J_t} = \omega(\cdot, J_t\cdot)$.

Now we introduce the main definition of this section. Let
$J' \in j_{(\phi,J_0)}$, $q \in M \setminus \text{Fix}\phi$
and $\alpha_H \in CF_n(H), \, \beta_{\widetilde H} \in
CF_n(\widetilde H)$ be as above. Opening up $u$ along $t = 0 \equiv 1$,
we consider
$$
v(\tau,t) = (\phi_H^t)^{-1}(u(\tau,t))
$$
on $\R \times [0,1]$ and extend to $\R \times \R$ by
$$
\phi(v(\tau, t+1)) = v(\tau,t)
$$
as before. We define
$$
\align A(\phi,J_0;J'; & \alpha_H, \beta_{\widetilde H}; q)
= \inf_{u}\Big\{E_J(u)
\mid u \, \text{
satisfies  and}\, (6.3), \, \, u(0,0) = q \Big\}\\
& = \inf_{v}\Big\{\int v^*\omega \mid v \, \text{ satisfies (1.14)
with} \\
& \hskip1.0in \exists [z^-,w^-] \in \alpha_H: z^-(t)
=\phi_H^t(v(-\infty)) \\
& \hskip1.0in \exists [\widetilde z^+,\widetilde w^+] \in
\beta_{\widetilde H} :
z^+(t)=\phi_H^t(v(+\infty)) \\
& \hskip1.0in \text{and } \, v(0,0) = q  \Big\}. \tag 6.4
\endalign
$$
As usual, we set $A(\phi,J_0;J'; \alpha_H, \beta_{\widetilde H}; q)
=\infty$, if (6.3) has no solution.
By definition, we have the upper bound
$$
A(\phi,J_0;J'; \alpha_H, \beta_{\widetilde H}; q) \leq
\lambda_H(\alpha_H) + \lambda_{\widetilde H}(\beta_{\widetilde H})
\tag 6.5
$$
if (6.3) has a solution. We then set
$$
\align
A(\phi, J_0;J';1;q) & =
\inf_{\alpha_H, \beta_{\widetilde H}}\Big\{
A(\phi,J_0;J';\alpha_H, \beta_{\widetilde H};q)\, \Big|\,
[\alpha_H] = 1^\flat =[\beta_{\widetilde H}]\Big\}.
\\
A(\phi, J_0;J';1) & = \sup_{q \in M\setminus \text{Fix }\phi}
A(\phi,J_0;J'; 1;q).
\endalign
$$
Theorem 5.4 implies
$$
A(\phi, J_0;J';1;q) \leq \rho(H;1) + \rho(\widetilde H;1)
= \gamma(H)
$$
for all $q \in M \setminus \text{Fix}\phi$, and hence
$$
A(\phi, J_0;J';1) \leq \gamma(H). \tag 6.6
$$
By definition, we also have
$$
A(\phi,J_0;J';1) \geq A(\phi,J_0;J').
$$
Finally we give the main definition.

\definition{Definition 6.2 (Homological area)} Let $\phi$, $J_0$ and
$J'$ as above. We define
$$
A(\phi,J_0;1) = \sup_{J'\in j_{(\phi,J_0)}} A(\phi,J_0;J';1),
$$
and
$$
A(\phi;1) = \sup_{J_0} A(\phi,J_0;1).
$$
We call $A(\phi;1)$ the {\it homological area} of $\phi$.
\enddefinition

With this definition, we have obtained the following theorem.

\proclaim{Theorem 6.3} Let $\phi$ be nondegenerate, and $J_0$ and
$J' \in j_{(\phi,J_0)}$ as before. Then we have
$$
A(\phi,J_0;J';1) \leq \gamma(\phi) \tag 6.7
$$
for all $J_0$ and $J'\in j_{(\phi,J_0)}$. In particular, we have
$$
0< A(\phi;1) \leq \gamma(\phi)< \infty. \tag 6.8
$$
\endproclaim
\demo{Proof} Let $H \mapsto \phi$. We take the supremum of the
right hand side of (6.6) over all $H \mapsto \phi$ which gives rise to
(6.7).
This  proves that $A(\phi,J_0;J';1)$ is finite.
We then take the infimum of (6.6) over all  $(J_0;J')$, which
gives rise to the proof of (6.8).
\qed\enddemo

Motivated by Theorem 6.3, we introduce the
following new capacities.

\definition{Definition 6.4} Let $B$ be any compact subset of $M$.
\roster \item {\bf (Spectral displacement energy)} We define
$$
e_\gamma(B) = \inf_{\phi}\{\gamma(\phi)\mid\phi(B) \cap B =
\emptyset\} \tag 6.9
$$
and call it the spectral displacement energy.

\item {\bf (Homological area capacity)} We define
$$
c_{hom}(B) : = \inf_{\phi}\Big\{ A(\phi;1)
\mid \phi(B) \cap B = \emptyset\Big\} \tag 6.10
$$
and call $c_{hom}(B)$ the
homological area capacity of $B$. As usual, we say $c_{hom}(B) =
\infty$ if $B$ cannot be displaced from itself by any Hamiltonian
diffeomorphism.
\endroster
\enddefinition

With these definitions, the following lower bound of the spectral
displacement energy $e_\gamma(B)$ immediately follows from Theorem 6.3.

\proclaim{Theorem 6.5} For any compact set $B \subset M$, we have
$$
e_\gamma(B) \geq c_{hom}(B). \tag 6.11
$$
\endproclaim

Now, to demonstrate that $A(\phi;1)$ is a much finer invariant
than $A(\phi)$, we provide an example of $\phi$ for which we
have the strict inequality
$$
0 < A(\phi) < A(\phi;1)
$$
on $S^2$, and a deformation of $\phi$ converging to a degenerate
Hamiltonian diffeomorphism such that $A(\phi) \to 0$ while
$A(\phi;1)$ is fixed.

\definition{Example 6.6} Consider the unit sphere $S^2$ with the
standard area form of area $4\pi$, and let $f$ be a Morse function such that
$-f$ has the unique maximum point $p^+$ and the minimum point
$p^-$ such that
$$
\max f - \min f = f(p^+) - f(p^-) < 2\pi, \quad \int_{S^2} f d\mu
= 0. \tag 6.12
$$
In addition, we make $f$ to have a pair consisting of local
maximum point $p_1$  and a saddle point $p_2$ such that there
exists a unique gradient trajectory connecting $p_2$ and $p_1$. We
can make the difference $(-f)(p_1) -(- f)(p_2) = f(p_2) - f(p_1)>
0 $ as small as we want. Now we consider the autonomous
Hamiltonian flow $\phi_f^t$ of $f$ and set $\phi = \phi_f^1$. Then
it is easy to see that the Morse complex of $-f$ is given by
$$
\part^{(-f)} p^+ = - p_2, \, \part^{(-f)} p_1 = p_2, \,
\part^{(-f)} p_2 = 0, \ \part^{(-f)} p^- = 0
\tag 6.13
$$
where $\part^{(-f)}$ is the boundary map in the Morse complex of $-f$.
Therefore we have
$$
\align
\ker \part^{(-f)} & = \hbox{span}_\Z \{ p^+ + p_1, p^-, p_2 \} \\
\hbox{im } \part^{(-f)} & = \hbox{span}_\Z \{ p_2\}.
\endalign
$$
Hence we have the Morse homology of $f$ given by
$$
HM_2(-f) = \{[p^+ + p^1]\}, \quad HM_0(-f) = \{[p^-]\}, \, \quad
HM_1(-f) = \{0\}.
$$
We claim
$$
\aligned
A(\phi;1) & = \max f - \min f  = f(p+) - f(p^-)\\
A(\phi) & = f(p_1) - f(p_2)
\endaligned
$$
and in particular we have $0 < A(\phi) < A(\phi;1)$.

By considering $f$ whose $C^2$-norm is sufficiently small, the
Floer complex
with respect to the standard complex structure $J_0=j$ is given by
$$
(CF(f), \part) \cong (CM(-f), \part^{(-f)}) \otimes \Lambda.
$$
It immediately follows that if we consider $J' = (\phi^t_f)^*J_0$,
then $J' \in j_{(\phi,J_0)}$ and the maps
$$
v_f(\tau,t) = \phi_f^t(\chi(\tau))
$$
are solutions of the mapping cylinder version of the
Cauchy-Riemann equation
$$
\cases {\part v \over \part \tau} + J'_t {\part v \over
\part t} = 0 \\
\phi(v(\tau,1)) = v(\tau,0), \quad \int_{\R \times [0,1]} |{\part
v \over \part \tau }|_{J'_t}^2 < \infty
\endcases
\tag 6.14
$$
for any gradient trajectory $\chi:\R \to M$ satisfying
$$
\dot\chi - \hbox{\rm grad } f(\chi) = 0.
$$
In particular the unique gradient trajectory connecting $p_1, \,
p_2$ gives rise to the unique Floer trajectory of (6.14) with the
same asymptotic condition. It is easy to check that this
trajectory is Fredholm regular. Therefore by the cobordism
argument, there exists at least one solution of (6.14) with the
same area as $v_f$ for all (generic) choices of $J_0$ and $J' \in
j_{(\phi,J_0)}$. This proves
$$
A(\phi,J_0;J') = f(p_1) - f(p_2)
$$
for all such $J_0, \, J'$ which in particular proves
$$
A(\phi) = f(p_1) - f(p_2). \tag 6.15
$$

Next since we have the fundamental class $1^\flat = [p^+ + p^1]$, it
follows that for each given point $p\in S^2 \setminus K$ where $K$
is the union of the trajectories connecting
$$
K = \overline W(p^+,p_2) \cup \overline W(p_1,p_2) \cup
\overline W(p_2,p^-)
$$
$p^+, \, p_1$ to $p_2$ and $p^2$ to $p^-$, there exists a unique
trajectory passing through $p$, issued at one of the local maximum
points $\{p^+, \, p_1\}$ to $p^-$ by the cap action in the Morse
homology of $-f$. (See [Oh4] for a proof of this fact.) Here we
denote by $W(p,q)$ the intersection
$$
W(p,q) = W^-(p) \cap W^+(q)
$$
and by $\overline{W}(p,q)$ its closure,
where $W^-(p), \, W^+(q)$ denote the unstable and stable manifolds
of $p$ and $q$ respectively.
Note that the closure $\overline{W}(p_2,p^-)$ divides $S^2$ into
two components
$$
S^2 \setminus \overline{W}(p_2,p^-) = \hbox{Int } W (p^+) \coprod
\hbox{Int } W(p_1).
$$
Now if we choose $p$ from $\hbox{Int }W(p^+)$, then the
corresponding $v_f$ has area given by $\max f - \min f = f(p^+) -
f(p^-)$.

It will be shown in section 10
that as long as the $C^2$ norm of $f$ is sufficiently small, {\it
the `thin' part of the Floer moduli space is separated from the
`thick' part thereof}. This then implies, again by the cobordism
argument, that there exists at least one Floer trajectory
connecting $p^+$ and $p^-$ for all generic choice of $J_0$ and $J'
\in j_{(\phi,J_0)}$ with area given by  $f(p^+) - f(p^-)$. This
then proves
$$
A(\phi;1) = f(p^+) - f(p^-). \tag 6.16
$$
Now (6.15) and (6.16) prove our claim.

Finally, by deforming $f$ so that it remains unchanged outside a
neighborhood of $\overline{W}(p_1,p_2)$ but letting the value $f(p_1) -
f(p_2) \to 0$ (see [Figure 5.2, Mi] for the description of such
a deformation of
Morse functions), we also see that $A(\phi)$ with $\phi=\phi_f$
can go to zero as $\phi$ converges to a degenerate Hamiltonian
diffeomorphism, which is not the identity, while $A(\phi;1)$ is
fixed.
\enddefinition
\medskip

The remaining section will be occupied by some discussion on
the relation between (6.11) and the optimal energy-capacity
inequality for the Hofer displacement energy $e(B) \geq c(B)$.
Recall that the Gromov capacity $c(B)$ is defined by
$$
\align
c(B)  & = \sup_{\lambda}\{ \lambda = \pi r^2 \mid g(B^{2n}(r))
\subset \hbox{Int }B, \\
& \quad g: B^{2n}(r) \to (M,\omega) \, \hbox{ is a symplectic
embedding}\}.
\endalign
$$
Another examination of the proof of Theorem 6.1
suggests the following general conjecture

\proclaim{[The Area Conjecture]} \roster \item Let $\phi$ be a
nondegenerate Hamiltonian diffeomorphism. Then we have
$$
A(\phi;1) \geq bArea(\phi). \tag 6.17
$$
\item The function $\phi \mapsto A(\phi;1)$ is continuous on
$Ham^{nd}(M,\omega)$ (in the $C^\infty$ topology), and can be
continuously extended to the set $Ham(M,\omega)$ as a
continuous function. In particular, the inequality (6.17) holds for
any arbitrary Hamiltonian diffeomorphism.
\endroster
\endproclaim

Combining Theorem 6.3 and the Area Conjecture, and using the
continuity of $\gamma$ and $bArea$ on $Ham(M,\omega)$ under the
$C^\infty$-topology, we now state several conjectures which are
immediate consequences of the Area Conjecture and Theorem 6.5.

\proclaim{Conjecture 6.7} Let $\phi$ be any Hamiltonian
diffeomorphism. Then we have
$$
\gamma(\phi) \geq bArea(\phi). \tag 6.18
$$
\endproclaim

\proclaim{Conjecture 6.8} Let $B\subset M$ be a closed subset.
Then we have $c_{hom}(B) \geq c(B)$.
\endproclaim

We note that we have the inequality
$$
e_\gamma(B) \leq e(B)
$$
where $e(B)$ is the Hofer displacement energy
$$
e(B) = \inf_{\phi}\{ \|\phi\|\mid \phi(B) \cap B = \emptyset \}.
$$
An immediate corollary of Theorem 6.5 and Conjecture 6.8
would be the optimal energy-capacity inequality for the
Hofer displacement energy $e(B) \geq c(B)$.
We hope to carry out further study of the homological area
$A(\phi;1)$ and the Area Conjecture elsewhere in the future.

\head{\bf \S 7.  The homological area and the generating function
of $\phi$}
\endhead

In this section, we will obtain some lower estimates, in the
spirit of [BP, Mc], for the Hofer norm of Hamiltonian
diffeomorphisms whose graphs in $(M,-\omega) \times (M,\omega)$
allow {\it generating functions} in the Darboux chart of the
diagonal $\Delta \subset M\times M$. As an application, we will
prove Theorem IV and its corollary stated in the introduction.

Suppose that the Hamiltonian diffeomorphism $\phi$ has the
property that its graph
$$
\Delta_\phi: = \text{graph } \phi \subset (M,-\omega) \times
(M,\omega)
$$
is close to the diagonal $\Delta \subset M\times M$ in the
following sense: Let $o_\Delta$ be the zero section of $T^*\Delta$
and
$$
\Phi: (\UU, \Delta) \subset M\times M \to (\VV, o_\Delta) \subset
T^*\Delta \tag 7.1
$$
be a Darboux chart such that \roster \item $\Phi^*\omega_0 =
-\omega \oplus \omega$ and

\item $\Phi|_\Delta = id_\Delta$ and $d\Phi|_{T\UU_\Delta}:
T\UU|_\Delta \to T\VV|_{o_\Delta}$ is the obvious symplectic
bundle map from $T(M\times M)_\Delta \cong N\Delta \oplus T\Delta$
to $T(T^*\Delta)_{o_\Delta} \cong T^*\Delta \oplus T\Delta$ which
is the identity on $T\Delta$ and the canonical map
$$
\widetilde \omega: N\Delta \to T^*\Delta; \quad X \mapsto X
\rfloor \omega
$$
on the normal bundle $N\Delta = T(M\times M)_\Delta/ T\Delta$.
\endroster

Consider the path
$$
t\mapsto \text{graph } t\,dS_\phi
$$
which is a Hamiltonian isotopy of the zero section and define
$\phi^t: M \to M$ to be the Hamiltonian  diffeomorphism satisfying
$$
\Phi(\Delta_{\phi^t}) = \text{graph } tdS_\phi \tag 7.2
$$
for $0 \leq t \leq 1$. By the requirement (1) and (2) for the
Darboux chart $\Phi$ of (7.1), it follows that the path $t \mapsto
\phi^t$ is a Hofer geodesic, i.e., a {\it quasi-autonomous
Hamiltonian path}: We recall a Hamiltonian path $t \mapsto \phi^t$
is called {\it quasi-autonomous} [BP], if its generating
Hamiltonian $H$ has a fixed maximum and a fixed minimum point over
$t \in [0,1]$. In addition, we note that the Hamiltonian path
$\phi^t$ provided by (7.2) {\it has no non-constant one periodic
orbits} at all. We denote by $H = H^\phi$ the corresponding
Hamiltonian with $\phi^t = \phi^t_{H^\phi}$.

We recall
$$
\rho(H;1) + \rho(\overline H;1) \leq \int_0^1 (\max H_t - \min
H_t) \, dt \tag 7.3
$$
On the other hand, if $H$ is the one obtained from (7.2), we also
have
$$
\int_0^t -\min H_u \, du = t \max S_\phi, \quad \int_0^t - \max
H_u \, du = t \min S_\phi. \tag 7.4
$$
and so
$$
\gamma(\phi) \leq \rho(H;1) + \rho(\overline H;1) \leq \max S_\phi
- \min S_\phi:= osc(S_\phi). \tag 7.5
$$

With this definition, we will prove the following theorem in
section 9, as a consequence of an existence theorem for solutions
of (1.14) satisfying the asymptotic conditions required in (6.4),
with the lower energy estimate.

All the results in the rest of this section can be generalized for
arbitrary {\it engulfable} Hamiltonian diffeomorphisms. For the
simplicity of exposition, we will restrict to the case where
$\phi$ is $C^1$-close to the identity, leaving the generalizations
to the engulfable case in [Oh10].

The following is the main theorem of this section. We will prove
the theorem in section 9 as a consequence of an existence theorem
of `thin' solutions of (1.14).

\proclaim{Theorem 7.1} Let $\phi$ be sufficiently $C^1$-small, and
$S_\phi$ be as above. Then we have
$$
A(\phi;1) \geq osc(S_\phi). \tag 7.6
$$
\endproclaim

Assuming this for the moment, we state several consequences of the
theorem. First we state the following corollary.

\proclaim{Corollary 7.2} Let $\phi$ be a $C^1$-small Hamiltonian
diffeomorphism and $S_\phi$ be its generating function as defined
above. Then we have
$$
\gamma(\phi) = osc(S_\phi) = A(\phi;1). \tag 7.7
$$
\endproclaim
\demo{Proof} Combining (6.8), (7.5) and (7.6), we derive the
following chain of inequalities
$$
A(\phi;1) \leq \gamma(\phi) \leq osc(S_\phi)\leq A(\phi;1).
$$
Therefore all the inequalities here are indeed the equalities.
This finishes the proof. \qed\enddemo

We now derive several consequences from this Corollary 7.2.

\proclaim{Theorem 7.3} Let $\Phi: \UU \to \VV$ be a Darboux chart
along the diagonal $\Delta \subset M \times M$, and $H=H^\phi$ is
the Hamiltonian generating $\phi^t: M \to M$ as in (7.2). Then
the path $t\in [0,1] \to \phi_H^t$ is length minimizing among all
paths from the identity to $\phi$. In this case, we also have
$\gamma(\phi) = \|\phi\|$.
\endproclaim
\demo{Proof} Let $K \mapsto \phi$ and $\phi^t_K$ be the
corresponding Hamiltonian path. By the canonical adjustment [Lemma
9.2, Oh5], we may assume that $K$ is one-periodic. It follows from
(3.2) that
$$
\gamma(K) \leq \|K\|. \tag 7.8
$$
Combining Corollary 7.2 and (7.4)-(7.5), we have obtained
$$
\|H^\phi\| = osc(S_\phi) = A(\phi,1) = \gamma(\phi).
$$
In particular, we have
$$
\gamma(\phi) = \|H^\phi\| \tag 7.9
$$
Therefore we derive  $\|H^\phi\| \leq \|K\|$ from (7.8) and (7.9)
which is exactly what we wanted to prove. The identity
$\gamma(\phi)=\|\phi\|$ follows from $\gamma(\phi) \leq \|\phi\|$,
(7.9) and the inequality $\|H^\phi\| \geq \|\phi\|$ by the
definition of $\|\phi\|$. This finishes the proof. \qed\enddemo

As in [BP, Mc], an immediate corollary of Theorem 7.3 is the
following $C^1$-flatness result of $Ham(M,\omega)$.

\proclaim{Corollary 7.4 [Proposition 1.8, Mc]} There is a path
connected neighborhood
$$
\NN \subset Ham(M,\omega)
$$
of the identity in the $C^1$-topology such that any element in
$\NN$ can be joined to the identity by a path that minimizes the
Hofer length. Moreover $(\NN,\|\cdot \|)$ is isometric to a
neighborhood of $\{0\}$ in a normed vector space, which is nothing
but the space $C^\infty_m(M)$ of normalized functions on $M$ with
the norm given by $osc(f) = \max f - \min f$.
\endproclaim

\bigskip

\centerline{\it PART II: ADIABATIC DEGENERATION AND}
\centerline{THICK-THIN DECOMPOSITION}
\medskip

\head{\bf \S 8. Proof of the fundamental existence theorem}
\endhead

In this section, we will give a complete self-contained proof of
Theorem 5.4, modulo some background materials presented in [En1],
[Oh9], restricting to the case when there is no discontinuity in
the Hamiltonian perturbation term and when the associated
Hamiltonian fibration is `flat'.

\subhead{\it 8.1. Pants product and Hamiltonian fibrations}
\endsubhead

We first recall the description of conformal structures on a
compact Riemann surface $\Sigma$ of genus 0 with three punctures
in terms of the minimal area metric [Z] which will be important
for the analysis involving the adiabatic degeneration of
pseudo-holomorphic curves relevant to the pants product.

We conformally identify $\Sigma$ with the union of three half
cylinders which we denote by $\Sigma_1$, $\Sigma_2$ and $\Sigma_3$
in the following way: the conformal structure on $\Sigma \setminus
\{x_1, x_2, x_3\}$ can be described in terms of the {\it minimal
area metric} [Z] which we denote by $g_\Sigma$. This metric makes
$\Sigma$ the union of three half-cylinders $\Sigma_i$ with flat
metric with each meridian circle having length $1$. The glued metric
on $\Sigma$ is smooth everywhere except at two points
$p, \overline p \in \Sigma$ lying on the boundary circles of $\Sigma_i$.
As in [FOh], we will use this metric for the
analytic estimates of pseudo-holomorphic curves implicit in the
argument. One important property of this {\it singular} metric is
that the metric is flat everywhere except at the two points $p, \overline
p$ where the metric has a {\it conical singularity}. Therefore
{\it its associated conformal structure smoothly extends across the
singularities $p, \, \overline p$.} (See [St] for such a
removable singularity theorem of the conformal structure.)
The resulting conformal
structure is the standard unique conformal structure on $\Sigma=
S^2 \setminus \{0, 1, \infty\}$.

As in [FOh], [En1], we identify $\Sigma$ as the union of $\Sigma_i$'s
$$
\Sigma = \cup_{i=1}^3 \Sigma_i
$$
in the following way: if we identify $\Sigma_i$ with $(-\infty, 0]
\times S^1$, then there are 3 paths $\theta_i$ of length ${1 \over
2}$ for $i = 1, \, 2,\, \, 3$ in $\Sigma$ connecting $p$ to
$\overline p$ such that
$$
\aligned
\part_1\Sigma & = \theta_1\circ \theta_3^{-1} \\
\part_2 \Sigma & = \theta_2 \circ \theta_1^{-1} \\
\part_3 \Sigma & = \theta_3 \circ \theta_2^{-1}.
\endaligned
\tag 8.1
$$
We fix a holomorphic identification of each $\Sigma_i, \, i= 1,2$
with $(-\infty, 0] \times S^1$ and $\Sigma_3$ with $[0,\infty)
\times S^1$ considering the decomposition (8.1).

We denote the identification by
$$
\varphi^+_i: \Sigma_i \to (-\infty, 0] \times S^1, \quad i =1, \,
2
$$
for positive punctures and
$$
\varphi^-_3: \Sigma_3 \to [0, \infty) \times S^1
$$
for the negative puncture. We denote by $(\tau,t)$ the standard
cylindrical coordinates on the cylinders.

We fix a cut-off function $\rho_-: (-\infty,0] \to [0,1]$ of the
type
$$
\rho_- = \cases 1  \quad & \text{for } \, \tau  \leq -r_{rig} \\
0 \quad  & \text{for }\, -r_{lef} \leq \tau \geq 0
\endcases
\tag 8.2
$$
with
$$
0< r_{lef} < r_{rig}, \tag 8.3
$$
and define $\rho_+: [0,\infty) \to [0,1]$ by $\rho_+(\tau) =
\rho_-(-\tau)$. We just denote by $\rho$ these cut-off functions when
there is no danger of confusion.

We now consider the (topologically) trivial {\it smooth} bundle
$P \to \Sigma$ with fiber isomorphic to $(M,\omega)$ and fix a trivialization
$$
\Phi_i: P_i:= P|_{\Sigma_i} \to \Sigma_i \times M
$$
on each $\Sigma_i$ for $i = 1,\, 2\, 3$. On each $P_i$,
we consider the closed two form of the type
$$
\omega_{P_i}:= \Phi_i^*(\omega + d(\rho H_t dt))
$$
for a time periodic Hamiltonian $H: [0,1] \times M \to \R$. Note
that this naturally extends to a closed two form $\omega_P$
defined on the whole $P \to \Sigma$ since the cut-off functions
vanish in the center region of $\Sigma$. This form is
nondegenerate in the fiber and restricts to $\Phi_i^*\omega$.

Such $\omega_P$ is called a {\it coupling form} and
induces a canonical symplectic connection $\nabla
= \nabla_{\omega_P}$ [GLS], [En1]. In addition it also fixes a
natural deformation class of symplectic forms on $P$ obtained by
those
$$
\Omega_{P,\lambda} := \omega_P + \lambda \omega_\Sigma
$$
where $\omega_\Sigma$ is an area form and $\lambda > 0$ is a
sufficiently large constant. We will always normalize
$\omega_\Sigma$ so that $\int_\Sigma \omega_\Sigma = 1$.

Next let $\widetilde J$ be an almost complex structure on $P$ such
that \roster \item $\widetilde J$ is $\omega_P$-compatible on each
fiber and so preserves the vertical tangent space and \item the
projection $\pi: P \to \Sigma$ is pseudo-holomorphic, i.e, $d\pi
\circ \widetilde J = j \circ d\pi$.
\endroster
When we are given three $t$-periodic Hamiltonians $H =
(H_1,H_2,H_3)$ and $J=(J_1,J_2,J_3)$, we say that $\widetilde J$
is {\it $(H,J)$-compatible}, if $\widetilde J$ additionally satisfies

(3) For each $i$, $(\Phi_i)_*\widetilde J = j\oplus J_{H_i}$ where
$$
J_{H_i}(\tau,t,x) = (\phi^t_{H_i})^*J_i
$$
on $M$ over a disc $D_i' \subset D_i$ in
terms of the cylindrical coordinates. Here $D_i' =
\varphi_i^{-1}((-\infty, -r_{i,rig}] \times S^1), i = 1,\, 2$ and
$\varphi^{-1}_3([r_{3,rig}, \infty) \times S^1)$ for some
$r_{i,rig} > 0$.

Condition (3) implies that the $\widetilde J$-holomorphic sections
$s$ over $\Sigma_i'$ for $i=1, \, 2\, \, 3$
are precisely the solutions of the equation
$$
{\part u \over \part \tau} + J_t^i\Big({\part u \over \part t} -
X_{H_i}(u)\Big) = 0 \tag 8.4
$$
if we write
$$
\Phi_i(s(\tau,t)) = (\tau,t, u_i(\tau,t))
$$
in the trivialization with respect to the cylindrical coordinates
$(\tau,t)$ on $D_i'$ induced by $\phi_i^\pm$ above.

Now we are ready to define the moduli space which will be relevant
to the definition of the pants product that we need to use. To
simplify the notation, we denote
$$
\widehat z =[z,w]
$$
in general and $\widehat z = (\widehat z_1, \widehat z_2, \widehat
z_3)$ where $\widehat z_i =[z_i,w_i] \in \text{Crit}\AA_{H_i}$ for
$i =1,2, 3$.

\definition{Definition 8.1} Consider the pair $(H,J)$ of
$H =\{H_i\}_{1\leq i \leq 3}$ and $J=\{J_i\}_{1 \leq i\leq 3}$, and let
$\widetilde J$ be a $(H,J)$-compatible almost complex structure.
Define $\MM(H, \widetilde J; \widehat z)$ to be the space of all
$\widetilde J$-holomorphic sections $s: \Sigma \to P$ that satisfy
\roster \item The maps $u_i: (-\infty, -r_{i,rig}] \times S^1 \to
M$ which are solutions of (8.4), satisfy
$$
\lim_{\tau \to -\infty}u_i(\tau,\cdot) = z_i, \quad i = 1,2
$$
where $\Phi_i(s(\tau,t)) = (\tau,t, u_i(\tau,t))$, and similarly for
$i=3$ changing $-\infty$ to $+\infty$.

\item The closed surface obtained by capping off $pr_M\circ\Phi_i
(s(\Sigma))$ with discs $w_i$ taken with the same orientation for
$i = 1,2$ and the opposite orientation for $i =3$ represents zero
in $\pi_2(M)$.
\endroster

We denote by $\widetilde\pi_2(M)$ the set of elements in
$\pi_2(M)$ modulo $\Z$-torsion elements.
\enddefinition

Note that $\MM(H, \widetilde J; \widehat z)$ depends only on the
equivalence class of $\widehat z$: we say that $\widehat z' \sim
\widehat z$ if
$$
z_i' = z_i, \quad w_i' = w_i \# A_i
$$
for $A_i \in \pi_2(M)$ and $\sum_{i=1}^3 A_i$ represents zero in
$\widetilde\pi_2(M)$. The dimension of $\MM(H, \widetilde J;
\widehat z)$ is given by
$$
\aligned \dim  \MM(H, \widetilde J; \widehat z) & = 2n -
(-\mu_{H_1}(z_1) + n) - (-\mu_{H_2}(z_2) +n) - (\mu_{H_3}(z_3) + n)\\
& = -n + (-\mu_{H_3}(z_3) + \mu_{H_1}(z_1) + \mu_{H_2}(z_2)).
\endaligned
\tag 8.5
$$
Note that when $\dim \MM(H,\widetilde J;\widehat z) = 0$, we have
$$
\mu_{H_3}(\widehat z_3) = \mu_{H_1}(\widehat
z_1)+\mu_{H_2}(\widehat z_2) -n.
$$
Now the pair-of-pants product $*$ for chains is defined by
$$
\widehat z_1 * \widehat z_2 =  \sum_{\widehat z_3} \#(\MM(H,
\widetilde J;\widehat z)) \widehat z_3 \tag 8.6
$$
for the generators $\widehat z_i$ and then by linearly extending
over the chains in $CF_*(H_1) \otimes CF_*(H_2)$. Our grading
convention makes this product have degree $-n$.

From now on, we will specialize our discussion to the triple
$(H_1, H_2, H_3)$ that allows a coupling form $\omega_P$ whose
associated symplectic connection becomes flat. For example,
the triple with $H_3 = H_1 \# H_2$ is one such choice.
For such a `flat' connection on $P \to \Sigma$, we have the identity
$$
\frac{1}{2} \int_\Sigma |(Ds)^v|_{\widetilde J}^2 = \int
s^*\omega_P. \tag 8.7
$$
for any $(H,J)$-compatible $\widetilde J$. {\it In general, such a
flat connection must have some pointwise singularities and it is
not possible to construct a {\it smooth} connection that is flat
everywhere, but one can choose a smooth connection whose curvature
becomes arbitrarily small both pointwise and in the integral
sense}. We refer readers to [section 7, Oh9] or to (10.9)-(10.10)
coming later for the $L^2$-estimates of the derivative
$|Ds|_{\widetilde J}$ in terms of the curvature of the general
Hamiltonian connection. We decompose $Ds = (Ds)^v + (Ds)^h$ the
sum of the vertical and horizontal parts with respect to the
chosen connection. Then we have the identity
$$
\int_\Sigma s^*\omega_P = \AA_{H_1}([z_1,w_1]) +
\AA_{H_2}([z_2,w_2]) - \AA_{H_1 \# H_2}([z_3,w_3]) \tag 8.8
$$
for any $\widetilde J$-holomorphic sections $s \in
\MM(H,\widetilde J;\widehat z)$. We refer to [\S 5, En1] or [Sc]
for its proof. We note that this, combined with the integral bounds
on the curvature, provides the uniform (total) energy bound for the
pseudo-holomorphic sections $s \in \MM(H,\widetilde J;\widehat
z)$. (See (10.9).) Therefore we will mainly concern the vertical
part of the energy.

Then (8.7) and (8.8) give rise to
$$
\frac{1}{2} \int_\Sigma |(Ds)^v|_{\widetilde J}^2
= -\AA_{H_1\# H_2} ([z_3,w_3])
+ \AA_{H_1}([z_1,w_1]) + \AA_{H_2}([z_2,w_2]) \tag 8.9
$$
if we choose the above mentioned {\it singular} flat connection.
Otherwise, (8.9) will be only an approximate identity whose error
can be made arbitrarily small.  We will ignore this small error
for the rest of discussion, preferring to use the singular
connection that is flat everywhere except at two points, $p, \,
\overline p$. However we emphasize that {\it the bundle $P \to \Sigma$
itself is smooth}.

This general identity can be also used the other way around.
More precisely, the following lower bound
will be a crucial ingredient for the limiting arguments
for the pseudo-holomorphic sections
with some asymptotic conditions which are allowed to vary
inside given Floer homology classes. The proof of this lemma
is an immediate consequence of (8.9) and the energy estimate
[Corollary 4.4, Oh9] which is omitted.

\proclaim{Lemma 8.2} Let $\widehat z_i = [z_i,w_i]$ for
$i=1,\, 2,\, 3$. Suppose $\AA_{H_3}([z_3,w_3]) \geq c$ for some
constant $c$. Then we have the lower bound
$$
\AA_{H_1}([z_1,w_1]) + \AA_{H_2}([z_2,w_2]) \geq c -(E^+(H_1) +
E^+(H_2)) + E^-(H_3). \tag 8.10
$$
\endproclaim

We will apply the above discussion to the triple
$$
H_1 = H, \, H_2 = \widetilde{H}\#(\e f), \, H_3 =  \e f \tag
8.11
$$
for which a flat connection can be explicitly constructed.
We give the construction of such a connection in Appendix 1.

Now (8.8) becomes
$$
\frac{1}{2} \int_\Sigma |(Ds)^v|_{\widetilde J}^2 = -\AA_{\e f}
([q,\widehat q]) + \AA_{H}([z_1,w_1])
+ \AA_{\widetilde{H}\#(\e f)}([z_2,w_2]). \tag 8.12
$$

\subhead{\it 8.2. Heuristic discussion of the proof of Theorem 5.4}
\endsubhead
\smallskip

To motivate our proof of Theorem 5.4, we first give the following
heuristic discussion imitating the proof of a similar existence
theorem given by the author [section 7, Oh4] for the purpose
of proving nondegeneracy of the invariant defined by (1.12)
on the cotangent bundle. We refer to [Oh7] for more detailed
explanations.

Morally we would like to apply the pants
product to the case
$$
H_1 = H, \, H_2 = \widetilde H, \, H_3 = 0 \tag 8.13
$$
for a pseudo-holomorphic section of an appropriate Hamiltonian
bundle $P \to \Sigma$ with boundary data
$$
\align \Phi\circ s|_{\part_1\Sigma} & =[z,w] \in \alpha_H, \,
\Phi\circ s|_{\part_2 \Sigma} =[\widetilde z',\widetilde w']
\in \beta_{\widetilde H}\\
\Phi\circ s|_{\part_3 \Sigma} & = [q,\widehat q] \in [M] \tag 8.14
\endalign
$$
with
$$
\AA_H([z,w])  + \AA_{\widetilde H}([\widetilde z',\widetilde w'])
\geq C(H),
$$
where $C(H)$ is a constant depending only on $H$. (See [Corollary
4.4, Oh9] for the precise value of $C(H)$.)

Note that because $H_3 = 0$  in the monodromy condition (8.13) and
the outgoing end with monodromy $H_2 = \widetilde H$ is equivalent
to the incoming end with monodromy $H_2 = H$, we can fill-up the
hole $z_3 \in \Sigma$ and consider the cylinder with one outgoing
and one incoming end with the same monodromy $H_1=H$. In other
words our Hamiltonian bundle $P \to \Sigma$ becomes the mapping
cylinder
$$
E_\phi:= \R \times \R \times M/ (\tau, t+1, \phi(x)) \sim (\tau,
t, x)  \to \R \times S^1
$$
of $\phi$ with the canonical Hamiltonian connection $\nabla_H$
induced by $H$. (See [En1] for a detailed explanation of the
relation between Hamiltonians and connections.) If we make the
almost complex structure $\widetilde J$ the pushforward of the
complex structure
$$
\widetilde J(\tau,t,x) = j\oplus J'(t,x)
$$
on $\R \times \R \times M$ under the covering projection
$$
\pi: \R \times \R \times M \to E_\phi =(\R \times \R \times
M)/(\tau,t+1,\phi(x)) \sim (\tau,t,x)
$$
and the natural connection induced by the Hamiltonian $H: [0,1]
\times \R \to M$, which is flat, then $\widetilde
J$-pseudo-holomorphic sections of $E_\phi$ become solutions of the
perturbed Cauchy-Riemann equation for $H$. The condition
$\Phi\circ s|_{\part_3\Sigma} = [q, \widehat q]$ in (8.14) leads to
the condition $q \in \text{Im }u$. Now the identity $1*1=1$
would lead to such a solution which in turn would
prove $\gamma(\phi) > 0$ via Proposition 5.3.

However in reality, $H_3=0$ is a degenerate Hamiltonian and
cannot be directly used for the construction of the above
pants product moduli space, which is the reason we have to
consider (8.11) instead and take the limiting argument as $\e \to 0$.
This requires us to study the singular degeneration problem, the
so called {\it adiabatic degeneration}, of the type studied in
[FOh] but with nontrivial quantum contributions around.

\subhead{\it 8.3. Construction of solutions I: analysis of the
thick part}
\endsubhead
\smallskip

In this subsection, we start the proof of the main
existence result, Theorem 5.4.

Suppose that we are given any $J' \in j_{(\phi,J_0)}$, and one
fixed $J_3' \in j_{(\phi_{\e f}^1,J_0)}$. Let $(\phi,J_0)$ be as
before and $f$ be a generic Morse function. For the sake of
simplicity, we will also assume that $f$ has a unique minimum
point.

We consider the admissible triple
$$
H_1 = H, \, H_2 = \widetilde{H}\#(\e f ), \, H_3 = \e f
$$
for the nondegenerate $H$. We note that if $H$ is nondegenerate
and $\e$ is sufficiently small, $H \#(\e f)$ is also
nondegenerate and there is a natural one-one correspondence
between $\text{Per}(H)$ and $\text{Per}(H \#(\e f))$ and the
critical points $\text{Crit}\AA_H$ and
$\text{Crit}\AA_{\widetilde{H}\#(\e f)}$. We again denote by
$[z',w']$ the critical point of $\AA_{H \#(\e f)}$ corresponding
to $[\widetilde z',\widetilde w'] \in
\text{Crit}\AA_{\widetilde{H}\#(\e f)}$.

Then we consider the triple of the families of almost complex
structures $J = (J_1, J_2, J_3)$ on the three ends $\Sigma_i\cong
(-\infty, -r_{i,rig}] \times S^1$ under the trivialization
$$
\Phi_i: P|_{\Sigma_i} \to \Sigma_i \times M
$$
defined by
$$
\align J^1_t & = (\phi_H^t)_*J'_t,\\
J^2_t & = (\phi^t_{\widetilde{H}\#(\e f)})_*
(\phi_{-\e f}^{1-t})^*J'_{(1-t)} \quad \text{and}\\
J^3_t & = (\phi_{\e f}^t)_*J'_{3,t}
\endalign
$$
respectively. Here the precise forms of $J^i_t$ are not important
except the point that they are chosen to be one-periodic.
Note that the family $(\phi_{-\e f}^{1-t})^*J'_{(1-t)} \in
j_{(\phi_{\widetilde H \#(\e f)},J_0)}$ and $J^2$ is still a one-periodic
family. As $\e \to 0$, $J^2$ converges to the
time-reversal family of $t\mapsto (\phi^t_H)_*J'_t$
$$
t \mapsto (\phi^{(1-t)}_H)_*J'_{(1-t)}.
$$
We also choose $J'_{3,t}$ so that
$$
J^3_t \equiv J_0 \tag 8.15
$$
at the end of $\Sigma_3$ and  the path
converges to the constant path $J_0$ as $\e \to 0$.
This is possible since $\phi^1_{\e f}$ converges to the identity map.

With these choices made, we study the limit as $\e \to 0$ of the
moduli space
$$
\MM(H^\e, \widetilde J^\e;\widehat z^\e)
$$
with the asymptotic boundary condition (8.14).
The main difficulty in studying this limit problem  is that
the limit is a {\it singular} limit and it is
possible that the image of the pseudo-holomorphic curves can
collapse to an object of Hausdorff dimension one. Our first
non-trivial task is to prove that there is a suitable limiting
procedure for any sequence $s^\e \in \MM(H^\e, \widetilde
J^\e;\widehat z^\e)$ as $\e \to 0$. We will show that there exists
a sub-sequence of $s^\e \in \MM(H^\e, \widetilde J^\e;\widehat
z^\e)$ that converge to the union
$$
u \cup \chi
$$
in the Hausdorff topology, where $\chi$ is a negative (cusp)
gradient trajectory of $-f$ landing at the critical point $q \in
\text{Crit }f$ and $u$ is a solution of (5.9) but satisfying
$$
u(0,0) \in \text{Im } \chi
$$
instead of $u(0,0) = q$. In particular, we will prove that the
piece of Hausdorff dimension one is always the image of a
(negative) gradient trajectory of $-f$.

This proven, it is easy to see that the above mentioned heuristic
discussion cannot produce a solution required in Theorem 5.4 {\it
unless the negative gradient trajectory of $-f$ landing at $q$ is
trivial} for some reason. At this stage, the condition
$$
\mu_{-\e f}^{Morse}(q) = 2n
$$
will play an essential role and enable us to conclude that the
only gradient trajectories of $-f$ {\it landing at} $q$ are the
constant map $q$ and so we are in a good position to start with.

Let $P \to \Sigma$, $\omega_P$ and $\Omega_{P,\lambda}$ be those
as defined in subsection 7.1. We equip $P$ with an
$(H,J)$-compatible almost
complex structure $\widetilde J$ such that
$$
\widetilde J = j \oplus J'_i
$$
on each $\Sigma_i$ where
$$
J'_1 \in j_{(\phi,J_0)}, \quad J'_2=\widetilde J'_1 \in
j_{(\phi^{-1}, J_0)}, \quad J'_3 \in j_{(\phi_{\e f}^1,J_0)}.
$$
More explicitly we define $\widetilde J$ by
$$
\widetilde J(\tau,t,x)(\alpha,\beta,\xi) = (-\beta, \alpha,
(\phi_{H_i^{\rho(\tau)}}^t)^*J_{i,t}(\xi - \beta
X_{\rho(\tau)H_i}) + \alpha X_{\rho(\tau)H}) \tag 8.16
$$
on each $\Sigma_i$
for each $i = 1,\, 2,\, 3$. (See [Oh9] for detailed usage of
the choice (8.16) in the derivation of $L^2$-bounds of the
derivative $Ds$ for a $\widetilde J$-holomorphic section $s$ of $P$.)

The $\widetilde J$-holomorphic sections $s$ over $\Sigma_i$ are
precisely the solutions of the equation
$$
{\part u \over \part \tau} + J^i_t\Big({\part u \over
\part t} - X_{H_i}(u)\Big) = 0 \tag 8.17
$$
if we write $\Phi_i(s(\tau,t)) = (\tau,t, u_i(\tau,t))$ in the
trivialization with respect to the cylindrical coordinates
$(\tau,t)$ on $\Sigma_i'$ induced by $\varphi_i^\pm$ as before.

By the construction, we have the identity
$$
1^\flat_{\e f} = h_{id}((1\cdot 1)^\flat) = [\alpha
*\beta]
$$
in homology for any cycles $\alpha
\in CF_*(H), \, \beta \in CF_* (\widetilde{H}\#(\e f))$
with
$$
[\alpha] = 1^\flat, \, [\beta] = 1^\flat, \, [\gamma] = 1^\flat,
$$.
In the level of cycles, we have
$$
\alpha*\beta = \gamma + \part_{\e f} (\eta) \tag 8.18
$$
for some $\eta \in CF_*(\e f)$. We recall the following
fundamental {\it tightness} result of the classical Mores
fundamental cycle in the Floer complex $CF_*(\e f)$.

\proclaim{Proposition 8.3 [Lemma 8.8, Oh5]} Consider the
autonomous Hamiltonian $\e f$ and  the cohomology class $1 \in
QH^0(M) \cong HQ^{2n}(-\e f)$. Then there is a {\it unique}
Novikov cycle $\gamma$ of $\e f $ of the form
$$
\gamma = \sum_j c_j [x_j,\widehat x_j] \tag 8.19
$$
with $x_j \in \text{Crit }_{2n}(- \e f)$ that represents the class
$1^\flat = [M]$. Furthermore the cycle $\gamma$ is tight in the
sense of Definition 2.13: for any Novikov cycle $\beta$ of $\e f$
in the class $[\beta] = [M]$, i.e, $\beta$ homologous to $\gamma$,
we have
$$
\lambda_{\e f}(\beta) \geq \lambda_{\e f}(\gamma). \tag 8.20
$$
\endproclaim

It follows from Proposition 8.3 and (8.18) that
$$
\lambda_{\e f}(\alpha*\beta) \geq -
{\delta\over 3} \tag 8.21
$$
with $\delta$ independent of $\e$ by choosing $\e$ sufficiently
small for any cycles $\alpha \in CF_*(H)$ and
$\beta \in CF_* (\widetilde{H}\#(\e f))$.

Now we specialize the above discussion to particularly chosen
$\alpha$ and $\beta$. By the definition of $\rho(\cdot; 1)$, we can find
$\alpha_H \in CF_n(H) $ and $\beta_{\widetilde H}\in CF_n(\widetilde H)$
satisfying
$$
\aligned
\lambda_{H}(\alpha_H)
& \leq \rho(H;1) + {\delta \over 2} \\
\lambda_{\widetilde H}(\beta_{\widetilde H})
& \leq \rho(\widetilde H;1) + {\delta \over 2}.
\endaligned
$$
We now transfer the cycle $\beta_{\widetilde H}$ to
a cycle $\beta_{\widetilde{H}\#(\e f)}\in CF_n(\widetilde{H}\#(\e f))$
defined by
$$
\beta_{\widetilde{H}\#(\e f)}: = h_\e(\beta_{\widetilde H})
\tag 8.22
$$
where $h_\e: CF_n(\widetilde H) \to CF_n(\widetilde{H}\#(\e f))$
is the Floer chain map $h_{(\LL,j)}$ over the linear homotopy
$\LL$ between $\widetilde H$ and $\widetilde{H}\#(\e f)$ and
$j$ close to the constant homotopy $J$.
Then we will have
$$
\rho(\widetilde H;1) - \e \max f \leq
\lambda_{\widetilde H}(\beta_{H}) \leq \rho(\widetilde H;1)
+ {\delta \over 2} + \e (-\min f).
\tag 8.23
$$
Here the upper bound immediately follows from (2.17).
However to get the lower bound, we use the choice made in
(5.17) for $\delta$ to apply the Handle sliding lemma [Oh5],
more specifically Corollary 6.4 [Oh5] :
By the choice (5.17) of $\delta$, we have
$$
\frac{\delta}{2} + \min (- \e f) < A_D(\phi,J_0;J') ( = A_{(H,J)})
$$
if $\e$ is sufficiently small. Here $A_{(H,J)}$ is the
constant for the equation (4.12) defined similarly as
$A_D(\phi,J_0;J')$. {\it This precludes the level of $\beta_{\widetilde H}$
going down too much under the chain map $h_\e$}, because we know
$$
|\rho(\widetilde{H}\#(\e f);1) - \rho(\widetilde H;1)|
\leq \e |f|.
$$
We refer readers to [section6, Oh5] for
the details of this argument to prove such a lower bound.

Furthermore by the continuity of $\rho(\cdot;1)$, we also have
$$
\lim_{\e \to 0}\rho(\widetilde{H}\#(\e f);1)  = \rho(\widetilde
H;1).
$$
Now we need to compare the values of $\lambda_{H}(\alpha_H)$ ,
$\lambda_{\widetilde{H}\#(\e f)}(\beta_{\widetilde{H}\#(\e f)})$
and $\lambda_{\e f}(\alpha_H*\beta_{\widetilde{H}\#(\e f)})$ .

For each given $\e > 0$, let $[z_\e,w_\e] \in \alpha_H$  and
$[\widetilde z'_\e,\widetilde w'_\e] \in \beta_{\widetilde{H}\#(\e f)}$
and $[q,\widehat q\# A] \in \gamma + \part_{\e f}(\eta
)$ such that the moduli space $\MM(H^\e,\widetilde J^\e;\widehat
z^\e)$ is non-empty.  Such a triple exists because of the relation
(8.18) and by the definition of the pants product `$*$'.  Using
(8.21),  we may choose $[q,\widehat q\# A] \in \gamma + \part_{\e
f}(\eta)$ so that
$$
\AA_{\e f}([q,\widehat q\# A]) \geq - {\delta \over 2} \tag 8.24
$$
Furthermore since $[\alpha_H*\beta_{\widetilde{H}\#(\e f)}] =
1^\flat$ and the {\it unique} maximum point of $-\e f$ is
homologically essential, we may also assume that $q$ is the unique
maximum point and $A=0$. For each given $\e >0$, let $s^\e$ be any
$\widetilde J$-holomorphic section of $P \to \Sigma$ with
asymptotic boundary conditions
$$
\aligned s|_{\part_1\Sigma} & =[z_\e,w_\e] \in \alpha_H, \,
s|_{\part_2 \Sigma} =[\widetilde z'_\e,\widetilde w'_\e] \in
\beta_{\widetilde{H}\#(\e f)}\\
s|_{\part_3 \Sigma} & = [q,\widehat q] \in \gamma
\endaligned
\tag 8.25
$$
written in the given trivialization $\Phi$. Then it follows from
(8.10) of Lemma 8.2 that we have
$$
\AA_H([z_\e,w_\e]) + \AA_{\widetilde{H}\#(\e f)}([\widetilde z'_\e,
\widetilde w'_\e])
\geq -\frac{\delta}{2} - ( E^+(H) + E^+(\widetilde{H}\#(\e f))) +
E^-(\e f).
$$
By choosing $\e > 0$ sufficiently small, we may assume that
$$
\AA_H([z_\e,w_\e]) + \AA_{\widetilde{H}\#(\e f)}([\widetilde
z'_\e, \widetilde w'_\e]) \geq -(E^+(H) + E^+(\widetilde{(-\e
f)\#H})) - \delta':= C(H) \tag 8.26
$$
where $\delta'$ can be made arbitrarily small by letting $\e \to
0$. In particular, combining (8.23) and (8.26), we derive the
lower bounds
$$
\aligned
C(H) - \rho(H;1) & \leq \AA_H([z_\e,w_\e]) \\
C(H) - \rho(\widetilde{H}\#(\e f);1) & \leq
\AA_{\widetilde{H}\#(\e f)}([z_\e,w_\e])\\
\endaligned
\tag 8.27
$$
for all $\e > 0$ by adjusting $C(H)$ slightly.
Therefore for each given $\e > 0$, (8.27) implies that there are
only finitely many possible
asymptotic periodic orbits among the generators $[z_\e,w_\e] \in
\alpha_H$ and $[z'_\e,w'_\e] \in \beta_{\widetilde{H}\#(\e f)}$
respectively such that the corresponding moduli space
$$
\MM(H^\e, \widetilde J^\e;\widehat z^\e)
$$
becomes non-empty for the asymptotic condition (8.25). On the
other hand, since we assume that $H$ is generic and nondegenerate,
$(-\e f )\# H$ is nondegenerate for all sufficiently small $\e
>0$, there is a canonical one-one correspondence between
$\text{Crit}\AA_H$ and $\text{Crit}\AA_{H\#(\e f)}$. Therefore
as we let $\e \to 0$, we may assume that the asymptotic orbits
$[z^\e,w^\e]\in \text{Crit}\AA_H$ converges to $[z,w]
\in \text{Crit}\AA_H$ and in turn $[\widetilde z'_\e, \widetilde w'_\e]
\in \text{Crit}\AA_{\widetilde{H}\#(\e f)}$ converges to
$[\widetilde z',\widetilde w'] \in \text{Crit}\AA_{\widetilde H}$ in the
$C^\infty$-topology. It remains to show that
$$
[\widetilde z',\widetilde w'] \in \beta_{\widetilde H}. \tag 8.28
$$
We now recall the definition (8.22) of
the cycles $\beta_{\widetilde{H}\#(\e f)}$. Because a Novikov cycle
is generated by an infinite number of generators, we can {\it not} say in
general that the whole cycle $\beta_{\widetilde{H}\#(\e f)}$
converges to $\beta_{\widetilde H}$ because the convergence
of the generators may not be uniform. However it is easy to check
by the compactness and transversality argument that
the generators in $\beta_{\widetilde{H}\#(\e f)}$ above any given level
indeed converge to those of $\beta_{\widetilde H}$,
using the fact that the linear homotopy from $\widetilde H$ to
$\widetilde{H}\#(\e f)$ converges to the constant homotopy
$\widetilde H$ as $\e \to 0$ and considering the corresponding
continuity equation (2.12) for the chain map $h_\e$. We leave the details
to the readers. This proves (8.28).

\definition{Remark 8.4}
Here we like to note
that convergence of $[z^\e,w^\e]\in \text{Crit}\AA_H$ is an immediate
consequence of the action lower bound (8.27) and
the Novikov condition of the cycle $\alpha_H$, while
convergence of $[\widetilde z'_\e, \widetilde w'_\e]
\in \text{Crit}\AA_{H\#(\e f)}$ is a consequence of this and
the {\it fixed} asymptotic condition
$$
s|_{\part_3 \Sigma}  = [q,\widehat q] \in \gamma
$$
in the outgoing end through a standard compactness argument on the
moduli space $\MM(H^\e, \widetilde J^\e;\widehat z^\e)$. The latter
is more subtle than the former because the Hamiltonian
$(-\e f) \# H$ varies as $\e \to 0$. Here the energy bound for
the sequence $s^\e$, convergence of the other ends and the fact that
$[\widetilde z'_\e, \widetilde w'_\e] \in \text{Crit}\AA_{H \#(\e f)}$
occurs as the asymptotic condition of the sequence $s^\e$, enter
the convergence proof in a crucial way.
\enddefinition

Having this convergence statement made for the asymptotic orbits,
for the simplicity of notations we omit the subscript $\e$ from
$[z_\e,w_\e]$ and $[z'_\e, w'_\e]$ and denote them by $[z,w]$ and
$[\widetilde z',\widetilde w']$ from now on. We also note that
because the degree of the quantum cohomology class $1$ is zero, we
have
$$
\mu_H([z,w]) = \mu_{\widetilde{H}\#(\e f)}([\widetilde z',
\widetilde w']) = \mu_{\e
f}([q,\widehat q]) = n
$$
and so
$$
\mu^{Morse}_{-\e f}(q) = n + \mu_{\e f}([q,\widehat q]) = 2n.
$$
In particular, any gradient trajectory $\chi: (K, \infty) \to M$
of $-\e f$ satisfying
$$
\dot\chi - \e \, \text{grad }f(\chi) = 0, \quad \lim_{\tau \to
\infty}\chi(\tau) = q
$$
must be the constant map $\chi \equiv q$, which is exactly what we
wanted to have in our heuristic discussion in the beginning of
this section.

We recall the uniform energy bound in (8.9) for the vertical
energy of the section $s$. Since this bound is uniform for all
sufficiently small $\e$ and for any choice of the cut-off
functions $\rho_i$ with $r_{i,rig} > r_{i,lef}$ mentioned in
(8.2), we can now carry out the adiabatic convergence argument for
$\Sigma_3$. To carry out this adiabatic convergence argument, we
will conformally change the metric on the base $\Sigma$ of the
fibration. We note that the vertical energy of the section $s$ is
invariant under the conformal change of the base metric. We will
realize this conformal change by a conformal diffeomorphism
$$
\psi_\e : S \setminus \{(0,0)\} \to (\Sigma\setminus
\{x_1, x_2, x_3\}, g_\e)
$$
where $S=\R \times S^1$ is the standard flat cylinder and the
metric $g_\e$ is constructed in a way similar to the minimal area
metric, but we change the lengths of the $\theta_i$'s in (8.1) in
the following way :
$$
\aligned
\text{length }\theta_1 & = 1 - {\e \over 2} \\
\text{length }\theta_2 & = \text{length }\theta_3 = {\e \over 2}.
\endaligned
\tag 8.29
$$
It is easy to see that we can choose the conformal diffeomorphism
$\psi_\e$ so that it restricts to a quasi-isometry
$$
\psi_\e: S \setminus D(\delta(\e)) \to (\Sigma_1^\e \cup \Sigma_2^\e,
g_\e)
$$
and to a conformal diffeomorphism
$$
\psi_\e : D(\delta(\e))\setminus \{0\} \to (\Sigma_3, g_\e)
$$
where $D(\delta)$ is the disc around $(0,0)$ and $\delta(\e) \to
0$ as $\e \to 0$ : $\psi_\e$ is the inverse of the conformal
diffeomorphism
$$
\psi_\e^{-1}(r, \theta) = \delta(\e) e^{-2\pi (r + \sqrt{-1}\theta)} ;
\quad [0, \infty) \times S^1 \to D(\delta(\e)) \setminus \{0\}. \tag 8.30
$$
We choose any sequences $\e_j \to 0$ and a
$(H,J)$-compatible $\widetilde J^\e$. Then choose any sequence
$s_j \in \MM(H^{\e_j}, \widetilde J^\e;\widehat z^{\e_j})$ that satisfy
(8.25). Consider the compositions
$$
\widetilde s_j = s_j\circ \psi_{\e_j}; \, S \setminus \{(0,0)\}
\to P.
$$
Since $\psi_\e$ is conformal, we have
$$
\int |(D\widetilde s_j)^v|^2_{\widetilde J^\e} = \int |(D
s_j)^v|^2_{\widetilde J^\e}
\tag 8.31
$$
and hence the uniform energy bound
$$
\align {1 \over 2} \int |(D\widetilde s_j)^v|^2_{\widetilde J^\e}
& = -\AA_{\e f}([q,\widehat q]) + \AA_{\widetilde{H}\#(\e f)}
([\widetilde z',\widetilde w']) +\AA_H([z,w])  \\
& \leq (\rho(H;1) + {\delta \over 2}) + (\rho(\widetilde{H}\#(\e f))
+ {\delta \over 2}) + {\delta \over 2}  \tag 8.32
\endalign
$$
provided $\e$ is sufficiently small. Here we used (8.9) for the
identity, (8.23)-(8.24)  for the inequality.

After choosing a subsequence, there are two alternatives to consider :
\roster
\item there exists some constant $c>0$ such that
$$
 {1 \over 2} \int_{D^2(\delta(\e_j))\setminus
\{(0,0)\}} |(D\widetilde s_j)^v|^2_{\widetilde J^\e} > c. \tag 8.33
$$
\item
$$
\limsup_{j \to \infty} {1 \over 2} \int_{D^2(\delta(\e_j))\setminus
\{(0,0)\}} |(D\widetilde s_j)^v|^2_{\widetilde J^\e} = 0. \tag 8.34
$$
\endroster

We will treat the case (8.33) in the end of subsection 8.5
and prove that $\widetilde s_j$ bubbles-off as $j \to \infty$.
We first study the case (8.34) in the next subsection.

\subhead{\it 8.4. Construction of solutions II; analysis of the
thin parts}
\endsubhead
\smallskip

We go back to the original sequence $s_j$. We can regard the
equation of pseudo-holomorphic sections with the asymptotic
conditions (8.25) as the following equivalent system of perturbed
Cauchy-Riemann equations:
$$
\align & \overline \part_{\widetilde J} s_i|_{C_i} = 0 \, \text{
for } \,
i = 1,2, 3 \tag 8.35-i \\
& s_1|_{\theta_3} = s_3|_{\theta_3},\, s_2|_{\theta_2} =
\widetilde s_1|_{\theta_2}, \, s_3|_{\theta_1} = s_2|_{\theta_1}.
\tag 8.36
\endalign
$$
We will study (8.35-i) for each $i = 1,2,3$ separately with the
matching condition (8.36) in our mind.

We start with (8.35-3). It follows from (8.34) and the conformal
invariance of the vertical energy that we have the energy bound
$$
{1 \over 2} \int_{\Sigma_3} |(D s_j)^v|^2_{\widetilde J} \leq
\delta_j \to 0 \tag 8.37
$$
as $j \to \infty$. On $[r_{3,rig}, \infty) \times S^1$, (8.35-3)
is equivalent to
$$
\dudtau + J_t^\e \Big(\dudt - X_{\e f}\Big) = 0 \tag 8.38
$$
for the section $s(\tau,t) = (\tau,t, u(\tau,t))$ where $J_t^\e =
(\phi_{\e f}^t)_*J'_3$ with $J'_3 \in j_{(\phi^1_{\e f},J_0)}$. We
fix $0 < r_{i,lef} < r_{i,rig}$ for the following discussion until
the end of the discussion on $\Sigma_3$. Only at the end, we will
let $r_{3,rig}$ go to zero.

We recall that $J^\e_t \to J_0$ as $\e \to 0$. The energy bound
(8.37) and the $\e$-regularity theorem (see Lemma 4.2) in the
context of pseudo-holomorphic curves, which can be easily modified
for (8.38) using the fact that $|X_{\e f}| \to 0$ in the
$C^\infty$-topology as $\e \to 0$, imply the following uniform
$C^1$-estimate. Here we would like to emphasize that we implicitly
use the fact that any disc of {\it a fixed size} in the cylinder
$\R \times S^1= \R \times (\R/\Z)$ is conformally isomorphic to
the standard disc $D = D^2(1)$. {\it This kind of (uniform)
estimate does not exist for the rescaled maps $\widetilde s_j$ as
$\e \to 0$, for which the conformal structures of the domain
cylinders degenerate as $\e \to 0$}.

\proclaim{Lemma 8.5} Let $\delta_j>0$ as in (8.37). Then there
exists $j_0$ such that for any $j \geq j_0$ and $u$, we have
$$
|Du_j(\tau, t)| \leq C (\delta_j + \e_j)
$$
on $\Sigma_3 = \R_+ \times S^1$ with $C > 0$ independent of $j$.
In particular, we have
$$
\operatorname{length}(t\mapsto u_j(\tau,t)) \leq C (\delta_j
+\e_j). \tag 8.39
$$
for any $\tau \in [0, \infty)$.
\endproclaim

\demo{Proof} We define
$$
v_j(\tau,t) = (\phi_{\e f}^t)^{-1}(u_j(\tau,t)) \tag 8.40
$$
on $[r_{3,rig},\infty) \times [0,1]$ and extend to $[1/2, \infty)
\times \R$ by the formula
$$
\phi_{\e_j f}^1(v_j(\tau, t+1)) = v_j(\tau,t)
$$
as before. Then it follows that $v_j$ satisfy
$$
\dvdtau + J_3' \dvdt = 0
$$
and the energy bound
$$
\frac{1}{2}\int_{[r_{3,rig},\infty) \times
[0,1]}|Dv_j|^2_{J'_{3,t}}\, dt \leq \delta_j.
$$
from (8.37). Furthermore by the remark right above Lemma 2.14, we
obtain
$$
\frac{1}{2}\int_{[r_{3,rig},\infty) \times
[-1/2,3/2]}|Dv_j|^2_{J'_{3,t}}\, dt \leq 2\delta_j. \tag 8.41
$$
The above equation is not precisely the equation of the type dealt
in the $\e$-regularity theorem, Lemma 4.2, {\it in that here
$J_3'$ is $t$-dependent}. However an examination of the proof of
[Proposition 3.3, Oh1] shows that the proof goes through for the
current case where $J_3'$ varies depending on the domain variable
$(\tau,t)$ of $v$, and the same $\e$-regularity theorem holds for
equation (8.38). Once we have the $\e$-regularity theorem at our
disposal, we cover $[2r_{3,rig},\infty) \times [0,1]$ by discs of
a fixed radius, say $r_{3,rig}$, whose center lie at
$[2r_{3,rig},\infty) \times [0,1]$. Note that any such disc is
contained in $[r_{3,rig},\infty) \times [-1/2,3/2]$. We choose
$j_0$ so large that $2\delta_j < \e$ for all $j \geq j_0$ where
$\e$ is the constant appearing in the $\e$-regularity theorem,
Lemma 4.2. Then, by applying Corollary 4.3 to each of these discs,
we obtain the uniform $C^1$-estimate
$$
|Dv_j(\tau,t)| \leq C\delta_j \to 0\tag 8.42
$$
on $[2r_{3,rig},\infty) \times [0,1]$ for a uniform constant $C>0$
as $j \geq j_0$. Since $u_j(\tau,t)=\phi_{\e_j f}^t(v_j(\tau,t))$
by the definition of $v_j$, we have
$$
\align \frac{\part u_j}{\part \tau} & = T\phi_{\e_j
f}^t\Big(\frac{\part v_j}{\part \tau}\Big) \\
\frac{\part u_j}{\part t} & = T\phi_{\e_j f}^t\Big(\frac{\part
u_j}{\part \tau}\Big) + X_{\e_j f}(u_j(\tau,t)) \\
& = T\phi_{\e_j f}^t\Big(\frac{\part u_j}{\part \tau}\Big) + \e_j
X_{f}(u_j(\tau,t)). \endalign
$$
Therefore we obtain
$$
\align |Du_j|^2 & = \Big|\frac{\part u_j}{\part \tau}\Big|^2 +
\Big|\frac{\part u_j}{\part t}\Big|^2  \\
& = \Big|T\phi_{\e_j f}^t\Big(\frac{\part v_j}{\part
\tau}\Big)\Big|^2 + \Big|T\phi_{\e_j f}^t\Big(\frac{\part
v_j}{\part \tau}\Big) + \e_j X_{f}(u_j(\tau,t))\Big|^2 \\
& \leq 2|T\phi_{\e_j f}^t|^2 \Big(\Big|\frac{\part v_j}{\part
\tau}\Big|^2 + |\frac{\part v_j}{\part \tau}\Big|^2\Big) +
2\e_j^2 |X_{f}(u_j(\tau,t))|^2.
\endalign
$$
The lemma immediately follows from this and (8.42) by adjusting
the constant $C> 0$ appropriately. \qed\enddemo

We reparameterize $u_j$ to define
$$
\overline u_j(\tau,t) = u_j\Big({\tau \over \e_j}, {t\over
\e_j}\Big)
$$
on $[2\e_j r_{3,rig} , \infty) \times \R /\e_j\Z$.  Then
$\overline u_j$ satisfies
$$
{\part \overline u_j \over \part \tau}(\tau,t) + J_t^{\e_j} {\part
\overline u_j \over \part t}(\tau,t) - \operatorname{grad
}_{g_{J_t^{\e_j}}}f(\overline u_j)(\tau,t) = 0. \tag 8.43
$$
We denote $\operatorname{grad}$ for
$\operatorname{grad}_{g_{J_0}}$ from now on. We choose the `center
of mass' for each circle $t \mapsto u_j(\tau,t)$ which we denote
by $\chi_j(\tau)$. More precisely, $\chi_j(\tau)$ is defined by
the following standard lemma. (See [K] for its proof.)

\proclaim{Lemma 8.6} Suppose a circle $z: S^1 \to M$ has diameter
less than $\delta$ with $\delta= \delta(M)$ sufficiently small and
depending only on $M$. Then there exists a unique point $x_z$,
which we call the center of mass of $z$, such that
$$
z(t) = \exp_{x_z}\xi(t),\quad \int_{S^1}\xi(t)\, dt = 0.
$$
Furthermore $x_z$ depends smoothly on $z$ but does not depend on
its parametrization.
\endproclaim

If $j_0$ is sufficiently large, (8.39) implies that the diameter
of each of the circles $t \mapsto u_j(\tau,t)$  is less than
$\delta$ given in  Lemma 8.5 for any $j \geq j_0$ and $\tau \in
[2\e_jr_{3,rig},\infty)$. Therefore we can write $\overline
u_j(\tau,t)$ as
$$
\overline u_j(\tau, t) = \exp_{\chi_j(\tau)} (\overline
\xi_j(\tau,t))
$$
with $\overline\xi_j(\tau,t) \in T_{\chi_j(\tau)}M$ for any $\tau
\in [\e_j r_{3,rig}, \infty)$ for all $j \geq j_0$, $j_0$
sufficiently large. Furthermore the $C^1$ estimate in Lemma 8.5
proves
$$
|\overline\xi(\tau,t)| \leq C \delta_j.
$$

We consider the exponential map
$$
\exp: \UU \subset TM \to M; \quad \exp(x, \xi):= \exp_x(\xi)
$$
and denote
$$
D_1\exp(x,\xi): T_xM \to T_xM
$$
the (covariant) partial derivative with respect to $x$ and
$$
d_2\exp(x, \xi): T_xM \to T_xM
$$
the usual derivative $d_2\exp(x,\xi):= T_{\xi}\exp_x: T_xM \to
T_xM$. We recall the property
$$
D_1\exp(x, 0) = d_2\exp(x, 0) = id
$$
which is easy to check. We now compute
$$
\aligned {\part \overline u\over \part \tau} & =
d_2\exp(\chi(\tau),\overline \xi(\tau,t))\Big({D \overline \xi
\over \part \tau}(\tau,t)\Big) + D_1\exp(\chi(\tau),
\overline \xi(\tau,t))(\dot\chi(\tau)) \\
{\part \overline u \over \part t} & = d_2\exp(\chi(\tau),
\overline \xi(\tau,t))\Big(
{\part \overline \xi \over \part t}(\tau,t)\Big) \\
\operatorname{grad} & f(\overline u)
=\operatorname{grad}f(\exp_{\chi(\tau)} \overline \xi(\tau,t)).
\endaligned
$$
Substituting these into (8.43) and multiplying
$(d_2\exp(\chi_j(\tau), \overline \xi_j(\tau,t)))^{-1}$
to the resulting equation, we get
$$
\align {D \overline \xi_j \over \part \tau}(\tau,t) & +
(d_2\exp(\chi_j(\tau), \overline
\xi_j(\tau,t)))^{-1}(D_1\exp(\chi_j(\tau),\overline
\xi_j(\tau,t))(\dot\chi_j(\tau)))\\
& + (\exp_{\chi_j(\tau)})^*J_t^\e(\chi_j(\tau))\Big({\part
\overline \xi_j \over \part t}\Big)(\tau,t) -
(\exp_{\chi_j(\tau)})^*(\operatorname{grad}f)(\overline\xi_j(\tau,t))
= 0
\endalign
$$
on $T_{\chi_j(\tau)}M$. Using the center of mass condition in
Lemma 8.6 and integrating over $t\in \R/\e_j\Z$,  we obtain
$$
\aligned \int_0^{\e_j} (d_2\exp(\chi_j(\tau), & \overline
\xi_j(\tau,t)))^{-1}
(D_1\exp(\chi_j(\tau),\overline \xi_j(\tau,t))(\dot\chi_j(\tau)))\, dt \\
& - \int_0^1 (\exp_{\chi_j(\tau)})^*(\operatorname{grad}f)
(\overline \xi_j(\tau,t)) \, dt= 0. \endaligned \tag 8.44
$$
Here we used the identities
$$
\align  & \int_{S^1} {\part \overline \xi_j \over \part t}(\tau,
t)\, dt = 0 \quad  \tag 8.45\\
& \int_{S^1} \overline \xi_j(\tau,t)\, dt \equiv 0 \equiv
\int_{S^1} {D \overline \xi_j \over \part \tau}(\tau, t)\, dt \tag
8.46
\endalign
$$
where the second identity of (8.46) is a consequence of the first.
On the other hand, it follows that there exists sufficiently small
$\delta_j > 0$ such that
$$
\align (1 - C\delta_j)|\dot\chi(\tau)| & \leq
|d_2\exp(\chi_j(\tau), \overline
\xi(\tau,t)))^{-1}(D_1\exp(\chi(\tau),\overline
\xi(\tau,t))(\dot\chi(\tau)))|\\
& \leq (1 + C\delta_j)|\dot\chi(\tau)| \tag 8.47
\endalign
$$
and
$$
|(\exp_{\chi(\tau)})^*(\operatorname{grad}f)(\xi(\tau,t))| \leq
(1+C\delta_j)|\operatorname{grad}f(\xi(\tau,t))| \leq (1+
C\delta_j)\|\operatorname{grad }f\|_{C^0} \tag 8.48
$$
hold. These follow from (8.44) and the following standard
estimates on the exponential map
$$
\aligned
|d_2\exp(x,\xi)(u)| & \leq C |\xi||u| \\
|D_1\exp(x,\xi)(u)| & \leq C |\xi||u|
\endaligned
$$
for $\xi,\, u \in T_xM$ where $C$ is independent of $\xi$, as long
as $|\xi|$ is sufficiently small, say smaller than injectivity
radius of the metric on $M$. (see [K].)

Combining (8.44), (8.47) and (8.48), we have obtained
$$
|\dot\chi_j(\tau)| \leq C\|\operatorname{grad }f\|_{C^0}.
$$
Therefore $\chi_j$ is equi-continuous on any given fixed interval
$[\eta, R] \subset (0, \infty)$ for $\eta$ taken arbitrarily small
and $R$ taken arbitrarily large. Here we recall that $2 \e_j
r_{3,rig} \to 0$ as $j \to \infty$. Therefore there exists a
subsequence of $\chi_j:[\eta,R] \to M$ uniformly convergent to
some $\chi_\infty:[\eta,R] \to M$. Furthermore it easily follows
from (8.48) and smoothness of the exponential map, we also have
the estimates
$$
\aligned |d_2\exp(\chi_j(\tau), \overline
\xi_j(\tau,t)))^{-1} & (D_1\exp(\chi_j(\tau), \overline
\xi_j(\tau,t))(\dot\chi_j(\tau)) - \dot\chi_j(\tau)| \\
& \leq C |\overline \xi_j(\tau,t)||\dot\chi_j(\tau)|
\endaligned
\tag 8.49
$$
and
$$
|\exp_{\chi_j(\tau)}^*(\operatorname{grad}f)(\overline\xi_j(\tau,t))
- \operatorname{grad}f(\chi_j(\tau))| \leq C
|\overline\xi_j(\tau,t)|. \tag 8.50
$$
where the constant $C$ depends only on $M$.  Since $\max|\overline
\xi_j|_{C^0} \leq C \delta_j \to 0$ uniformly, the equation (8.44)
converges uniformly to
$$
\dot\chi_\infty - \operatorname{grad }f(\chi_\infty) = 0.
$$
Therefore $\chi_\infty$ is a gradient trajectory of $f$ defined
on $[\eta, R] \subset (0, \infty)$. Recalling
the $C^1$-estimate in Lemma 8.5 and
$$
\overline u_j(\tau,t) = \exp_{\chi_j(\tau)} (\overline
\xi(\tau,t))
$$
we have proven that $\overline u_j|_{[\eta,R]\times
(\R/{\e_j\Z)}}$ uniformly converges to $\chi_\infty$. By a
boot-strap argument by differentiating (8.44), this convergence
can be turned into a $C^\infty$-convergence.

Therefore letting $\eta \to 0$ and $R \to \infty$, by a standard
argument of local convergence as in [Fl2], we prove that the
sequence $u_j$ converges to a {\it connected} finite union of
gradient trajectories of $-f$
$$
\chi_0 \# \chi_1 \#\cdots \# \chi_N
$$
for some $N \in \Z_+$ where $\chi_0$ is defined on $[0,\infty)$
and all other $\chi_i$'s are on $\R$, and $\chi_N$ satisfies
$$
\chi_N(\infty) = q.
$$
However since $\mu^{Morse}_{-\e f}(q) = 2n$, i.e, $q$ is a (local)
maximum of $-f$, $\chi_N$ must be constant which in turn implies
that all $\chi_k$ must be the constant map $q$ for all $k = 1,
\cdots, N$. Going back to the original variable $(\tau,t) \in
[2r_{3,rig}, \infty) \times \R/\Z$ and recalling
$$
u_j(\tau,t) = \overline u_j(\e_j\tau, \e_j t),
$$
we have proven that $u_j|_{[2r_{3,rig}, \infty) \times S^1}$
converges to the constant map $q$ in the (fine)
$C^\infty$-topology.

Finally, we repeat the above whole process for a sequence of
$r_{3,rig} \to 0$ and conclude that
$u_j|_{\Sigma_3 \times S^1}$ converges to the constant map $q$
in the (fine) $C^\infty$-topology, by the standard argument of
taking a diagonal subsequence. From now on, we assume this.

\medskip

\subhead{\it 8.5. Construction of solutions III: wrap-up}
\endsubhead
\smallskip

We have established that $u_j|_{\Sigma_3^\e \times S^1}$ converges
to the constant map $q$ in the (fine) $C^\infty$-topology under the
condition (8.34).
In particular $s|_{D(\delta(\e_j)\setminus \{(0,0)\}}$ converges to
the constant point $((0,0),q) \in \Sigma \times M = P$
and its restriction $s|_{\part D(\delta(\e_j)}$. Now we
consider the restriction
$$
s_j: \Sigma \setminus D(\delta(\e_j)) \to P|_{\Sigma
\setminus D(\delta(\e_j))}
$$
of the section $s_j: \Sigma \to P$. Since $s_j|_{\part D(\delta(\e_j))}$
converges to $((0,0),q)$ in the $C^\infty$ topology,
we can smoothly extend $s_j|_{\part D(\delta(\e_j))}$ into the
interior of $D(\delta(\e_j))$ : Apply Lemma 8.6 and write
$\Phi(s_j(t)) = (t, u_j(t))$ on $t \in \part D(\delta(\e_j))$ as
$$
u_j(t) = \exp_q\xi_j(t), \quad \int_{S^1}\xi_j(t)\, dt = 0
$$
and then we extend this to a smooth map
$u'_j: D(\delta(\e_j)) \to M$ this by the formula
$$
u'_j(z) = \exp_q(r \xi_j(\theta))
$$
$(r,\theta)$ is the polar coordinates of $z \in D(\delta(\e_j))$,
i.e., $z = r + \sqrt{-1}\theta$. Now we consider a modified sequence
$s'_j: \Sigma \to P$ of the sections $\widetilde s_j$ defined by
$$
s'_j(z) = \cases s(z) \quad & z \in \Sigma \setminus D(\delta(\e_j))\\
u'_j(z) \quad & z \in D(\delta(\e_j)).
\endcases
$$
By flattening the boundary of $u'$ reparameterizing
the $r$ variable, we can make $s'_j$ a smooth section for
each $j$. It follows that $s'_j$ is $\widetilde J$-holomorphic
on $\Sigma \setminus D(\delta(\e_j))$ and $|D s'_j|_{D(\delta(\e_j))}
\to 0$ as $j \to \infty$ by construction.

As we mentioned before the lower
bound (8.27) implies, by the definition of the Novikov chains,
that there are only finitely many elements $[z,w] \in \alpha_H$
and $[\widetilde z',\widetilde w'] \in \beta_{\widetilde H}$ above
the lower bound. Therefore we can now apply Gromov-Floer type compactness
arguments on $\Sigma$ as $j \to \infty$ to the approximately
$\widetilde J$-holomorphic sections $s'_j$.

We note that by the choice of minimal area metric representation of
$\Sigma$ given in (8.29), the union $\Sigma_1^\e \cup \Sigma_2^\e$
converges to $\Sigma \setminus \{(0,0)\}$ as $\e \to 0$ in the metric sense.

When $s'_j$ bubbles-off,
it follows from the $(H,J)$-compatibility of $\widetilde J$
via the maximum principle that the bubble $u$ must be contained in
a fiber. Here we recall the relation $J_t = (\phi_J^t)_*J'_t$.
If we define
$$
v(z):= (\phi_H^t)^{-1}(u(z)),
$$
$v$ is a $J_t'$-holomorphic sphere satisfying
$$
\omega(v) = \omega(u) \geq A_S(\phi,J_0;J')
$$
where the inequality comes from the definition (1.10) of
$A_S(\phi,J_0;J')$. Therefore we have
$$
\limsup_{j \to \infty}{1 \over 2} \int |(D s_j)^v|^2_{\widetilde
J^\e} = \limsup_{j \to \infty}{1 \over 2} \int |(D\widetilde
s_j)^v|^2_{\widetilde J^\e} \geq A_S(\phi,J_0;J')
$$
and are in the second alternative of Theorem 5.4.

Therefore, by choosing a subsequence if necessary,
it remains to consider the case where
there exists a constant $C > 0$ such that
$$
|Ds'_j|_{\widetilde J} \leq C. \tag 8.51
$$
Here we  point out that by the definition of the data
$(P, \widetilde J^\e, H^\e)$ and by the choice of minimal area
metrics $g_\e$ as in (8.1), the bundle data $(P, \widetilde
J^\e)$ on $\R \times S^1 \setminus \{(0,0)\}$ smoothly converge to
$(P, \widetilde J)$ on compact subsets of $\R \times S^1 \setminus
\{(0,0)\}$. Note that the latter datum is smooth.

Using the $C^1$-bound (8.51) and this smooth convergence of the
data $(H^\e,\widetilde J^\e)$ to $(H,\widetilde J)$
we can bootstrap to extract a broken-trajectory limit
$s_\infty$ on
$\Sigma=\Sigma^- \cup \Sigma^+$ such that
$$
s_\infty = (s^-_1 \# s^-_2 \# \cdots \# s^-_N)\# (s^+_1 \# \cdots
\# s^+_L)
$$
where $s^-_i$ for $1 \leq i \leq N-1$ and $s^+_j$ for $2 \leq i
\leq L$ are  maps from $\R \times [0,1]$ to $M$ but $s^-_N$ a map
from $\Sigma^-=(-\infty, 0] \times S^1$ and $s^+_1$ from $
\Sigma^+= [0,\infty) \times S^1$.
Now the uniform $C^1$-bound (8.51) of $\widetilde s_j$ implies
that $s^-_N \# s^+_1$ becomes indeed smooth even across $\tau = 0$
by the elliptic regularity. Hence we have produced a smooth
broken-trajectory solution of the equation (8.52).

If we write $s|_{\Sigma_i}$ in terms of the mapping
cylinder of $\phi$, all of $v_i^\pm$ satisfy
$$
\cases {\part v \over \part \tau} + J'_t {\part v \over
\part t} = 0 \\
\phi(v(\tau,1)) = v(\tau,0), \quad \int |{\part v \over \part \tau
}|_{J'_t}^2 < \infty.
\endcases
\tag 8.52
$$
In particular, the join $s^-_N \# s^+_1$ also satisfies (8.52) on
$\R \times [0,1]$.

Now we consider $u_\infty: \R \times S^1 \to M$ defined by
$$
u_\infty(\tau,t) =  (\phi_H^t)(v_\infty(\tau,t)). \tag 8.53
$$
The topological hypothesis
$$
[u \cup (\bigcup_{i =1}^3 w_i)] = 0 \quad \text{in } \, \widetilde\pi_2(M)
$$
in the definition of the pants product (see Definition 8.1) gives
rise to the basic topological condition (2.11), after attaching all
the possible bubbles to $u_\infty$ that we removed in the
process. Furthermore
$$
v(0,0) = u(0,0) = q \in B(\lambda)
$$
while $v(\pm\infty) \in \text{Fix }\phi$ and hence $v$ cannot be
constant. This gives rise to the alternative (1) of Theorem 5.4.

Finally we consider the case (8.33). We fix a sufficiently small
$\delta_0 >0$ so that (A.1) in Appendix holds. In particular,
since $\widetilde s_j$ have the uniform energy bound (8.32),
there exists $C > 0$ such that
$$
{1 \over 2} \int_{D^2(\delta_0)\setminus
\{(0,0)\}} |(D\widetilde s_j)^v|^2_{\widetilde J^\e} < C
\tag 8.54
$$
for all $j$, where $D^2(\delta_0) \setminus \{(0,0)\}
\subset C \setminus \{(0,0)\}$. Under the conformal map
$$
\align
\psi_j: \Big[\frac{1}{2\pi}\ln (\delta(\e_j)/
\delta_0)),\infty\Big)\times S^1 & \to
D^2(\delta_0)\setminus \{(0,0)\} ; \\
\psi_j(\tau,t) & = \delta(\e_j)e^{-2\pi(\tau + it)}
\endalign
$$
$D^2(\delta(\e_j))\setminus \{(0,0)\}$
corresponds to $[0,\infty) \times S^1$. We denote
$$
\overline s_j = \widetilde s_j \circ \psi_j.
$$
If we write
$$
\Phi_{\e_j}(\overline s_j)(\tau,t)= (\tau,t,u_j(\tau,t))
$$
in the trivialization $\Phi_\e: P \to S \times M$ defined as before,
the equation $\overline\partial_{\overline J_j} \overline s_j = 0$
on $D^2(\delta(\e_j))\setminus \{(0,0)\}$ corresponds to
$$
\frac{\partial u_j}{\partial \tau}
+ J^3_{\e_j,t} \Big(\frac{\partial u_j}{\partial t}
- \e_j X_f \Big) = 0 \tag 8.55
$$
on $[0,\infty)\times S^1$ similarly to (8.38).
Here we recall from (8.15) that $J^3_{\e_j,t} \equiv J_0$
at the end of $[0,\infty)\times S^1$.
Due to the asymptotic condition $\lim_{\tau \to \infty}u(\tau) =
q$, a constant orbit, the solution $u_j$ closes up at
$\infty$ and so $\widetilde s_j: D^2(\delta_0) \setminus \{(0,0)\}
\to P$ extends continuously over $(0,0)$ for each $j$.
Since $q$ is a critical point of a Morse function
$f$, $|Du_j|$ decays exponentially with the decay rate $b \e_j$
for some uniform constant $b>0$ depending only on $f$.
Given this exponential decay,
a straightforward calculation proves H\"older continuity
of $\widetilde s_j$ at $(0,0)$ : there exist constants
$C = C(u_j) > 0$ and $\delta_1 = \delta_1(u_j)$ such that
$$
|D\widetilde u|(z) \leq C |z|^{-1+b \e_j} \tag 8.56
$$
for all $z \neq 0$ with $|z| \leq \delta_1$. Here $\widetilde u_j$
is the vertical component of the local representation
$\widetilde s_j (z) = (z, \widetilde u_j)$ in
the trivialization $\psi_j^*\Phi_{\e_j}$.

On the other hand, the energy concentration (8.33) implies
the derivative blow-up
$$
|(D\widetilde s_j)^v(z_j)|  \to \infty \tag 8.57
$$
for some sequence $z_j \in D^2(\delta(\e_j))\setminus \{(0,0)\}$
as $j \to \infty$.

Now we will imitate the argument from [Fl2], [Oh1]
to produce a non-constant bubble. More specifically, we refer
readers to the argument from p 129 - p 130 of [Oh1].
We fix $\delta_0>0$ as before
and consider a pair of discs $D' \subset D$ given by
$$
D=D^2(\delta_0), \quad D' = D^2\Big(\frac{\delta_0}{2}\Big).
$$
We then consider the constants
$$
\aligned
\eta_j:=\inf\{\eta > 0 \mid & \quad \text{there exists
$x \in D^2\Big(\frac{\delta_0}{2}\Big)$ such that } \\
& \|(D\widetilde s_j)^v\|_{p,B_\eta(x)} \geq \eta^{2/p-1} \}
\endaligned
\tag 8.58
$$
for a fixed constant $p > 2$.
Here $B_\eta(x)$ is the closed $\eta$-ball in $D^2(\delta_0)$ with
the center lying at $x$.

Now we claim that $\eta_j \to 0$, after choosing a subsequence
if necessary. Suppose to the contrary that
$\eta_j$ is bounded away from zero. Then
we can choose $\eta_0$ so that
$$
0 < \eta_0 < \min\Big\{\inf_j\{\eta_j\}, \frac{\delta_0}{2}\Big\}.
$$
and so cover $D^2(\frac{\delta_0}{2})$ by a finite number, which
is independent of $j$, of balls $B_{\eta_0}(x) \subset D^2(\delta_0)$
such that  $x \in D^2(\frac{\delta_0}{2})$ and
$\|(D\widetilde s_j)^v\|_{p,B_{\eta_0}(x)} \leq \eta_0^{2/p -1}$. This
would then imply the $L^p$-bound
$$
\|(D\widetilde s_j)^v\|_{p,D^2(\delta_0/2)} \leq C(p, \eta_0, \delta_0)
\tag 8.59
$$
Using the facts that $\partial_{\widetilde J^{\e_j}}\widetilde s_j = 0$
and $J_3^{\e_k} \to j \oplus J_0$ on $D^2(\delta_0) \setminus \{(0,0)\}$
(see (A.1) Appendix), (8.59) implies the uniform derivative bound
$$
|(D\widetilde s_j)^v|_{\infty, D^2(\delta_0/3)} \leq C
$$
which then contradicts to (8.57). This proves $\eta_j \to 0$
as $j \to \infty$.

Therefore there must be some $x_j \in D^2(\delta_0/2)$ such that
$$
\int_{B_{\eta_j}(x_j)} |(D\widetilde s_j)^v|^p \geq
\frac{1}{2}\eta_j^{2-p}, \quad \eta_j \to 0. \tag 8.60
$$
Without loss of any generality, we may assume that $x_j \to x_0 \in
D^2(\delta_0/2)$. Obviously we have
$$
r_j: = \eta_j^{-1}\text{dist}(x_j,\partial D^2(\delta_0)) \to \infty.
$$
Then for every $R > 0$ and $j$ large enough we have rescaled maps
$$
s_j': \C \supset B_R(0) \to P ; s_j'(w) = \widetilde s_j(\eta_j
(w - x_j))
$$
satisfying
$$
\align
\|(Ds_j')^v\|_{p,B_1(0)} & \geq \frac{1}{2} \\
\|(Ds_j')^v\|_{p,B_1(x)} & \leq 1 \quad \text{for all }\,
x\in B_{R-1}(0).
\endalign
$$
Furthermore $s_j'$ satisfies the equation
$\overline\partial_{J'_j} s_j' = 0$ on $B(R)\setminus \{0\}$
where $J'_j$ is the pull-back family of $\widetilde J^{\e_j}$
under the above rescaling map
$$
w \in B_R(0) \to \eta_j(w - x_j) \in B_{\eta_j}(x_j) \subset
D^2(\delta_0).
$$
In particular $J'_k \to j\oplus J_0$ on $B_R(0) \times M$
in the limit as $k \to \infty$. Now letting $R \to \infty$ and taking
a diagonal subsequence, we have obtained a non-constant
$J_0$-holomorphic map $v_\infty': \C \setminus \{0\} \to M$
with energy less than $C$. Applying the removal singularity theorem,
we then obtain a nonconstant $J_0$-holomorphic bubble.

Combining all the above, we have
finally finished the proof of Theorem 5.4.

\definition{Remark 8.7} Here we have developed our adiabatic convergence
scheme only up to the level that enables us to prove that the adiabatic sequence
converges to a gradient trajectory of the Morse function
$f$, {\it after renomalization}, when there
occurs no non-constant bubble. We hope
to develop a full-scale convergence scheme of the
adiabatic sequence considered in this paper elsewhere in the future.
\enddefinition

\head{\bf \S 9. Local Floer complex and its homology}
\endhead

In this section, we will prove Theorem 7.1. We first note that
proving the inequality $A(\phi;1) \geq \hbox{\rm osc}(S_\phi)$
requires to prove an existence of solution $v$ of (1.12)
satisfying
$$
E_{J'}(v) \geq \hbox{\rm osc}(S_\phi).
$$
for some $J_0$ and for some $J' \in j_{(\phi;J_0)}$. For this
purpose, we will exploit some distinguished features of the Floer
boundary operator and the pants product, when the relevant
Hamiltonians $H = H^\phi$ as defined in (7.1), (7.2) and Theorem
7.1 are $C^1$-small. This distinguished feature was previously
exploited by the present author in the context of Lagrangian
intersection Floer theory [Oh2], which we call {\it the thick and
thin decomposition} of the Floer moduli space.

We first recall the
following standard invariant of $(M,\omega)$
$$
\align A(\omega; J_0) & = \inf\{\omega([u])\mid u \,
\text{ is non-constant $J_0$-holomorphic} \} \\
A(\omega) & = \sup_{J_0} A(\omega; J_0).
\endalign
$$

\subhead {\it 9.1.  The thin part of the Floer boundary operator}
\endsubhead
\smallskip

The following proposition will be proved in section 10.

\proclaim{Proposition 9.1} Let $J_0$ be an almost complex
structure on $M$ and $H$ have no non-constant one-periodic orbit.
Let $\alpha_i, \, i=1, \, 2$ be any two positive constants
satisfying
$$
0 < \alpha_i < A(\omega; J_0), \quad \alpha_1 + \alpha_2 <
A(\omega; J_0)
$$
and any $\alpha_3 > 0$ be given. Let $J_t' = (\phi_H^t)^*J_0$ as
before. Let $u$ be any cusp-solution of (5.9). Then there exists a
constant $\delta = \delta(J_0,\alpha)> 0$ such that if $\|H\| <
\delta$, the following alternative holds: \roster \item either $u$
is `very thin'
$$
\int \Big|{\part u \over \part \tau} \Big|^2_{J_0} = \int
\Big|{\part v \over \part \tau} \Big|^2_{J'_t } = \int v^*\omega <
\alpha_1, \tag 9.1
$$
\item or it is `thick'
$$
\int \Big|{\part u \over \part \tau} \Big|^2_{J_0}= \int
\Big|{\part v \over \part \tau} \Big|^2_{J'_t}
> A(\omega; J_0) - \alpha_2.
\tag 9.2
$$
\endroster
In addition,  by choosing $\delta > 0$ smaller if necessary, we
can achieve
$$
\sup_{\tau \in \R}\hbox{ \rm diam }(t \mapsto u(t, \tau)) <
\alpha_3 \tag 9.3
$$
for all thin solutions, which in turn implies
$$
\int \Big|{\part u \over \part \tau} \Big|^2_{J_0} \leq \|H\| =
\int_0^1 (\max_x H_t - \min_x H_t)\, dt \tag 9.4
$$
\endproclaim
Note that {\it any $C^2$-small Hamiltonian or any $C^1$-small
autonomous Hamiltonian functions satisfy the hypotheses of this
proposition}. In the following discussion, we will assume one of
these two conditions.

Proposition 9.1 enables us to decompose the Floer boundary
operator, with respect to the Floer regular pair $(H,J)$,
$$
\part = \part_0 + \part'
$$
where $\part_0$ is the contribution from those satisfying (9.1)
and $\part'$ the one from those satisfying (9.2). We fix
$$
\alpha_2 < \frac{1}{4} A(\omega; J_0).
$$
As in [Oh2], by
comparing the areas of the terms in the identity
$$
\part \part = 0
$$
it follows that $\part_0$ satisfies the identity
$$
\part_0 \part_0 = 0. \tag 9.5
$$
We call the complex $(CF(H), \part_0)$ with $\part_0 =
\part_{(J_0,H),0}$ the {\it little Floer complex}
of the pair $(J_0,H)$ and the corresponding homology,
denoted by
$$
HF^{little}(H) := H_*(CF(H), \part_0)
$$
the {\it little Floer homology} of $(H,J_0)$ or simply of $H$.

A similar decomposition result also holds for the continuation
equation (2.12) as long as the homotopy $\HH = \{H^s\}_{0 \leq
s \leq 1}$ is $C^1$-small. This gives rise to a decomposition of
the chain map
$$
h_\HH = h_{\HH, 0} + h_\HH'.
$$
Again the area comparison of the thick and thin parts of the chain
map identity $\part_{H^1} \circ h_\HH = h_\HH \circ
\part_{H^0}$,
we have
$$
\part_{H^1,0} \circ h_{\HH,0} = h_{\HH,0} \circ
\part_{H^0,0}. \tag 9.6
$$

Next we recall that if $\|H\|_{C^1}$ is sufficiently small, all
the one-periodic orbits of $H$ are constant orbits corresponding
to the fixed point set of the time-one map $\phi_H^1$. They carry
the obvious constant bounding discs. This provides a distinguished
subset
$$
\text{Crit}^{fix}\AA_H: = \{[p, \widehat p] \mid p \in
\text{Fix}\phi_H^1 \} \subset \text{Crit}\AA_H \tag 9.7
$$
of $\text{Crit}\AA_H$. We consider the $\Q$-subspace
$$
CF^{fix}(H): = \text{span}_\Q \Big\{\text{Crit}^{fix}\AA_H
\Big\}\subset CF(H). \tag 9.8
$$
It follows from Proposition 9.1 that the subspace $CF^{fix}(H)$ is
preserved under the action of $\part_0$ and hence
$$
(CF^{fix}(H), \part_0)
$$
forms a sub-complex of the little Floer complex $(CF(H),\part_0)$
and satisfies
$$
(CF(H),\part_0) = (CF^{fix}(H), \part_0)\otimes_\Q \Lambda_\omega
\tag 9.9
$$
as a chain complex. We call $(CF^{fix}(H), \part_0)$ the {\it
local Floer complex} of $(H,J_0)$ or simply of $H$, and its
associate homology the {\it local Floer homology} which we denote
by $HF^{loc}(H;\Q)$. By a similar discussion to the above, under
the $C^1$-small homotopy $\HH$ between two  $H$'s, $h_{\HH,0}$
induces a chain isomorphism between the two local Floer complexes.

We summarize the above discussions into the following theorem.

\proclaim{Theorem 9.2} Let $\delta > 0$ be a sufficiently small
constant and $H_i$ for $i =0, \, 1$ satisfy $\|H_i\| < \delta$.
Consider the linear homotopy
$$
\HH: s \mapsto (1-s)H^0 + s H^1
$$
and let $h_\HH = h_{\HH, 0} + h_\HH'$ be the above decomposition.
Then both maps
$$
\align h_{\HH,0} & : (CF(H^0), \part_0) \to (CF(H^1), \part_0) \\
h_{\HH,0} & : (CF^{fix}(H^0), \part_0) \to (CF^{fix}(H^1),
\part_0)
\endalign
$$
are chain isomorphisms and so induce isomorphisms in their
homologies
$$
\align h_{\HH,0} & : HF^{little}(H^0) \to HF^{little}(H^1),\\
h_{\HH,0} & : HF^{loc}(H^0;\Q) \to HF^{loc}(H^1;\Q)
\endalign
$$
by choosing $\delta > 0$ smaller, if necessary. In particular, we
have the canonical isomorphisms
$$
\align  HF^{loc}(H;\Q) & \cong H_*(M;\Q)  \tag 9.10 \\
HF^{little}(H) & \cong HF^{loc}(H;\Q) \otimes_\Q \Lambda_\omega
\cong H(M;\Q) \otimes_\Q \Lambda_\omega.
\endalign
$$
\endproclaim
We refer to [Ch] or [Oh2,5] for some similar results and for more
details of the proof of the theorem.

\subhead{\it 9.2. The thin part of the pants product}
\endsubhead
\smallskip

Next we consider the pants product. We will be particularly
interested in the case
$$
H_1 = H^\phi, \, H_2 = \widetilde H^\phi,  \, H_3 = 0
$$
which we regard as the limit of the case (8.11).

In this case, again the relevant moduli space is decomposed into
the thin part  and the thick part, but this time we need to
incorporate the setting of the Hamiltonian fibration. A similar
decomposition exists for arbitrary numbers of ends, of course.

To motivate our discussion, we recall the quantum product of $a,
\, b \in H^*(M)$ can be written as
$$
a \circ b = a\cap b + \sum_{A \neq 0} (a,b)_A q^{-A}
$$
where $(a,b)_A \in H^*(M)$ is the cohomology class defined by
$$
(a,b)_A = PD ([\MM_3(J;A) \times_{(ev_1,ev_2)}(Q_a\times Q_b),
ev_0]).
$$
Here $\MM_3(J;A)$ is the set of stable maps with three marked
points $(x_1, x_2, x_3) \subset S^2$ and
in class $A \in H_2(M)$. We denote by $ev_i:\MM_3(J;A) \to
M$ for $i =0, \, 1, \, 2$ the evaluation maps. Then
$$
\MM_3(J;A) \times_{(ev_1,ev_2)}(Q_a\times Q_b)
$$
is the fiber product of $\MM_3(J;A)$ and $Q_a, \, Q_b$ via the
evaluation maps $(ev_1,ev_2)$. $[\MM_3^A
\times_{(ev_1,ev_2)}(Q_a\times Q_b), ev_0]$ is the homology class
of the fiber product as a chain via the map
$$
ev_0:\MM_3(J;A) \times_{(ev_1,ev_2)}(Q_a\times Q_b) \to M.
$$
Geometrically it is provided by the image (by $ev_0$) of the
holomorphic spheres intersecting the cycles $Q_a$ and $Q_b$ at the
first and second marked points. The term corresponding to $A=0$
provides the classical cup produce $a\cup b$.

We recall the pants product $\gamma * \delta$ of the Novikov
cycles $\gamma$ and $\delta$ of $H_1$ and $H_2$ is defined by
replacing the above $\MM_3(J;A)$ by the moduli space
$$
\MM(H,\widetilde J;\widehat z)
$$
for $\widehat z=(\widehat z_1, \widehat z_2, \widehat z_3)$ with
$\widehat z_i = [z_i,w_i]$. We note that  when all the
one-periodic orbits of $H_1, \, H_2,\, H_3$ are constant, any
periodic orbit $z= p$ has the canonical lifting $[p,\widehat p]$
in $\widetilde \Omega_0(M)$. Therefore we can write the generator
$[z,w]$ as
$$
p \otimes q^A:=[p, \widehat p\#A].
$$
We fix a trivialization $P = \Sigma \times M$ and write $u =
pr_2\circ s$.

For each given $p=(p_1,p_2,p_3)$ with $p_i$ (constant) periodic
orbits of $H_i$ respectively, we denote by $\MM_3(H,\widetilde J;
p)$ the set of all $\widetilde J$-holomorphic sections $s$ over
$\Sigma$ with the obvious asymptotic conditions
$$
u(x_1) = p_1, u(x_2) = p_2, \, u(x_3) = p_3 \tag 9.11
$$
where $u$ is the fiber component of $s$ under the trivialization
$\Phi: P \to \Sigma \times M$ and the $x_i$'s are the given
punctures in $\Sigma$ as described in section 8.1. If we denote by
$\MM_3(H,\widetilde J)$ the set of all {\it finite energy}
$\widetilde J$-holomorphic sections, then it has the natural
decomposition
$$
\MM_3(H,\widetilde J) = \cup_{p} \MM_3(H,\widetilde J;\widehat p).
$$
Furthermore, because of the asymptotic condition (9.7) for any
element $s \in \MM_3(H,\widetilde J)$, $u = pr_M\circ s$ naturally
compactifies and so defines a homotopy class $[u] \in \pi_2(M)$.
In this way, $\MM_3(H,\widetilde J;\widehat p)$ is further
decomposed into
$$
\MM_3(H,\widetilde J;\widehat p) = \cup_{A \in
H_2(M)}\MM_3^A(H,\widetilde J;\widehat p)
$$
where $\MM_3^A(H,\widetilde J;\widehat p)$ is the subset of
$\MM_3(H,\widetilde J;\widehat p)$ in class $A$. We denote
$$
\MM'_3(H,\widetilde J;\widehat p)=\cup_{A\neq
0}\MM_3^A(H,\widetilde J;\widehat p) \tag 9.12
$$
and then have
$$
\MM_3(H,\widetilde J;\widehat p) = \MM_3^0(H,\widetilde J;\widehat
p) \coprod \MM'_3(H,\widetilde J;\widehat p).
$$
In the following proposition, the quantities for the triple $H
=(H_1,H_2;H_3)$ defined by
$$
\aligned E^-(H) & : = E^-(H_1) + E^-(H_2) + E^+(H_1 \# H_2)\\
 E^+(H) & : = E^+(H_1) + E^+(H_2) + E^-(H_1 \# H_2)\\
 E(H) &: = E^+(H) + E^-(H)
 \endaligned
\tag 9.13
$$
play a similar role for the pants products as (1.2) does for the
boundary map. More elaborate discussion of these quantities will
be carried out in [Oh10]. Related study of these kind of
quantities is also given in [En1,2].

\proclaim{Proposition 9.3} Suppose that $H = (H_1, H_2, H_3)$
consisting of $H_i$ that have no non-constant periodic orbits
and let $J$ be given. Fix
any constants $\alpha_i, \, i=1, \, 2$ with
$$
0 < \alpha_i < A(\omega; J_0), \quad \alpha_1 + \alpha_2 <
A(\omega; J_0). \tag 9.14
$$
Consider the Hamiltonian fibration $P \to \Sigma$ of $H$ equipped
with $(H,J)$-compatible almost complex structure $\widetilde J$ with
respect to the symplectic form
$$
\Omega_{P,\lambda} = \omega_P + \lambda \omega_\Sigma.
$$
We denote by $s$ an element in  $\MM_3(H,\widetilde J)$. Then
there exists a constant $\delta = \delta(J_0,\alpha)$ such that
for any $H=(H_1,H_2,H_1 \# H_2)$ satisfying
$$
E(H), \quad \max_{i}\|H_i\|_{C^1} < \delta,
$$
we can find a coupling form $\omega_P$  a
sufficiently small constant $\lambda>0$ for which the following
alternative holds: \roster \item all those $s$ with $\int
s^*\omega_P = 0$ are `very thin'
$$
\int |Ds|^2_{\widetilde J}  < \alpha_1 \tag 9.15
$$
and its  fiber class is $[u] = 0$ or, \item all the elements $s$
with $\int s^*\omega_P \neq 0$ are `thick' i.e.,
$$
\int |Ds|^2_{\widetilde J} > A(\omega; J_0) - \alpha_2. \tag 9.16
$$
\endroster
\endproclaim

For the trivial generators $\widehat p_i = [p_i,\widehat p_i], \,
i= 1, \, 2$, the pants product is given by
$$
\widehat p_1 * \widehat p_2 = \widehat p_1 *_0 \widehat p_2 +
\sum_{A\neq 0} \widehat p_1 *_A \widehat p_2 \tag 9.17
$$
where $\widehat p_1 *_A \widehat p_2$ is defined as
$$
\widehat p_1 *_A \widehat p_2  = \sum_{p_3 \in Per(H_3)}
(p_1,p_2;p_3)_A \, \widehat p_3\otimes q^A.
$$
Here
$$
(p_1,p_2;p_3)_A = \#(\MM_3^A(H,\widetilde J;\widehat p)). \tag
9.18
$$
(9.17) induces the formula for arbitrary generators $[z,w] = [p,
\widehat p\#A]$ in an obvious way.

The main point of Proposition 9.3 is that the thin part
$\MM_3^0(H, \widetilde J)$ of the  moduli space
$\MM_3(H,\widetilde J)$ is `far' from the thick part $\MM_3'(H,
\widetilde J)$, and isolated under the continuation of the whole
moduli space $\MM_3(H,\widetilde J)$ via $C^1$-small deformations
of Hamiltonians as long as the thick solutions $s$ described in
Proposition 9.3 do not appear during the deformations. This is
guaranteed if we require (9.11) as the asymptotic boundary
condition, and if $H_3 = H_1 \# H_2$ and so the corresponding
curvature of the fibration $P$ is `flat'. (See section 8 and the
proof of Proposition 9.3 in section 10 for more discussion
relevant to this remark.)

This in particular implies that the term in the pants product
induced by $\widehat p_1*_0 \widehat p_2$, which is the
contribution of the thin part of the moduli space, is invariant in
homology under the $C^2$-small continuation of $H$ and induces the
classical cup product
$$
\cup: H^k(M) \times H^\ell(M) \to H^{k+\ell}(M).
$$
In particular the product by the identity $1$ induces an
isomorphism $1 :H^0(M) \to H^0(M)$.

A more systematic discussions on the local complex, the pants
product and their applications will be carried out in [Oh10].

\subhead{\it 9.3. The canonical Floer fundamental cycle}
\endsubhead
\smallskip

In this subsection, we recall  the construction and some
distinguished properties of {\it canonical Floer fundamental
cycles} from [Oh7]. We choose a Morse function $f$ such that $f$
has the unique global minimum point $x^-$ and
$$
f(x^-) = 0, \quad f(x^-) < f(x_j) \tag 9.19
$$
for all other critical points $x_j$. Then we choose a fundamental
Morse cycle
$$
\alpha =\alpha_{\e f} = [x^-,\widehat x^-] + \sum_j a_j
[x_j,\widehat x_j] \tag 9.20
$$
where $x_j \in \text{Crit }_{2n} (-f)$.

Considering Floer's homotopy map $h_\LL$ over the linear path
$$
\LL: \, s\mapsto (1-s) \e f +s H
$$
for sufficiently small $\e > 0$, we transfer the above fundamental
Morse cycle $\alpha$ and define a fundamental Floer cycle of $H$
by
$$
\alpha_H: = h_{\LL}(\alpha) \in CF(H). \tag 9.21
$$
This is  what we call {\it the canonical fundamental Floer cycle}
of $H$.  We would like to emphasize that {\it the linear homotopy
$\LL$ used for this fundamental Floer cycle can be made
$C^2$-small if $H$ is so}, and so $h_\LL$ restricts to a
well-defined chain isomorphism between the Morse complex of $-\e
f$ and the local Floer complex $(CF(H);\part_0)$. We recall the
following lemma from [Oh7]. We like to mention that
any sufficiently $C^1$-small autonomous $H$ will satisfy the
hypotheses in this lemma.

\proclaim{Lemma 9.4  [Proposition 4.2, Oh7]} Suppose that $H$ is
a generic one-periodic such that $H_t$ has the unique
nondegenerate global minimum $x^-$ which is fixed
and under-twisted for all $t \in [0,1]$.
Suppose that $f: M \to \R$ is
a Morse function such that $f$ has the unique global minimum point
$x^-$ and $f(x^-)=0$. Then the canonical fundamental cycle has the
expression
$$
\alpha_H = [x^-, \widehat x^-] + \beta \in CF(H)\tag 9.22
$$
for some Floer Novikov chain $\beta \in CF(H)$ with the inequality
$$
\lambda_H(\beta) < \lambda_H([x^-,\widehat x^-]) = \int_0^1 - H(t,
x^-)\, dt. \tag 9.23
$$
In particular its level satisfies
$$
\align \lambda_H(\alpha_H) & = \lambda_H([x^-,\widehat x^-])\\
& = \int_0^1 - H(t, x^-)\, dt = \int_0^1 -\min H\, dt.
\endalign
$$
\endproclaim
The same construction applies to $\beta_{\widetilde H}$ by replacing
the maximum and minimum points in the above construction and considering the
homotopy
$$
s \mapsto (1-s)  (\e f) + s \widetilde H.
$$
The following lemma shows that the above canonical fundamental
cycles $\alpha_H$ and $\beta_{\widetilde H}$ are tight in the
sense of Definition 2.13 {\it inside the local complex}. The proof
of this was given in [Proposition 4.5, Oh5] which is stated in
terms of {\it homological essentialness} of the critical point
$x^\pm$. (See [Definition 13.2, Po1].) The two statements are
however equivalent.

\proclaim{Lemma 9.5 [Proposition 4.5, Oh5]} Suppose that $H$ is
nondegenerate and quasi-autonomous with the unique fixed maximum
point $x^+$ and a fixed minimum point $x^-$. Suppose that
$\|H\|_{C^2} < \delta$ with $\delta >0$ is sufficiently small.
Then both $\alpha_H$ (respectively $\beta_{\widetilde H}$) are
tight in the local complex $(CF^{fix}(H), \part_{H,0})$
(respectively in $(CF^{fix}(\widetilde H), \part_{\widetilde
H,0}$); that is, for any $\part_{H,0}$-cycle $\alpha$ satisfying
$\alpha = \alpha_H + \part_{H,0}(\delta)$ for a chain $\delta \in
CF^{fix}(H)$, then
$$
\lambda_H(\alpha) \geq \lambda_H(\alpha_H) = \int_0^1 -H(x^-)\,
dt. \tag 9.24
$$
\endproclaim

\medskip
\subhead{\it 9.4. Proof of Theorem 7.1}
\endsubhead
\smallskip

In this subsection, we will use the isomorphism (9.10) of the
local Floer complex to prove an existence theorem of `thin'
solutions of (1.12).

Let $J_0$ be a given almost complex structure and $H$ be given.
Define $J'$ to be the family defined by $J'_t = (\phi_{H}^t)^*J_0$
as before. We first note that
$$
A_S(\phi_H^1,J_0;J') = A(\omega;J_0). \tag 9.25
$$
This is because there is a one-one correspondence
$$
w \mapsto \phi\circ w
$$
between the set of $J_0$-holomorphic spheres and that of
$\phi^*J_0$-holomorphic spheres.

We go back to the study of $H^\phi$. Now for given $J_0$, let us
choose $\phi$ sufficiently $C^1$-close to the identity so that the
Hofer norm $\|H^\phi\|$ satisfies
$$
\|H^\phi\| < \frac{1}{2} A(\omega; J_0) \tag 9.26
$$
and so $\gamma(\phi) \leq \frac{1}{2}A(\omega; J_0)$. By the
definition of $H^\phi$ in section 7, $d_{C^1}(\phi_{H^\phi},id)$
can be made arbitrarily small and so satisfies the hypotheses of
Proposition 9.1. Then applying Proposition 9.1 we derive from (9.26)
that any solution $v$ of (1.12) satisfying
$$
E_{J'}(v) < \frac{1}{2} A(\omega;J_0) \tag 9.27
$$
must indeed have the energy which is precisely the same as
$$
E_{J'}(v) = \|H^\phi\| = osc(S_\phi).
$$
Therefore to prove Theorem 7.1 it is enough to prove that
(1.14) has a solution $v$ satisfying (9.27).

Now we run the existence scheme used in the proof of Theorem 5.4
in section 8 for the cycles $\alpha_{H^\phi}$ and
$\beta_{\widetilde H^{\phi}}$ above. Without loss of any
generality, we may assume that $H^\phi$ satisfies the hypothesis
in Lemma 9.4 if we choose $\phi$ is sufficiently $C^1$-small.

In particular, the canonical fundamental cycle $\alpha_{H^\phi}$
has the form
$$
\alpha_{H^\phi} = [x^-, \widehat x^-] + \beta \tag 9.28
$$
with $\lambda_{H^\phi}(\beta) < \AA_{H^\phi}([x^-,\widehat x^-])$
and is tight in the local Floer complex. Similarly we have
$$
\beta_{\widetilde H^\phi} = [x^+,\widehat x^+] + \widetilde \beta
\tag 9.29
$$
with $\lambda_{\widetilde H^\phi}(\widetilde \beta) <
\AA_{\widetilde H^\phi}([x^+,\widehat x^+])$ and tight in the
corresponding local Floer complex. In practice, to run our existence
scheme as in section 8, we need to consider the admissible triple
$$
H_1 = H^\phi, \, H_2 = \widetilde{H^\phi}\#(\e f), \, H_3 = \e
f
$$
for a fixed Morse function $f$, and the corresponding cycles
$\alpha_{H^\phi}$ and $\beta_{\widetilde{H^\phi}\#(\e f)}$. It follows
that by choosing $\e$ sufficiently small, $\beta_{\widetilde{H^\phi}\#(\e f)}$
is also tight in the local Floer complex
and have the similar structure as in (9.29).

Since both $\alpha_H$ and $\beta_{\widetilde{H^\phi}\#(\e f)}$
are tight  in the local Floer complex, the pants product
$$
\alpha_H * \beta_{\widetilde{H^\phi}\#(\e f)}
$$
has a {\it nontrivial} contribution coming from the thin moduli
space
$$
\MM_3^0(H,\widetilde J;\widehat p)
$$
defined as in subsection 9.2, for the asymptotic orbits $\widehat
p$ consisting of
$$
\widehat p = ([x^-,\widehat x^-], [x^+,\widehat x^+]; [q, \widehat
q])
$$
with $q \in \text{Crit}(-\e f)$. Since this is the case for all
small $\e >0$, we can now carry out the same proof as in section
8. {\it The upshot of the above discussion for
$H^\phi$ is that case (2) of Theorem 5.4 cannot occur for the
above choices of $\phi$ and of the canonical Floer fundamental
cycles $\alpha_H$ and $\beta_{\widetilde{H^\phi}\#(\e f)}$.}
This is thanks to the inequality (5.18) imposed for such a solution
in Theorem 5.4 : In the current case, for given $J_0$,
we choose $\phi$ sufficiently $C^1$ close to the identity.
Such a choice makes $H^\phi$ sufficiently $C^0$ small so that
$\|H^\phi\| < \frac{1}{3}A(\omega;J_0)$. Then (3.2), applied to
$b=1$, implies that
$$
\rho(H^\phi;1) + \frac{\delta}{2}, \,\,
\rho(\widetilde H^\phi;1) + \frac{\delta}{2} < \frac{2}{3}A(\omega;J_0)
$$
if we choose $\delta < \frac{1}{3}A(\omega;J_0)$.
Now we choose $J'\in j_{(\phi,J_0)}$ sufficiently $C^1$-close to
the constant path $J_0$ which gives rise to
$$
A(\omega;J'):= \inf_{t \in [0,1]}\{A(\omega;J'_t)\} \geq
\frac{2}{3} A(\omega;J_0),
$$
by the lower semicontinuity property of the function
$J' \mapsto A(\omega;J')$: this lower semicontinuity
can be easily proven by a similar argument as the proof
of Proposition 4.4. This precludes the possibility (2) in Theorem 5.4.

Once this is the case, a simple examination of the existence proof
of Theorem 5.4 in section 8 immediately gives rise to the
following existence theorem of thin (broken-trajectory) solutions.

\proclaim{Proposition 9.6} Let $J_0$ be any compatible almost
complex structure. Assume that $\phi$ is sufficiently $C^1$-small
and that the graph, $\text{graph }\phi$, is contained in the
Darboux neighborhood $\UU$ as in (7.2). Let $H^\phi: [0,1] \times
M \to \R$ be the quasi-autonomous Hamiltonian that generates the
Hamiltonian path $\phi^t$ determined by (7.2). Let $\alpha_H \in
CF_n(H^\phi), \, \beta_{\widetilde H^\phi}\in CF_n(\widetilde
H^\phi)$ be as before.  Then there must be a `thin'
broken-trajectory solution
$$
u = u_1\# u_2 \# \cdots \# u_N
$$
of
$$
\cases
\dudtau + J_t \Big(\dudt - X_H(u)\Big) = 0 \\
u_1(-\infty) = [x^-,\widehat x^-] \in \alpha_H, \,
u_N(\infty) = [x^+,\widehat x^+],  \hbox{with }\,
[\widetilde x^+, \widetilde{\widehat x^+}] \in {\beta_{\widetilde H}}, \\
u_j (0,0) = q \in B(u) \quad \text{for some $j=1, \cdots, N$ }
\endcases
\tag 9.30
$$
with $J_t' = (\phi_H^t)^*J_0$, that satisfy
$$
\aligned \AA_H(u(-\infty)) & =\int_0^1 -\min H_t \, dt = \int_0^1
- H_t(x^-) \, dt
= \max S_\phi\\
\AA_H(u(\infty)) & = \int_0^1 -\max H_t \, dt = \int_0^1 -H_t(x^+)
\,dt =  \min S_\phi.
\endaligned
\tag 9.31
$$
\endproclaim

Finally we are ready to finish the proof of Theorem 7.1.

\demo{Wrap-up of the proof of Theorem 7.1} We choose $\alpha_i <
\frac{1}{3} A(\omega;J_0)$ in Proposition 9.1. We will prove that
for the choice of $J'$ given by $J_t' = (\phi_{H^\phi}^t)^*J_0$,
we have
$$
A(\phi,J_0;J';1) \geq osc(S_\phi).
$$
This will then finish the proof of Theorem 7.1.

Proposition 9.6 implies the existence of a {\it thin}
broken-trajectory solution that satisfies
$$
u(-\infty) = [x^-,\widehat x^-], \quad u(\infty) = [x^+,\widehat
x^+]
$$
and hence
$$
\AA_H(u(-\infty)) - \AA_H(u(\infty)) = \hbox{\rm osc}(S_\phi).
$$
Furthermore it follows from Proposition 9.1 that we have
$$
\sup_{\tau \in \R}\hbox{ diam }(t \mapsto u_j(\tau,t)) < \alpha_3
$$
for all $j =1, \cdots, N$. We note that the asymptotic orbits of
$u$ are constant loops and so each $u_j$ defines a cycle in $M$. In
addition, Lemma 8.6 implies that the image of $u$ is {\it
homologous to a one dimensional path} that is the join
$$
\chi_1 \# \cdots\# \chi_N
$$
where $\chi_j: \R \to M$ is the `center of mass' path of the loops
$$
t \mapsto u_j(\tau,t), \quad j = 1, \cdots, N
$$
defined for each $\tau \in \R$ as in Lemma 8.6. This latter
topological property holds for any thin cusp-solutions $u'$. Now
the same calculation as the one carried out in the proof of Lemma 4.6 gives rise to
$$
E_J(u') = \int (u')^*\omega + \int_0^1(-H_t(x^-) + H_t(x^+)) \,
dt.
$$
Since $u'$ is homologous to a one dimensional path and so
homologous to zero as a two chain, $\int (u')^*\omega = 0$. Therefore we
obtain
$$
E_J(u') = \int_0^1( -H_t(x^-) + H_t(x^+)) \, dt = \hbox{\rm osc}(S_\phi)
$$
where the second identity follows from  (7.4). This and (9.27)
imply
$$
\hbox{\rm osc}(S_\phi) =  E_{J}(u') < \frac{1}{2}A(\omega;J_0),
\tag 9.32
$$
 for any thin broken-trajectory solution $u'$
connecting $[x^-,\widehat x^-]$ and $[x^+,\widehat x^+]$. On the
other hand, for any {\it thick} solution $u''$, we have
$$
E_J(u'') \geq A(\omega; J_0) - \alpha_2 >
\frac{2}{3}A(\omega;J_0). \tag 9.33
$$
Combining (9.32) and (9.33), we have indeed proved
$$
A(\phi,J_0;J';1) = osc(S_\phi)
$$
for the choice $J'= \{(\phi_{H^\phi}^t)^*J_0\}$. By the definition
of $A(\phi;1)$, Definition 6.2,  we have proved $A(\phi;1) \geq
\hbox{\rm osc}(S_\phi)$. This finishes the proof of Theorem 7.1.
\qed\enddemo

\head{\bf \S 10. Thick and thin decomposition of the Floer moduli
space}\endhead

In this section, we will prove Proposition 9.1 and 9.3.
Similar theorem can be proven for the higher Massey product
moduli spaces with arbitrary $k$ number of punctures on $\Sigma$.
The cases of our main interest correspond to $k = 2, \, 3$. We start
with the case of $k=2$. We recall a similar result from
[Proposition 10.1, Oh2] in the context of Lagrangian intersection
Floer theory for $k=2$.

\demo{Proof of Proposition 9.1} The proof is by contradiction as
in [Oh2,5]. Suppose to the contrary that for some $0 < \alpha_1,
\, \alpha_2 < A(\omega; J_0)$ there exists a sequence $\delta_j\to
0$, $H_j$ with $\|H_j\| \leq \delta_j$ and $u_j$ satisfying (9.2)
for $H_j$ but with the bound
$$
\alpha_1 \leq \int\Big|{\part u_j \over \part \tau}
\Big|^2_{J_0} = \int\Big|{\part v_j \over \part \tau}
\Big|^2_{J'_t} \leq A(\omega; J_0) - \alpha_2 \tag 10.1
$$
Because of the energy upper bound in (10.1), we can apply Gromov's
type of compactness argument. Since we assume that any
one-periodic trajectory is constant, the homotopy class $[u] \in
\pi_2(M)$ is canonically defined for any finite energy solution
$u$ and in particular for $u_j$. A straightforward computation
shows
$$
\int \Big|{\part u_j \over \part\tau}\Big|^2_{J_0} = \int
\omega\Big({\part u_j \over \part \tau}, {\part u_j \over \part
t}\Big)\, d\tau\, dt - \int_0^1 (H(u_j(\infty)) -
H(u_j(-\infty)))\, dt \tag 10.2
$$
and so for the finite energy solution $u_j$, we derive
$$
\int \Big|{\part u_j \over \part \tau}\Big|^2_{J_0} =
\omega([u_j]) - \int_0^1 (H(u_j(\infty)) - H(u_j(-\infty))\, dt
\leq \omega([u_j]) + \|H\|. \tag 10.3
$$
Due to the energy bound (10.1) and by a standard compactness
argument, we may assume that the homotopy class $[u_j] =A$ fixed
by choosing a subsequence if necessary. If $A= 0$,
then (10.3) implies
$$
\int \Big|{\part u_j \over \part \tau}\Big|^2_{J_0} \leq \|H\|
\leq \delta_j
$$
which contradicts the lower bound in (10.1) if $j$ is sufficiently
large so that $\delta_j < \alpha_1$. Now assume that $A \neq 0$.
In this case, there must be some $C >0$ and $\tau_j \in \R$ with
$$
\text{diam} (t\mapsto u_j(t,\tau_j)) \geq C > 0 \tag 10.4
$$
where $C$ is uniform over $j$. By translating $u_j$ along the
direction of $\tau$, we may assume that $\tau_j = 0$. Now if
bubbling occurs, we just take the bubble to produce a non-constant
$J_0$-holomorphic sphere.  If bubbling does not occur, we take a
local limit around $\tau = 0$ using the energy bound (10.1), to
produce a map $u: \R \times S^1 \to M$ satisfying
$$
\dudtau + J_0 \dudt = 0. \tag 10.5
$$
Furthermore if there occurs no bubbling and so if the sequence
converges uniformly to the local limit around $\tau=0$, then the
local limit cannot be constant because of (10.4) and so it must be
a non-constant $J_0$-holomorphic cylinder with the energy bound
$$
{1 \over 2} \int |Du|_{J_0}^2 \leq A(\omega; J_0) - \alpha_2 <
\infty.
$$
By the removable singularity theorem, $u$ can be compactified to a
$J_0$-holomorphic sphere that is non-constant.

Therefore whether the sequence bubbles off or not, we have
$$
\limsup_j \int \Big|{\part u_j \over \part \tau}\Big|_{J_0}^2 \geq
A(\omega; J_0).
$$
This then contradicts the upper bound of (10.1). This finishes the
proof of the alternative (9.1)-(9.2). Again by arguing by
contradiction, we derive (9.3) from (10.4) and (9.4) from (10.3).
This finishes the proof. \qed\enddemo

Now we consider the case with $k=3$ and give the proof of
Proposition 9.3. The main idea of the proof for this case is
essentially the same as that of $k=2$, except that we need set-up
the appropriate Hamiltonian fibration in the proof. In the proof,
we are content to give the proof of Proposition 9.3 and refrain
from trying to provide the optimal form of the exposition on the
various computations related to the curvature and the Hofer type
quantities. We will postpone an optimal form of the exposition to
[Oh10].

\demo{Proof of Proposition 9.3}  For given $H=(H_1,H_2,H_3)$ with
$H_1 \# H_2 = H_3$, we choose a Hamiltonian fibration $P \to
\Sigma$ with connection $\nabla$ whose monodromy at the ends given
by the $H_i$ for $i = 1,2, \, 3$ as constructed in section 8. (See
[Oh9] for more details.)

In the minimal metric representation as in section 8, we have the
formula of the curvature two form
$$
K(s(\tau,t)) d\tau \wedge dt
$$
of $\nabla$ where $K$ is the function
$$
K(s(\tau,t)) = \sum_{i=1}^2 \rho_-'(\tau) H_i(t, u(\tau,t)) -
\rho_+'(\tau)(H_1\# H_3)(t, u(\tau,t)) \tag 10.6
$$
on each $\Sigma_i$ where $\Phi\circ s(\tau,t) =
(\tau,t,u(\tau,t))$ in the trivialization $\Phi$. (See [section 3,
Oh9] for the derivation of this formula, more specifically the
formula (3.16) of [Oh9].) We note that the curvature is
concentrated near the transition regions of the cut-off functions
and in particular compactly supported. Furthermore we have the
integral bound for the curvature integral
$$
-E^+(H) \leq \int K\circ s \leq E^-(H). \tag 10.7
$$
and the $L^\infty$-bound for $K$
$$
-E_\infty^+(H) \leq K(s(z)) \leq E^-_\infty(H) \tag 10.8
$$
where
$$
\align E^-_\infty(H) &:= \max\{-\min_{t,x} H_1(t,x),
\, -\min_{t,x}H_2(t,x); \, \max_{t,x}H_3(t,x) \} \\
E^+_\infty(H) &:= \max\{\max_{t,x} H_1(t,x),
\, \max_{t,x}H_2(t,x); \, -\min_{t,x}H_3(t,x) \}\\
E_\infty(H) &:= E^-_\infty(H) + E^+_\infty(H).
\endalign
$$
We note that for the triple $H = (H_1, H_2; H_1 \# H_2)$ it was
shown in [GLS], [En1], [Oh9] that the coupling two form $\omega_P$
of the connection $\nabla$ has the property that
$$
\Omega_{P,\lambda} = \omega_P + \lambda\omega_\Sigma
$$
is non-degenerate for all $\lambda \in [\delta, \infty)$, provided
$E_\infty(H)< \delta$. Let $\widetilde J$ be an $(H,J)$-compatible
almost complex structure which is also compatible to the
symplectic form $\Omega_{P,\lambda}$. For the simplicity of
notation, we denote $| \cdot | = |\cdot |_{\widetilde J}$. For a
$\widetilde J$-holomorphic section, we have the identity
$$
{1 \over 2} \int |Ds|^2 = \int s^*(\Omega_{P,\lambda}) = \int
s^*\omega_P + \lambda \tag 10.9
$$
We decompose $Ds = (Ds)^v + (Ds)^h$ into vertical and horizontal
parts and write
$$
|Ds|^2 = |(Ds)^v|^2 + |(Ds)^h|^2 + 2 \langle (Ds)^v, (Ds)^h
\rangle.
$$
Now it is straightforward to prove
$$
|(Ds)^h|^2d\tau\wedge dt = 2(K(s)d\tau \wedge dt +
\lambda\omega_\Sigma) \tag 10.10
$$
by evaluating
$$
\align \sum_{i=1}^2|(Ds)^h(e_i)|^2 & =
\sum_{i=1}^2\Omega_{P,\lambda}((Ds)^h(e_i),
\widetilde J (Ds)^h(e_i)) \\
& = \sum_{i=1}^2 (\omega_P + \lambda \omega_{\Sigma})((Ds)^h(e_i),
\widetilde J (Ds)^h(e_i))
\endalign
$$
for an orthonormal frame $\{e_1, e_2\}=\Big\{\frac{\part}{\part
\tau}, \frac{\part}{\part t}\Big\}$ of $T\Sigma$. We refer to
[Oh9] for detailed derivation of these formulae.

Combining (10.9)-(10.10), we have derived
$$
\aligned
\frac{1}{ 2} \int |(Ds)^v|^2 & = \int s^*\omega_P + \lambda -
\Big(\int K(s)d\tau\wedge dt + \int \lambda \omega_\Sigma\Big) \\
& = \int s^*\omega_P - \int K(s)
\endaligned
\tag 10.11
$$
which is equivalent to the formula [(4.13), Oh9]. We now divide
our discussion into two.

Firstly, if $\int s^*\omega_P = 0$, we derive
$$
\frac{1}{2} \int |(Ds)^v|^2 \leq E(H) \tag 10.12
$$
from (10.6) and (10.11). Therefore
$$
\frac{1}{2}\int |(Ds)^v|^2 \leq \delta. \tag 10.13
$$
Combining (10.11)-(10.13), we have established
$$
\int |Ds|^2 \leq \alpha_1
$$
if $\delta> 0$ is chosen sufficiently small.

Secondly, if $\int s^*\omega_P \neq 0$, then $pr_\Sigma\circ s:
\Sigma \to \Sigma$ has degree one and the fiber homotopy class
$[u]$ of $s$ satisfies
$$
[u] =: A \neq 0 \in \pi_2(M). \tag 10.14
$$
Furthermore noting that as $\|H\|_{C^1} \to 0$, the connection can
be made closer and closer to the trivial connection in the trivial
fibration $P = \Sigma \times M$ and the $(H,J)$-compatible $\widetilde
J$ also converges to the product almost complex
structure $j \oplus J_0$ and hence the image of any $\widetilde
J$-holomorphic section cannot be completely contained in the
neighborhood of one of the obvious horizontal sections
$$
\Sigma \times \{q\}
$$
for any one fixed $q \in M$. Now consider a sequence $H_j$ with
$\|H_j\|_{C^1} \to 0$, and $(H_j,J)$-compatible almost complex
structure $\widetilde J_j$, and let $s_j$ be a sequence of
$\widetilde J_j$-holomorphic sections in the fixed fiber class
(10.14). In other words, if we write
$$
\Phi\circ s_j(z) = (z, u_j(z))
$$
in the trivialization $\Phi:P \to \Sigma \times M$, then we have
$[u_j] = A \neq 0$. Since we assume $A \neq 0$, there is a
constant $C > 0$ such that
$$
\text{diam}(u_j) \geq C > 0
$$
for all sufficiently large $j$ after choosing a subsequence. By
applying a suitable conformal transformation on the domain, either
by taking a bubble if bubble occurs or by choosing a limit when
bubbling does not occur, we can produce at least one non-constant
$J_0$-holomorphic map
$$
u_\infty: S^2 \to M
$$
out of the $u_j$'s as in the case $k=2$ before. Furthermore we
also have the energy bound
$$
\limsup_{j \to \infty} \int |(Ds_j)^v|^2_{\widetilde J_j} \geq
\int |Du_\infty|^2_{J_0}.
$$
Therefore we have
$$
\int |Ds_j|^2 \geq \int |(Ds_j)^v|^2 \geq \int |Du_\infty|^2_{J_0}
\geq A(\omega; J_0) - \alpha_2
$$
for a sufficiently large $j$. This finishes the proof of
Proposition 9.1. \qed\enddemo

\definition{Remark 10.1} In the above proof of Proposition 9.3,
the readers might be wondering why we are short of stating
\medskip

``By the Gromov compactness theorem, the sequence $s_j$ converges
to
$$
s_h + \sum_{k=1}^n w_k, \quad n\neq 0 \tag 10.15
$$
as $j \to \infty$ or $\|H_j\| \to 0$, where $s_h$ is an obvious
horizontal section and each $w_j$ is a $J_0$-holomorphic sphere
into a fiber $(M,\omega)$ of $P = \Sigma \times M$.''
\medskip

The reason is such a convergence result {\it fails} in general for
two reasons: First unless we specify how the limiting sequence
$H_j$ converges to 0, the sequence $s_j$ cannot have any limit in
any reasonable topology. This is because the case $H=0$ is a
singular situation in the study of the Floer moduli space
$\MM(H,\widetilde J;\widehat z)$. Secondly even if we specify a
good sequence, e.g., consider the `adiabatic' sequence
$$
H_{1,j} = \e_j f_1, \, H_{2,j} = \e_j f_2, \, H_{3,j} = \e_j f_3
$$
for Morse functions $f_1,\, f_2, \, f_3$ with the same sequence
$\e_j \to 0$, we still have to deal with the degenerate limit,
i.e. the limit that contains components of Hausdorff dimension one
as studied in section 8. What we proved in the proof is that we
can always produce at least one non-constant $J_0$-holomorphic
sphere as $j \to \infty$ without using such a strong convergence
result, when the homotopy class of the $s_j$ is not trivial (in
the fiber direction).
\enddefinition

\definition{Remark 10.2} One can easily see that there are many
rooms for improvements in the above proof and in the statements of
the results in section 9. We summarize several points of
improvement: \roster

\item The statement of Proposition 9.4 can be further improved
using more systematic usage of  Hamiltonian fibrations and
associated Hamiltonian connections in the spirit of [En1] and
[Oh9]. (See the proof of Proposition 9.4 above for some flavor
thereof.)

\item All the arguments used in section 9 and 10 can be
generalized to the case where $\phi$ is engulfable. This requires
a more systematic discussion on the local Floer complexes and more
precise estimates involving the curvature of the Hamiltonian
fibration.

\item The inequalities (9.26), (9.27) are obviously not optimal.
\endroster
\enddefinition

In the sequel [Oh10], we will further elaborate the above points
and derive various consequences of them in relation to Hofer's
geometry of the Hamiltonian diffeomorphism group.

\head{\bf Appendix : Construction of flat connections}
\endhead
\smallskip

In this appendix, we will construct a flat Hamiltonian connection
associated to the triple (8.11)
$$
H_1 = H, \, H_2 = \widetilde{H}\#(\e f), \, H_3 =  \e f.
$$
We refer to [Oh9] for a detailed calculation of the curvature
associated to a {\it deformed mapping cylinder}. In our case,
such a flat connection in our minimal area metric representation
of the fibration $P \to \Sigma$ can be constructed gluing three
two-parameter family
$$
\phi_i: (s,t) \in [0,1] \times [0,1] \to Ham(M,\omega)
$$
of Hamiltonian diffeomorphisms and then elongating the family
in the $s$-direction via a
cut-off function $\rho$ of the type (8.2). We now give the
explicit formulae of them that glue smoothly along the
boundary of $\Sigma_i$ according to the gluing rule (8.1).
This flat connection has singularities at the two points $p, \,
\overline p$ mentioned before. Without loss of any generality,
we assume that $H$ is boundary
flat, i.e., that there exists $\delta_0 > 0$ such that
$$
H_t\equiv 0 ; \quad t \in [0,\delta_0] \cup [1-\delta_0,1].
\tag A.1
$$

Here come the formulae for $\phi_i$s :
$$
\align
\phi_1(s,t) & = \cases id & \quad \hskip0.6in
 0 \leq t \leq \frac{1-s}{2} \\
\phi_H^{\frac{t-\frac{1-s}{2}}{\frac{1+s}{2}}}
& \quad \hskip0.6in \frac{1-s}{2} \leq t \leq 1
\endcases \\
\phi_2(s,t) & = \cases \phi_H^{\frac{\frac{1+s}{2}-t}{\frac{1+s}{2}}}
& \quad  0 \leq t \leq \frac{1-s}{2}\\
\phi_H^{\frac{\frac{1+s}{2}-t}{\frac{1+s}{2}}} \circ
\phi_{\e f}^{\frac{t-\frac{1-s}{2}}{\frac{1+s}{2}}}
& \quad \frac{1-s}{2} \leq t \leq 1 \\
\endcases\\
\phi_3(s,t) & = \cases \phi_{\e f}^{\frac{t}{\frac{1+s}{2}}}
& \quad \hskip0.7in 0 \leq t \leq \frac{1+s}{2}\\
id &\quad
\hskip0.7in \frac{1+s}{2}\leq t \leq 1.
\endcases
\endalign
$$
Compatibility to the gluing rule (8.29) can be achieved by
rescaling the domain of $\phi_i$ in the $s$-direction
in the obvious way, and smoothness of the family away from
$p, \, \overline p$ can be checked easily. Flatness of the
associated connection
obviously follows from the fact that it comes from the above
smooth two parameter family of Hamiltonian diffeomorphisms.
(See [Ba] for the relevant argument.)

We like to note that due to the boundary flatness assumption
(A.1) on $H$ the associated connections converges to the trivial
connection on $D^2(\delta_0) \setminus \{(0,0)\}$,
or equivalently the associated Hamiltonians converges to
the zero Hamiltonian on $D^2(\delta_0) \setminus \{(0,0)\}$
with respect to the standard metric of $S = \R \times S^1$.

Now we elongate the pants in the direction of $s$ using cut-off
functions $\rho$ as defined in (8.2), which gives rise to a flat
connection on $\Sigma$ that is smooth everywhere except two points
$p, \, \overline p$.  One important point
to note here is that as $\e \to 0$, the above flat connection on
$\Sigma_1 \cup \Sigma_2$ converges to the connection on the mapping
cylinder induced by the Hamiltonian $H$ as described in subsection
8.2.

\head {\bf References}
\endhead

\enddocument